\documentclass[leqno,english]{amsart}
\usepackage{hyperref}
\usepackage{amssymb}
\usepackage[square,numbers]{natbib}
\usepackage[all]{xy}

\numberwithin{subsubsection}{section}

\renewcommand{\theenumi}{\arabic{enumi}}
\renewcommand{\labelenumi}{(\theenumi)}

\theoremstyle{plain}

\newtheorem*{mainthm}{Theorem}

\swapnumbers
\newtheorem{thm}[subsubsection]{Theorem}
\newtheorem{lemm}[subsubsection]{Lemma}
\newtheorem{prop}[subsubsection]{Proposition}
\newtheorem*{claim}{Claim}
\newtheorem{fact}[subsubsection]{Fact}
\newtheorem{obsv}[subsubsection]{Observation}

\DeclareMathOperator{\Map}{Map}

\DeclareMathOperator{\Mor}{Mor}

\DeclareMathOperator{\Hom}{\mathit{Hom}}

\DeclareMathOperator{\Op}{\mathcal{O}}
\DeclareMathOperator{\Bij}{\mathit{Bij}}

\DeclareMathOperator{\unit}{\mathbf{1}}

\DeclareMathOperator{\Sh}{shuffle}

\DeclareMathOperator{\NN}{\mathbb{N}}
\DeclareMathOperator{\ZZ}{\mathbb{Z}}
\DeclareMathOperator{\kk}{\Bbbk}
\DeclareMathOperator{\FF}{\mathbb{F}}
\DeclareMathOperator{\QQ}{\mathbb{Q}}

\DeclareMathOperator{\AOp}{\mathtt{A}}

\DeclareMathOperator{\COp}{\mathtt{C}}
\DeclareMathOperator{\DOp}{\mathtt{D}}
\DeclareMathOperator{\EOp}{\mathtt{E}}
\DeclareMathOperator{\FOp}{\mathtt{F}}
\DeclareMathOperator{\GOp}{\mathtt{G}}
\DeclareMathOperator{\KOp}{\mathtt{K}}
\DeclareMathOperator{\LOp}{\mathtt{L}}
\DeclareMathOperator{\POp}{\mathtt{P}}
\DeclareMathOperator{\QOp}{\mathtt{Q}}
\DeclareMathOperator{\ROp}{\mathtt{R}}
\DeclareMathOperator{\SOp}{\mathtt{S}}
\DeclareMathOperator{\Free}{\mathtt{F}}

\DeclareMathOperator{\Sym}{S}
\DeclareMathOperator{\Tens}{T}
\DeclareMathOperator*{\TFact}{T^c\Sigma}
\DeclareMathOperator*{\BFact}{\mathit{B}}

\DeclareMathOperator{\C}{\mathcal{C}}

\DeclareMathOperator{\E}{\mathcal{E}}
\DeclareMathOperator{\K}{\mathcal{K}}
\DeclareMathOperator{\Simp}{\mathcal{S}}
\DeclareMathOperator{\M}{\mathcal{M}}

\DeclareMathOperator{\eset}{\underline{\mathsf{e}}}
\DeclareMathOperator{\fset}{\underline{\mathsf{f}}}
\DeclareMathOperator{\gset}{\underline{\mathsf{g}}}
\DeclareMathOperator{\hset}{\underline{\mathsf{h}}}
\DeclareMathOperator{\rset}{\underline{\mathsf{r}}}

\DeclareMathOperator{\uset}{\underline{\mathsf{u}}}
\DeclareMathOperator{\vset}{\underline{\mathsf{v}}}

\DeclareMathOperator{\Tr}{Tr}

\DeclareMathOperator{\Id}{Id}
\DeclareMathOperator{\id}{id}
\DeclareMathOperator{\sgn}{sgn}
\DeclareMathOperator{\coker}{coker}
\DeclareMathOperator{\im}{im}
\DeclareMathOperator{\colim}{colim}
\DeclareMathOperator{\cst}{\mathit{cst}}
\DeclareMathOperator{\Ab}{Ab}

\DeclareMathOperator{\Tor}{Tor}
\DeclareMathOperator{\Indec}{Indec}

\DeclareMathOperator{\gr}{\mathit{gr}}

\DeclareMathOperator{\Ho}{Ho}

\DeclareMathOperator{\HGamma}{\mathit{H\Gamma}}
\DeclareMathOperator{\CXi}{\mathit{C\Xi}}
\DeclareMathOperator{\NXi}{\mathit{N\Xi}}

\DeclareMathOperator{\dercirc}{\triangleright}

\title[Iterated bar complexes of E-infinity algebras]{Iterated bar complexes of E-infinity algebras\\and\\homology theories}
\author{Benoit Fresse}
\date{31 March 2010 (current version -- preliminary version in October 2008)}
\address{UMR 8524 de l'Universit\'e Lille 1 - Sciences et Technologies - et du CNRS\\
Cit\'e Scientifique -- B\^atiment M2\\
F-59655 Villeneuve d'Ascq C\'edex (France)}
\email{Benoit.Fresse@math.univ-lille1.fr}
\urladdr{http://math.univ-lille1.fr/\~{ }fresse}
\thanks{Research supported in part by ANR grant JCJC06-0042 OBTH}

\begin{document}

\begin{abstract}
We proved in a previous article that the bar complex of an $E_\infty$-algebra
inherits a natural $E_\infty$-algebra structure.
As a consequence, a well-defined iterated bar construction $B^n(A)$
can be associated to any algebra over an $E_\infty$-operad.
In the case of a commutative algebra $A$,
our iterated bar construction reduces to the standard iterated bar complex of $A$.

The first purpose of this paper is to give a direct effective definition of the iterated bar complexes of $E_\infty$-algebras.
We use this effective definition to prove that the $n$-fold bar construction admits an extension
to categories of algebras over $E_n$-operads.

Then we prove that the $n$-fold bar complex determines the homology theory
associated to the category of algebras over an $E_n$-operad.
In the case $n = \infty$,
we obtain an isomorphism between the homology of an infinite bar construction
and the usual $\Gamma$-homology with trivial coefficients.
\end{abstract}

\maketitle

\section*{Contents}{\parindent=0cm

\medskip\textbf{Introduction, \pageref{Intro}}

\medskip\textbf{Conventions and background, \pageref{Background}}

\medskip\textbf{The construction of iterated bar complexes, \pageref{IteratedComplexDefinition}}

{\leftskip=1cm

\textit{\S\ref{IteratedBarModules}. Iterated bar modules,~\pageref{IteratedBarModules}}

\textit{\S\ref{QuasiFreeLifting}. Iterated bar modules as quasi-free modules,~\pageref{QuasiFreeLifting}}

\textit{\S\ref{CompleteGraphOperad}. Interlude: operads shaped on complete graph posets,~\pageref{CompleteGraphOperad}}

\textit{\S\ref{KStructure}. The structure of iterated bar modules,~\pageref{KStructure}}

\textit{\S\ref{EnDefinition}. The restriction of iterated bar modules to $E_n$-operads,~\pageref{EnDefinition}}

}

\medskip\textbf{Iterated bar complexes and homology theories,~\pageref{HomologyTheories}}

{\leftskip=1cm

\textit{\S\ref{UsualBarComplexes}. Prelude: iterated bar complexes of usual commutative algebras,~\pageref{UsualBarComplexes}}

\textit{\S\ref{OperadicHomology}. Homology of algebras over operads and operadic $\Tor$-functors,~\pageref{OperadicHomology}}

\textit{\S\ref{EnHomology}. Iterated bar complexes and homology of algebras over $E_n$-operads,~\pageref{EnHomology}}

\textit{\S\ref{EinfinityHomology}. Infinite bar complexes,~\pageref{EinfinityHomology}}

}

\medskip\textbf{Afterword: applications to the cohomology of iterated loop spaces,~\pageref{Conclusion}}

\medskip\textbf{Bibliography,~\pageref{Biblio}}

}

\part*{Introduction}\label{Intro}

The standard reduced bar complex $B(A)$
is defined basically as a functor from the category of associative differential graded algebras
(associative dg-algebras for short)
to the category of differential graded modules
(dg-modules for short).
In the case of a commutative dg-algebra,
the bar complex $B(A)$ inherits a natural multiplicative structure
and still forms a commutative dg-algebra.
This observation implies that an iterated bar complex~$B^n(A)$
is naturally associated to any commutative dg-algebra $A$, for every $n\in\NN$.
In this paper,
we use techniques of modules over operads to study extensions of iterated bar complexes
to algebras over $E_n$-operads.
Our main result asserts that the $n$-fold bar complex~$B^n(A)$
determines the homology theory associated to an $E_n$-operad.

\medskip
For the purpose of this work,
we take the category of dg-modules
as a base category
and we assume that all operads belong to this category.
An $E_n$-operad refers to a dg-operad equivalent to the chain operad of little $n$-cubes.
Many models of $E_n$-operads
come in nested sequences
\begin{equation}
\EOp_1\subset\dots\subset\EOp_n\subset\dots\subset\colim_n\EOp_n = \EOp
\end{equation}
such that $\EOp = \colim_n\EOp_n$
is an $E_\infty$-operad,
an operad equivalent in the homotopy category of dg-operads
to the operad of commutative algebras $\COp$.
Recall that an $E_1$-operad is equivalent to the operad of associative algebras $\AOp$
and forms, in another usual terminology, an $A_\infty$-operad.
The structure of an algebra over an $A_\infty$-operad
includes a product $\mu: A\otimes A\rightarrow A$
and a full set of homotopies
that make this product associative.
The structure of an algebra over an $E_\infty$-operad
includes a product $\mu: A\otimes A\rightarrow A$
and a full set of homotopies
that make this product associative and commutative.
The intermediate structure of an algebra over an $E_n$-operad
includes a product $\mu: A\otimes A\rightarrow A$
which is fully homotopy associative,
but homotopy commutative up to some degree only.

\medskip
Throughout the paper,
we use the letter $\C$ to denote the category of dg-modules
and the notation ${}_{\POp}\C$,
where $\POp$ is any operad,
refers to the category of $\POp$-algebras in dg-modules.
The category of commutative dg-algebras,
identified with the category of algebras over the commutative operad $\COp$,
is denoted by ${}_{\COp}\C$.

Recall that an operad morphism $\phi: \POp\rightarrow\QOp$
gives rise to a restriction functor $\phi^*: {}_{\QOp}\C\rightarrow{}_{\POp}\C$
since the restriction of a $\QOp$-action through~$\phi$ provides any $\QOp$-algebra
with a natural $\POp$-action.
The category of $\POp$-algebras is also equipped
with an obvious forgetful functor to the base category $\C$.
For a nested sequence of $E_n$-operads,
we have a chain of restriction functors
\begin{equation*}
\xymatrix{ {}_{\EOp_1}\C &
\ar[l] & \ar@{.}[l] & {}_{\EOp_n}\C\ar[l] &
\ar[l] & \ar@{.}[l] & \ar[l]{}_{\EOp}\C & {}_{\COp}\C\ar[l] }
\end{equation*}
starting from the category of commutative algebras.

The bar construction is defined basically as a functor from the category of associative dg-algebras
to the category of dg-modules.
This usual bar construction has an extension to any category of algebras over an $A_\infty$-operad,
equivalently to any category of algebras over an $E_1$-operad.
For a nested sequence of $E_n$-operads,
the restriction of the usual bar construction to commutative dg-algebras
and the extended bar complex $B: {}_{\EOp_1}\C\rightarrow\C$
fit a commutative diagram
\begin{equation*}
\xymatrix{ {}_{\EOp_1}\C\ar@{.>}[drr]_{B} & \ar[l] & \ar@{.}[l] &
{}_{\EOp}\C\ar[l] & {}_{\COp}\C\ar[l]\ar[dll]^{B} \\
&& \C && }.
\end{equation*}

\medskip
Recall that the bar construction
of commutative dg-algebras
factors through the category of commutative dg-algebras:
\begin{equation*}
\xymatrix{ {}_{\COp}\C\ar[dr]\ar@{.>}[rr]^{B} && {}_{\COp}\C\ar[dl]^{\text{forgetful}} \\
& \C & }.
\end{equation*}
Suppose that the $E_\infty$-operad $\EOp$ is cofibrant as a dg-operad.
In~\cite{Bar1}
we prove that the bar construction of $\EOp$-algebras
admits a factorization through the category of $\EOp$-algebras
so that the functor $B: {}_{\COp}\C\rightarrow{}_{\COp}\C$
admits an extension
\begin{equation*}
\xymatrix{ {}_{\EOp}\C\ar@{.>}[r]^{B} & {}_{\EOp}\C \\
{}_{\COp}\C\ar[u]\ar[r]^{B} & {}_{\COp}\C\ar[u] }.
\end{equation*}
As a byproduct,
we can form a composite functor
\begin{equation*}
\xymatrix{ {}_{\EOp}\C\ar@{.>}[r]^{B} & \cdots\ar@{.>}[r]^{B} & {}_{\EOp}\C \\
{}_{\COp}\C\ar[r]^{B}\ar[u] & \cdots\ar[r]^{B} & {}_{\COp}\C\ar[u] }
\end{equation*}
to extend the iterated bar construction of commutative dg-algebras
to the category of $\EOp$-algebras.
The definition of the iterated bar construction $B^n: A\mapsto B^n(A)$
can be generalized to any $E_\infty$-operad $\EOp$ (not necessarily cofibrant as a dg-operad),
because we can pick a cofibrant replacement $\epsilon: \QOp\xrightarrow{\sim}\EOp$
and use the associated restriction functor $\epsilon^*: {}_{\EOp}\C\rightarrow{}_{\QOp}\C$
to map $\EOp$-algebras into the category of algebras over~$\QOp$.

\medskip
The cochain complex $C^*(X)$
of a topological space (or simplicial set)~$X$
forms an instance of an algebra over an $E_\infty$-operad.
By theorems of~\cite{Bar1},
the iterated bar complex $B^n(C^*(X))$ is equivalent to $C^*(\Omega^n X)$,
the cochain algebra of the $n$-fold loop space of $X$,
provided that the space $X$ satisfies mild finiteness and completeness assumptions.
This result gives the original motivation for the extension of iterated bar complexes
to categories of algebras over $E_\infty$-operads.

The topological interpretation of iterated bar complexes (the actual objective of~\cite{Bar1})
uses the full multiplicative structure of the bar complex,
but the definition of the iterated bar construction as a composite functor $B^n: {}_{\EOp}\C\rightarrow{}_{\EOp}\C$
involves much more information than necessary for the determination of the differential of~$B^n(A)$.
Moreover,
the structure of a minimal cofibrant $E_\infty$-operad
models all secondary operations
which occur on the homology of $E_\infty$-algebras (see~\cite{ChataurLivernet}).
For this reason,
the determination of the iterated bar complex
by naive iterations of the bar construction
carries deep difficulties.

\medskip
The first purpose of this paper is to give a direct construction of the iterated bar complex,
within the framework of~\cite{Bar1},
but so that we avoid the iteration process of the definition
and the use of cofibrant replacements of $E_\infty$-operads.
Roughly,
we show that the definition of the iterated bar complex~$B^n(A)$
can be reduced to a construction of linear homological algebra
in the context of operads.

Let $\ROp$ be any operad.
In~\cite{Bar0},
we show that a functor $\Sym_{\ROp}(M,-): {}_{\ROp}\C\rightarrow\C$
is naturally associated to any right $\ROp$-module $M$
and all functors on $\ROp$-algebras
which are defined by composites of colimits and tensor products
have this form.

Let $\ROp = \EOp$ or $\ROp = \COp$.
In~\cite{Bar1},
we check that the bar construction is an instance of a functor of this form $B(-) = \Sym_{\ROp}(B_{\ROp},-)$
for a certain right $\ROp$-module $B_{\ROp}$.
In fact,
we prove the existence of multiplicative structures on the bar construction
at the module level.
An assertion of~\cite{Bar0} implies that the iterated bar construction $B^n(A)$,
defined by a composite of functors associated to right modules over operads,
forms itself a functor determined by a right $\ROp$-module.
This observation gives the starting point of the construction of this article.

Throughout the paper,
we use the notation $B^n_{\ROp}$
for the right $\ROp$-module which represents the iterated bar complex $B^n: {}_{\ROp}\C\rightarrow\C$.
We study the structure of~$B^n_{\ROp}$.
We check that $B^n_{\ROp}$ is defined by a pair $B^n_{\ROp} = (T^n\circ\ROp,\partial_{\ROp})$,
where:
\begin{itemize}
\item
the composite $T^n\circ\ROp$ represents a free right $\ROp$-module,
whose associated functor $\Sym_{\ROp}(T^n\circ\ROp): {}_{\ROp}\C\rightarrow\C$
is the iterated tensor coalgebra underlying the iterated bar construction;
\item
the term $\partial_{\ROp}$ refers to a twisting homomorphism of right $\ROp$-modules
$\partial_{\ROp}: T^n\circ\ROp\rightarrow T^n\circ\ROp$;
the differential of the iterated bar construction~$B^n(A)$
is defined by the addition of the twisting homomorphism
\begin{equation*}
\underbrace{\Sym_{\ROp}(T^n\circ\ROp,A)}_{= B^n(A)}
\xrightarrow{\Sym_{\ROp}(\partial_{\ROp},A)}\underbrace{\Sym_{\ROp}(T^n\circ\ROp,A)}_{= B^n(A)}
\end{equation*}
induced by $\partial_{\ROp}: T^n\circ\ROp\rightarrow T^n\circ\ROp$
to the natural differential
of the iterated tensor coalgebra.
\end{itemize}
For the commutative operad $\ROp = \COp$,
the definition of the twisting homomorphism $\partial_{\COp}$ arises from the standard definition
of the iterated bar construction
of commutative algebras.
For an $E_\infty$-operad $\ROp = \EOp$,
the twisting homomorphism $\partial_{\EOp}$
solves a lifting problem
\begin{equation*}
\xymatrix{ T^n\circ\EOp\ar[d]\ar@{.>}[r]^{\partial_{\EOp}} & T^n\circ\EOp\ar[d] \\
T^n\circ\COp\ar[r]_{\partial_{\COp}} & T^n\circ\COp }.
\end{equation*}

We prove the homotopy uniqueness of a solution of this lifting problem.
We also prove the existence of a solution by effective arguments.
For this purpose,
we use an easy generalization of classical techniques of linear homological algebra
in the context of right modules over operads.
The desired direct definition of the functor $B^n: {}_{\EOp}\C\rightarrow\C$
is obtained in this way.

In this construction we lose multiplicative structures,
but we can apply the argument of~\cite[Theorem 2.A]{Bar1} to redefine an $E_\infty$-multiplicative structure on $B^n(A)$,
in full or in part,
at any stage of iteration of the bar construction.
The uniqueness argument of~\cite{Bar1} ensures that we retrieve the good structure anyway.

\medskip
The main objective of the paper is to prove that the $n$-fold bar complex $B^n: {}_{\EOp}\C\rightarrow\C$
extends to the category of $\EOp_n$-algebras (but not further)
\begin{equation*}
\xymatrix{ {}_{\EOp_1}\C & \ar[l] & \ar@{.}[l] &
{}_{\EOp_n}\C\ar[l]\ar@{.>}@/_/[drrr]_{B^n} & \ar[l] &  \ar@{.}[l] &
\ar[l]{}_{\EOp}\C\ar[d]^{B^n} \\
&&&&&& \C }
\end{equation*}
for certain $E_\infty$-operads $\EOp$
equipped with a filtration of the form (*).
We obtain as a main resut that the $n$-fold desuspension of the iterated bar complex $\Sigma^{-n} B^n(A)$
determines the homology of~$\EOp_n$-algebras,
the homology theory $H^{\EOp_n}_*(A)$
defined abstractly as the homology of a derived indecomposable functor $L\Indec: \Ho({}_{\EOp_n}\C)\rightarrow\Ho(\C)$.

To define the iterated bar complex of~$\EOp_n$-algebras,
we observe that the twisting homomorphism $\partial_{\EOp}: T^n\circ{\EOp}\rightarrow T^n\circ{\EOp}$
admits a restriction
\begin{equation*}
\xymatrix{ T^n\circ\EOp_n\ar[d]\ar@{.>}[r]^{\partial_{\EOp_n}} & T^n\circ\EOp_n\ar[d] \\
T^n\circ\COp\ar[r]_{\partial_{\COp}} & T^n\circ\COp }
\end{equation*}
for a good choice of~$\partial_{\EOp}$.
Hence,
we can form a twisted right $\EOp_n$-module $B^n_{\EOp_n} = (T^n\circ\EOp_n,\partial_{\EOp_n})$
and the associated functor $\Sym_{\EOp_n}(B^n_{\EOp_n})$
gives the desired extension of the iterated bar complex.
Then we prove that $B^n_{\EOp_n} = (T^n\circ\EOp_n,\partial_{\EOp_n})$
forms a cofibrant resolution of a unit object $I$
in right $\EOp_n$-modules
and we use this observation to conclude that the functor $\Sym_{\EOp_n}(B^n_{\EOp_n})$
determines the $\EOp_n$-homology $H^{\EOp_n}_*(-)$.

\medskip
The iterated bar complexes are connected by suspension morphisms
\begin{equation*}
\sigma: \Sigma^{1-n} B^{n-1}(A)\rightarrow\Sigma^{-n} B^n(A)
\end{equation*}
and we can perform a colimit
to associate an infinite bar complex
\begin{equation*}
\Sigma^{-\infty} B^{\infty}(A) = \colim_n\Sigma^{-n} B^n(A)
\end{equation*}
to any $\EOp$-algebra $A$.
The relationship $H^{\EOp_n}_*(A) = H_*(\Sigma^{-n} B^n(A))$ also holds in the case $n = \infty$.

The $\Gamma$-homology of~\cite{Robinson} and the $E_\infty$ Andr\'e-Quillen homology of~\cite{Mandell}
are other definitions of the homology theory associated to an $E_\infty$-operad.
Our result implies that the $\Gamma$-homology with trivial coefficient
agrees with the homology of the infinite bar complex $\Sigma^{-\infty} B^{\infty}(A)$,
for any $E_\infty$ algebra $A$ (provided that $A$ is cofibrant as a dg-module),
and similarly as regards the dg-version of the $E_\infty$ Andr\'e-Quillen homology of~\cite{Mandell}.
This relationship between $\Gamma$-homology and iterated bar complexes
does not seem to occur in the literature
even in the case of commutative algebras,
for which we can apply the classical definition of iterated bar constructions.

\medskip
The identity $H^{\EOp_n}(A) = H_*(\Sigma^{-n} B^n(A))$
enables us to deduce the $\EOp_n$-homology of usual commutative algebras (polynomial algebras, exterior algebras, divided power algebras, \dots)
from results of~\cite{Cartan}.
In the case $n = \infty$,
this approach could be used to give explicit representatives of $\Gamma$-homology classes
and to improve on results of~\cite{Richter}.

In general,
we have a natural spectral sequence
\begin{equation*}
E^1 = \Sigma^{-n} B^n(H_*(A))\Rightarrow H_*(\Sigma^{-n} B^n(A)) = H^{\EOp_n}_*(A),
\end{equation*}
whose $E^1$-term reduces to the usual $n$-fold bar construction of a commutative algebra for any $n>1$
(the homology of an $\EOp_n$-algebra forms a commutative algebra for $n>1$, an associative algebra for $n=1$).
In the case of the cochain algebra $A = C^*(X)$ of a space $X$,
we conjecture that this spectral sequence agrees, from $E^2$-stage,
with a spectral sequence of~\cite{AhearnKuhn}
which is defined with Goodwillie's calculus of functors.

On one hand,
one might gain quantitative information on the cohomology of iterated loop spaces
from the study of such spectral sequences,
arising from filtrations of iterated bar complexes.
On the other hand,
the connection between the homology of iterated bar complexes and the homology of $E_n$-algebras
could be used to gain qualitative information
of a new nature on the cohomology of iterated loop spaces
-- notably, our result implies that certain groups of homotopy automorphisms of $E_n$-operads
act on this cohomology and we conceive, from insights of~\cite{Kontsevich},
that this gives an action of higher versions of the Grothendieck-Teichm\"uller group
on the cohomology of iterated loop spaces.
Besides,
for a sphere $X = S^{n-m}$,
we have a chain complex, defined purely algebraically (in terms of characteristic structures of $E_n$-operads),
computing the cohomology $H_*^{\EOp_n}(C^*(S^{n-m}))$ (see~\cite{Bar5}).
Thus
spheres are first examples of spaces for which our approach seems appropriate
and for which we plan to study the applications of our results.

\medskip
Throughout the article,
we study the application of our constructions to the Barratt-Eccles operad,
a nice combinatorial $E_\infty$-operad $\EOp$
equipped with a filtration
of the form (*).
In fact,
we will observe that this $E_\infty$-operad
has all extra structures needed for the constructions of the article.
Besides,
we proved in~\cite{BergerFresse}
that the Barratt-Eccles acts on cochain complexes of spaces.
Thus the Barratt-Eccles operad
is also well suited for the topological applications
of our results.

\section*{Plan}
This paper is a sequel of~\cite{Bar1},
but the reader is essentially supposed to have some familiarity with the theory of modules over operads,
which gives the background of our constructions,
and with the Koszul duality of operads (over a ring),
which is used in homology calculations.
The overall setting is reviewed in the preliminary part of the paper, \emph{``Conventions and background''},
at least to fix conventions.

The next part, \emph{``The construction of iterated bar complexes''},
is devoted to the definition of the iterated bar complex
of algebras over $E_{\infty}$-operads
and to the extension of the construction to algebras over $E_n$-operads.
In~\emph{``Iterated bar complexes and homology theories''},
we prove that the $n$-fold bar complex determines the homology of algebras
over an $E_n$-operad.
In the concluding part, \emph{``Applications to the cohomology of iterated loop spaces''},
we explain the applications of our results
to iterated loop space cohomology.
We refer to the introduction of each part for a more detailed outline.

References to the appendix \emph{``Iterated bar modules and level tree posets''}
of a former version of this paper have to be redirected to the addendum~\cite{Bar2Appendix}.
The purpose of this addendum is to explain the relationship between Batanin's categories of pruned trees (see~\cite{BataninHigherOperads,BataninConfigurations})
and iterated bar complexes
and to revisit our constructions in this formalism.
The reader can use this reference~\cite{Bar2Appendix}
as an informal introduction to the constructions of the article.

\section*{Main theorems}

{\parindent=-1cm\leftskip=1cm

\textbf{Theorem~\ref{IteratedBarModules:IteratedBarEquivalence} (p.~\pageref{IteratedBarModules:IteratedBarEquivalence}):}
Definition of the iterated bar construction of algebras over an $E_\infty$-operad $\EOp$
by a certain cofibrant replacement in the category of right $\EOp$-modules.

\textbf{Theorem~\ref{QuasiFreeLifting:BarTwistingCochainLifting} (p.~\pageref{QuasiFreeLifting:BarTwistingCochainLifting}):}
Definition of the iterated bar module as a lifting of quasi-free modules.

\textbf{Theorems~\ref{EnDefinition:TwistingCochainRestriction}-\ref{EnDefinition:Result} (p.~\pageref{EnDefinition:TwistingCochainRestriction}):}
Extension of the $n$-fold bar construction to algebras over an $E_n$-operad.

\textbf{Theorems~\ref{EnHomology:IteratedBarModuleHomotopy}-\ref{EnHomology:Conclusion} (p.~\pageref{EnHomology:IteratedBarModuleHomotopy}):}
The $n$-fold bar complex determines the homology theory associated to $E_n$-operads $H^{\EOp_n}_*(-)$.

\textbf{Theorem~\ref{EinfinityHomology:Result} (p.~\pageref{EinfinityHomology:Result}):}
The infinite bar complex determines the homology theory associated to $E_{\infty}$-operads,
equivalently the $\Gamma$-homology with trivial coefficients.

\textbf{In the concluding part (p.~\pageref{Conclusion:TopologicalApplications}):}
For the cochain algebra of a space $X$,
the homology $H^{\EOp_n}_*(C^*(X))$ determines, under mild finiteness and completeness assumptions,
the cohomology of the $n$-fold loop space~$\Omega^n X$.

}

\part*{Conventions and background}\label{Background}

The structure of a module over an operad is used to model functors on algebras over operads.
The purpose of this part is to review this overall setting with the aim of fixing our framework.

The point of view adopted in this paper is borrowed from the book~\cite{Bar0}
to which we refer for a comprehensive account of the background of our constructions.
The r\'esum\'e of this part would be sufficient to make the conceptual setting of the paper accessible to readers
who are already familiar with usual definitions of the theory of operads.

In~\cite{Bar0,Bar1},
we only use the standard definition of an operad,
according to which the elements of an operad model operations
with $r$ inputs indexed by integers $i\in\{1,\dots,r\}$, for any $r\in\NN$.
But there is another usual definition of the structure of an operad
in which the inputs of operations
are indexed by any finite set $\eset = \{e_1,\dots,e_r\}$.
The indexing by finite sets is more natural for certain constructions of the article.
Therefore,
we revisit a few definitions of~\cite{Bar0,Bar1}
in this formalism.

But, first of all,
we recall the categorical settings of~\cite{Bar1}
that we keep for this article.

\subsubsection{Categorical settings}\label{Background:CategoricalSettings}
A commutative ground ring $\kk$ is fixed once and for all.
In applications,
we take $\kk = \QQ$, the field of rationals,
or $\kk = \FF_p$, a finite primary field,
or $\kk = \ZZ$, the ring of integers,
but no assumption on the ground ring is really required (outside commutativity).

We take the category of differential graded modules over $\kk$
as a base symmetric monoidal category
and we use the notation $\C$ to refer to this category.
For us,
a differential graded module (a dg-module for short) consists of a $\ZZ$-graded module $C$
equipped with a differential $\delta: C\rightarrow C$
that lowers degrees by $1$.

The letter $\E$ refers to a symmetric monoidal category over~$\C$,
whose structure includes a unit object $\unit\in\E$,
an external tensor product $\otimes: \C\times\E\rightarrow\E$,
and an internal tensor product $\otimes: \E\times\E\rightarrow\E$,
which satisfy the usual relations of symmetric monoidal categories.
In this paper,
we take either $\E = \C$, the category of dg-modules itself,
or $\E = \M$, the category of $\Sigma_*$-objects in $\C$,
or $\E = \M{}_{\ROp}$, the category of right modules over an operad $\ROp$.
The definition of the second-mentioned categories is recalled next.

Other examples include the categories of dg-modules
over a graded $\kk$-algebra $R$
(see~\S\ref{EinfinityHomology:GammaHomology}).

In the paper,
we use an external hom-functor $\Hom_{\E}(-,-): \E^{op}\times\E\rightarrow\C$
characterized by the adjunction relation $\Mor_{\E}(C\otimes E,F) = \Mor_{\C}(C,\Hom_{\E}(E,F))$,
where $C\in\C$, $E,F\in\E$.
The elements of the dg-module $\Hom_{\E}(E,F)$ are called homomorphisms
to be distinguished from the actual morphisms of~$\E$.
In the case $\E = \C$,
the dg-module $\Hom_{\E}(E,F)$ is spanned in degree $d$ by the morphisms of $\kk$-modules $f: E\rightarrow F$
which raise degree by $d$.
This explicit definition has a straightforward generalization
for the other categories $\E = \M,\M{}_{\ROp}$
whose definition is recalled next.

\subsubsection{Functors on finite sets and $\Sigma_*$-modules}\label{Background:SetIndexing}
The category of $\Sigma_*$-modules $\M$
consists of collections $M = \{M(r)\}_{r\in\NN}$
whose term $M(r)$
is an object of the base category $\C$ (for us, the category of dg-modules)
equipped with an action of the symmetric group on $r$ letters $\Sigma_r$.

In certain applications,
we use that the category of $\Sigma_*$-modules $\M$ is equivalent to the category of functors $F: \Bij\rightarrow\C$,
where $\Bij$ refers to the category formed by finite sets as objects
and bijective maps as morphisms
(this equivalence is borrowed from~\cite{GetzlerJones,GinzburgKapranov},
see also the surveys of~\cite[\S 1.1.8]{FressePartitions} and~\cite[\S 1.7]{MarklShniderStasheff}).

In one direction,
for a functor $F: \Bij\rightarrow\C$,
the dg-module $F(\{1,\dots,r\})$ associated to the set $\eset = \{1,\dots,r\}$
inherits an action of the symmetric group $\Sigma_r$
since a permutation $w\in\Sigma_r$ is equivalent to a bijection $w: \{1,\dots,r\}\rightarrow\{1,\dots,r\}$.
Hence,
the collection $F(r) = F(\{1,\dots,r\})$, $r\in\NN$,
forms a $\Sigma_*$-module naturally associated to $F: \Bij\rightarrow\C$.

In the other direction,
for a given $\Sigma_*$-module $M$,
we set:
\begin{equation*}
M(\eset) = \Bij(\{1,\dots,r\},\eset)\otimes_{\Sigma_r} M(r),
\end{equation*}
for any set with $r$ elements $\eset = \{e_1,\dots,e_r\}$.
The tensor product $S\otimes C$ of a dg-module $C$ with a finite set $S$ is defined as the sum of copies of $C$
indexed by the elements of $S$.
The quotient over $\Sigma_r$ makes the natural $\Sigma_r$-action on $M(r)$
agree with the action of permutations by right translations on $\Bij(\{1,\dots,r\},\eset)$.
The map $M: \eset\mapsto M(\eset)$
defines a functor naturally associated to $M$.

Intuitively,
an element $x\in M(r)$, where $r\in\NN$,
can be viewed as an operation with $r$ inputs indexed by $(1,\dots,r)$;
an element $x\in M(\eset)$, where $\eset$ is a finite set $\eset = \{e_1,\dots,e_r\}$,
represents an operation $x = x(e_1,\dots,e_r)$
whose inputs are indexed by the elements of~$\eset$.
In the definition of the functor associated to a $\Sigma_*$-module $M$,
we use a bijection $u\in\Bij(\{1,\dots,r\},\eset)$
to reindex the inputs of an operation $x\in M(r)$
by elements of~$\eset$.

Now,
the subtlety is that we may use morphisms $f: M(r)\rightarrow N(r)$, like chain-homotopies,
which do not preserve symmetric group actions.
In this context,
we have to assume that the finite set $\eset$ comes equipped with a bijection $u: \{1,\dots,r\}\rightarrow\eset$
in order to define a morphism $f: M(\eset)\rightarrow N(\eset)$
associated to $f$.
In applications,
the bijection $u: \{1,\dots,r\}\rightarrow\eset$
is determined by an ordering $\eset = \{e_1<\dots<e_r\}$
that $\eset$ naturally inherits from a larger set $\eset\subset\fset = \{f_1<\dots<f_s\}$.

In many usual situations,
we specify the bijection $u: \{1,\dots,r\}\rightarrow\eset$
by the sequence of values $(u(1),\dots,u(r))$.

\subsubsection{Indexing by finite sets and the tensor product}\label{Background:TensorProduct}
The definition of the tensor product of $\Sigma_*$-modules
is recalled in~\cite[\S\S 2.1.5-2.1.7]{Bar0}.
The expansion given in this reference has a nice reformulation in terms of finite set indexings:
the functor equivalent to the tensor product $M\otimes N$ of $\Sigma_*$-modules $M,N\in\M$
is defined on any finite set $\eset$
by the direct sum
\begin{equation*}
(M\otimes N)(\eset) = \bigoplus_{\uset\amalg\vset = \eset} M(\uset)\otimes N(\vset),
\end{equation*}
where the pair $(\uset,\vset)$ runs over partitions of~$\eset$.
In this light,
the elements of $(M\otimes N)(\eset)$ can be viewed as tensors $x(u_1,\dots,u_r)\otimes y(v_1,\dots,v_s)$
of elements $x\in M$ and $y\in N$
together with a sharing $\eset = \{u_1,\dots,u_r\}\amalg\{v_1,\dots,v_s\}$
of the indices of~$\eset$.
Throughout this paper,
we adopt this representation of elements in tensor products,
but we usually drop indices to simplify notation: $x\otimes y = x(u_1,\dots,u_r)\otimes y(v_1,\dots,v_s)$.

In~\cite{Bar0},
we represent the elements of $M\otimes N$ by point-tensors $w\cdot x\otimes y\in(M\otimes N)(r+s)$,
where $w\in\Sigma_{r+s}$, $x\otimes y\in M(r)\otimes N(s)$
(see more specifically~\S 0.5 and~\S 2.1.9 of~\emph{loc. cit.}).
In the formalism of this paragraph,
the point-tensor $w\cdot x\otimes y$
is equivalent to the tensor
\begin{equation*}
x(u_1,\dots,u_r)\otimes y(v_1,\dots,v_s)\in(M\otimes N)(\{1,\dots,r+s\})
\end{equation*}
such that $(u_1,\dots,u_r) = (w(1),\dots,w(r))$ and $(v_1,\dots,v_s) = (w(r+1),\dots,w(r+s))$.

Note that each summand of a partition $\eset = \{u_1,\dots,u_r\}\amalg\{v_1,\dots,v_s\}$
inherits a canonical ordering
if $\eset$ forms itself an ordered set.

The tensor product of $\Sigma_*$-modules
is equipped with a symmetry isomorphism
\begin{equation*}
\tau: M\otimes N\rightarrow N\otimes M,
\end{equation*}
which can be defined componentwise
by the symmetry isomorphism
$\tau: M(\uset)\otimes N(\vset)\rightarrow N(\vset)\otimes M(\uset)$
inherited from the category of dg-modules.
The tensor product of $\Sigma_*$-modules is also obviously associative
and has the $\Sigma_*$-module such that
\begin{equation*}
\unit(r) = \begin{cases} \kk, & \text{if $r = 0$}, \\ 0, & \text{otherwise}, \end{cases}
\end{equation*}
as a unit object.
Besides,
we have an exterior tensor product $\otimes: \C\times\M\rightarrow\M$
such that $(C\otimes M)(r) = C\otimes M(r)$
for any $C\in\C$, $M\in\M$.
Thus
we finally obtain that the category of $\Sigma_*$-modules forms a symmetric monoidal category
over the base category of dg-modules.

\subsubsection{The composition structure of $\Sigma_*$-modules and operads}\label{Background:CompositionStructure}
Let $\E$ be any symmetric monoidal category over $\C$.
Each $\Sigma_*$-module $M\in\M$
gives rise to a functor $\Sym(M): \E\rightarrow\E$
which maps an object $E\in\E$
to the module of symmetric tensors with coefficients in $M$:
\begin{equation*}
\Sym(M,E) = \bigoplus_{r=0}^{\infty} (M(r)\otimes E^{\otimes r})_{\Sigma_r}.
\end{equation*}
In this construction,
the coinvariants $(-)_{\Sigma_r}$ identify the natural action of permutations on tensors
with the natural $\Sigma_r$-action on $M(r)$.
The map $\Sym: M\mapsto\Sym(M)$ defines a functor from the category of $\Sigma_*$-objects $\M$
to the category of functors $F: \E\rightarrow\E$.

The tensor product of the category of $\Sigma_*$-modules $\otimes: \M\times\M\rightarrow\M$
reflects the pointwise tensor product of functors $\Sym(M\otimes N,E) = \Sym(M,E)\otimes\Sym(N,E)$,
and similarly as regards the tensor product over dg-modules.

The category of $\Sigma_*$-modules is also equipped with a composition product $\circ: \M\times\M\rightarrow\M$,
characterized by the relation $\Sym(M\circ N) = \Sym(M)\circ\Sym(N)$, for $M,N\in\M$,
and we have a unit object
\begin{equation*}
I(r) = \begin{cases} \kk, & \text{if $r = 1$}, \\ 0, & \text{otherwise}, \end{cases}
\end{equation*}
such that $\Sym(I) = \Id$.
The composition product can be defined by the formula $M\circ N = \bigoplus_{r=0}^{\infty} (M(r)\otimes N^{\otimes r})_{\Sigma_r}$,
where we form the tensor power of~$N$ within the category of $\Sigma_*$-modules
and we use the tensor product of $\Sigma_*$-modules over dg-modules
to form the tensor product with $M(r)$.
Equivalently,
we have a formula $M\circ N = \Sym(M,N)$,
where we apply the definition of the functor $\Sym(M): \E\rightarrow\E$
to the category $\Sigma_*$-modules $\E = \M$.

The structure of an operad $\POp$
is defined by a composition product $\mu: \POp\circ\POp\rightarrow\POp$
together with a unit morphism $\eta: I\rightarrow\POp$
that satisfy the usual relations of monoid objects.

Modules over operads are defined naturally by using the composition structure of $\Sigma_*$-modules:
a left module over an operad $\POp$
consists of a $\Sigma_*$-module $N$
equipped with a left $\POp$-action determined by a morphism $\lambda: \POp\circ N\rightarrow N$;
a right module over an operad $\ROp$
is defined symmetrically as a $\Sigma_*$-module $M$
equipped with a right $\ROp$-action determined by a morphism $\rho: M\circ\ROp\rightarrow M$;
a bimodule over operads $(\POp,\ROp)$
is a $\Sigma_*$-module $N$
equipped with both a right $\ROp$-action $\rho: N\circ\ROp\rightarrow N$
and a left $\POp$-action $\lambda: \POp\circ N\rightarrow N$
that commute to each other.

We refer to~\cite{Bar0} for a comprehensive study of modules over operads
and further bibliographical references on recollections of this paragraph.

\subsubsection{Pointwise composition products}\label{Background:TermwiseComposites}
In the original definition of the structure of an operad~\cite{May},
the composition product $\mu: \POp\circ\POp\rightarrow\POp$
is determined by a collection of morphisms
\begin{equation*}
\mu: \POp(r)\otimes\POp(n_1)\otimes\dots\otimes\POp(n_r)\rightarrow\POp(n_1+\dots+n_r),
\end{equation*}
where $r,n_1,\dots,n_r\in\NN$.
This definition is also used in the paper.
The equivalence between the latter definition and the definition of~\S\ref{Background:CompositionStructure}
comes from an explicit expansion of the composition product of $\Sigma_*$-modules
(see for instance~\cite[\S 2.2]{Bar0}, brief recollections are also given in~\S\ref{Background:TensorCompositeDistribution}).

The image of $p\in\POp(r)$ and $q_1\in\POp(n_1),\dots,q_r\in\POp(n_r)$
under the termwise composition product
is usually denoted by $p(q_1,\dots,q_r)$.
In the paper,
we also use a generalization of the definition of the composite $p(q_1,\dots,q_r)$
for elements $q_1\in\POp(\eset_1),\dots,q_r\in\POp(\eset_r)$,
where $(\eset_1,\dots,\eset_r)$ are any finite sets.
In this situation,
the composite $p(q_1,\dots,q_r)$
returns an element of $\POp(\eset_1\amalg\dots\amalg\eset_r)$.

Similar definitions apply to modules over operads.

\subsubsection{Categories of modules and operads}\label{Background:OperadModules}
Throughout the paper,
we use the notation $\M{}_{\ROp}$ for the category of right $\ROp$-modules,
the notation ${}_{\POp}\M$ for the category of left $\POp$-modules,
and the notation ${}_{\POp}\M{}_{\ROp}$ for the category of $(\POp,\ROp)$-bimodules.

The usual definitions of linear algebra (free objects, extension and restriction functors)
have a natural extension in the context of modules over operads (see relevant sections of~\cite{Bar0}).
In particular, a relative composition product $M\circ_{\ROp} N$
is associated to any pair $(M,N)$
such that
$M$ is equipped with a right $\ROp$-action
and
$N$ is equipped with a left $\ROp$-action.
The object $M\circ_{\ROp} N$
is defined by the reflexive coequalizer of the morphisms $d_0,d_1: M\circ\ROp\circ N\rightrightarrows M\circ N$
induced by the right $\ROp$-action on $M$
and the left $\ROp$-action on $N$ (see for instance~\cite[\S 5.1.5, \S 9.2.4]{FressePartitions}
or~\cite[\S 2.1.7]{Bar0}, we also refer to the bibliography of~\cite{FressePartitions}
for further references on this definition).

The extension and restriction functors associated to an operad morphism $\phi: \POp\rightarrow\QOp$,
respectively $\psi: \ROp\rightarrow\SOp$,
are denoted by:
\begin{equation*}
\phi_!: {}_{\POp}\M\rightleftarrows{}_{\QOp}\M :\phi^*,
\quad\text{respectively}
\quad\psi_!: \M{}_{\ROp}\rightleftarrows\M{}_{\SOp} :\psi^*.
\end{equation*}
In the context of bimodules,
we have extensions and restrictions of structure
on the left and on the right.
These extension and restriction functors are also denoted by:
\begin{equation*}
\phi_!: {}_{\POp}\M{}_{\ROp}\rightleftarrows{}_{\QOp}\M{}_{\ROp} :\phi^*,
\quad\text{respectively}
\quad\psi_!: {}_{\POp}\M{}_{\ROp}\rightleftarrows{}_{\POp}\M{}_{\SOp} :\psi^*.
\end{equation*}
The extension functors are given by relative composition products of the form:
\begin{equation*}
\phi_! N = \QOp\circ_{\POp} N,
\quad\text{respectively}
\quad\psi_! M = M\circ_{\ROp}\SOp.
\end{equation*}
In formulas,
we usually omit specifying structure restrictions and we use the expression of the relative composite
to denote structure extensions rather than the functor notation.
Nevertheless
we do use the functor notation $(\phi_!,\phi^*)$ to refer to the extension and restriction of structure
as functors between modules categories.

For a $\Sigma_*$-module $K$,
the composite $\POp\circ K$ inherits a left $\POp$-action
and represents the free object generated by $K$ in the category of left $\POp$-modules.
The symmetrical composite $K\circ\ROp$
represents the free object generated by $K$ in the category of right $\ROp$-modules.

Recall that the composition product is not symmetric and does not preserves colimits on the right.
For that reason,
categories of left modules
differ in nature from categories of right modules
over operads.

\subsubsection{Algebras over operads}\label{Background:OperadAlgebras}
By definition,
an algebra over an operad $\POp$ consists of an object $A\in\E$
together with an evaluation product $\lambda: \Sym(\POp,A)\rightarrow A$
which satisfies natural associativity and unit relations
with respect to the composition product and the unit of $\POp$.
The evaluation product $\lambda: \Sym(\POp,A)\rightarrow A$
is equivalent to a collection
of morphisms $\lambda: \POp(r)\otimes A^{\otimes r}\rightarrow A$.

In the context of dg-modules,
the evaluation morphism associates an actual operation $p: A^{\otimes r}\rightarrow A$
to every element of the operad $p\in\POp(r)$,
and we use the notation $p(a_1,\dots,a_r)$
to refer to the image of a tensor $a_1\otimes\dots\otimes a_r\in A^{\otimes r}$
under this operation $p: A^{\otimes r}\rightarrow A$.

Note that the definition of a $\POp$-algebra
makes sense in any symmetric monoidal category $\E$ over $\C$,
and not only in the category of dg-modules itself.
We adopt the notation ${}_{\POp}\E$
to refer to the category of $\POp$-algebras in $\E$.

In the next paragraphs,
we explain that this definition of the category of $\POp$-algebras in a symmetric monoidal category over $\C$
applies, besides the category of dg-modules itself, to the category of $\Sigma_*$-objects
and to categories of right modules over an operad.
In this paper,
we only use these particular symmetric monoidal categories over $\C$,
but we need the overall idea of a $\POp$-algebra in a symmetric monoidal category over the base category
to make constructions more conceptual.

Recall that an operad morphism $\phi: \POp\rightarrow\QOp$
determines a pair of adjoint extension and restriction functors $\phi_!: {}_{\POp}\E\rightleftarrows{}_{\QOp}\E :\phi^*$,
defined like the extension and restriction functors of modules over operads (see~\cite[\S 3.3]{Bar0}).

\subsubsection{Tensor products and composition structures}\label{Background:TensorCompositeDistribution}
We recall in~\S\ref{Background:TensorProduct}
that the category of $\Sigma_*$-modules comes equipped with a tensor product
and forms a symmetric monoidal category over the base category of dg-modules.
As a consequence,
we can associate a category of algebras in $\Sigma_*$-modules to any operad $\POp$.
In fact,
since we have an identity $\Sym(M,N) = M\circ N$
between the functor $\Sym(M): \E\rightarrow\E$ associated to $M$
and the composition with $M$ in $\Sigma_*$-modules,
we have a formal equivalence between the structure of a $\POp$-algebra in $\Sigma_*$-modules,
determined by a morphism $\lambda: \Sym(\POp,N)\rightarrow N$,
and the structure of a left $\POp$-module,
which is itself determined by a morphism $\lambda: \POp\circ N\rightarrow N$.

The tensor product of $\Sigma_*$-modules
satisfies a distribution relation $(M\otimes N)\circ P\simeq(M\circ P)\otimes(N\circ P)$
with respect to the composition product.
From this observation,
we deduce that a tensor product of right $\ROp$-modules
inherits a natural right $\ROp$-action (see~\cite[\S 6]{Bar0})
and
we obtain that the category of right modules $\M{}_{\ROp}$
forms a symmetric monoidal category over $\C$,
like the category of $\Sigma_*$-modules $\M$.
As a consequence,
we can also apply the ideas of~\S\ref{Background:OperadAlgebras}
to the category right $\ROp$-modules $\E = \M{}_{\ROp}$
and we can associate a category of algebras in right modules $\ROp$
to any operad $\POp$.
Actually,
the composite $M\circ N$
inherits a natural right $\ROp$-action when $N$ is a right $\ROp$-modules
and the identity $\Sym(M,N) = M\circ N$ holds in the category of right $\ROp$-modules,
where we apply the definition of the functor $\Sym(M): \E\rightarrow\E$
to the category of right $\ROp$-modules $\E = \M{}_{\ROp}$.
From these observations,
we deduce readily that the category of $\POp$-algebras in right $\ROp$-modules
is formally equivalent to the category of $(\POp,\ROp)$-bimodules.

The associativity of the composition product of $\Sigma_*$-modules
is equivalent to the distribution relation $\Sym(M,N)\circ P\simeq\Sym(M,N\circ P)$
on the module of symmetric tensors~$\Sym(M,N)$.

The definition of $\POp$-algebras in terms of evaluation morphisms $\lambda: \POp(r)\otimes A^{\otimes r}\rightarrow A$
and pointwise operations
applies to the context of $\Sigma_*$-modules and right-modules over operads
too.
In the case of a $\Sigma_*$-module $N$,
the evaluation morphisms $\lambda: \POp(r)\otimes N^{\otimes r}\rightarrow N$
are formed by using the tensor product of $\Sigma_*$-modules,
whose definition is recalled in~\S\ref{Background:TensorProduct}.
Thus,
the operation $p: N^{\otimes r}\rightarrow N$
associated to an element $p\in\POp(r)$
consists of a collection of dg-module homomorphisms
\begin{equation*}
p: N(\eset_1)\otimes\dots\otimes N(\eset_r)\rightarrow N(\eset_1\amalg\dots\amalg\eset_r)
\end{equation*}
indexed by partitions $\eset = \eset_1\amalg\dots\amalg\eset_r$.
We also use the notation $p(a_1,\dots,a_r)$
to represent the evaluation of operations $p\in\POp(r)$
on point-tensors in $\Sigma_*$-modules.

In the case of a right $\ROp$-module $N$,
the operation $p: N^{\otimes r}\rightarrow N$
associated to an element $p\in\POp(r)$
is simply assumed to preserve right $\ROp$-actions.

The equivalence between $\POp$-algebras in $\Sigma_*$-modules, respectively in right $\ROp$-modules,
and left $\POp$-modules, respectively of a $(\POp,\ROp)$-bimodules,
is pointed out in~\cite[\S\S 3.2.9-3.2.10,\S 9]{Bar0}.
Depending on the context,
we use either the idea of algebras in symmetric monoidal categories
or the language of modules of operads,
because each point of view has its own interests.

An operad $\POp$ forms obviously an algebra over itself in the category of right modules over itself.
In this case,
the evaluation operation $p: \POp(\eset_1)\otimes\dots\otimes\POp(\eset_r)\rightarrow\POp(\eset_1\amalg\dots\amalg\eset_r)$
is nothing but the composition product of~\S\ref{Background:TermwiseComposites}.

\subsubsection{Free algebras over operads}\label{Background:FreeOperadAlgebras}
The object $\Sym(\POp,E)$ associated to any $E\in\E$
inherits a natural $\POp$-algebra structure
and represents the free object generated by $E$
in the category of $\POp$-algebras ${}_{\POp}\E$.
In the paper,
we use the notation $\POp(E) = \Sym(\POp,E)$ to refer to the object $\Sym(\POp,E)$
equipped with this free $\POp$-algebra structure.

In the case $\E = \M$,
we have an identity between the free $\POp$-algebra
and the free left $\POp$-module generated by a $\Sigma_*$-module $M$,
and similarly in the context of right modules over an operad.
In the paper,
we use both representations $\POp(M) = \POp\circ M$
for these free objects
since each representation has its own interests.

The associativity of composition products
is equivalent to a distribution relation $\POp(M)\circ N = \POp(M\circ N)$.

\subsubsection{Modules over operads and functors}\label{Background:OperadFunctors}
Recall that a $\Sigma_*$-module $M$
determines a functor $\Sym(M): \E\rightarrow\E$.
For a right $\ROp$-module $M$,
we have a functor $\Sym_{\ROp}(M): {}_{\ROp}\E\rightarrow\E$
from the category of $\ROp$-algebras ${}_{\ROp}\E$
to the underlying category $\E$.
The object $\Sym_{\ROp}(M,A)\in\E$
associated to an $\ROp$-algebra $A\in{}_{\ROp}\E$
is defined by the reflexive coequalizers of morphisms $d_0,d_1: \Sym(M\circ\ROp,A)\rightrightarrows\Sym(M,A)$
induced by the right $\ROp$-action on $M$
and the left $\ROp$-action on $A$ (see~\cite[\S 5.1]{Bar0}).

For a left $\POp$-module $N\in{}_{\POp}\M$,
the objects $\Sym(N,E)$ inherit a natural left $\POp$-action
so that the map $\Sym(N): E\mapsto\Sym(N,E)$ defines a functor $\Sym(N): \E\rightarrow{}_{\POp}\E$
from the underlying category $\E$ to the category of $\POp$-algebras ${}_{\POp}\E$ (see~\cite[\S 3.2]{Bar0}).
For a $(\POp,\ROp)$-bimodule $N\in{}_{\POp}\M{}_{\ROp}$,
we have a functor $\Sym_{\ROp}(N): {}_{\ROp}\E\rightarrow{}_{\POp}\E$,
from the category of $\ROp$-algebras ${}_{\ROp}\E$
to the category of $\POp$-algebras ${}_{\POp}\E$ (see~\cite[\S 9.2]{Bar0}).

The relative composition product $M\circ_{\ROp} N$
reflects the composition of functors associated to modules over operads:
for a $(\POp,\ROp)$-bimodule $M$ and an $(\ROp,\QOp)$-bimodule~$N$,
we have a natural isomorphism
\begin{equation*}
\Sym_{\QOp}(M\circ_{\ROp} N)\simeq\Sym_{\ROp}(M)\circ\Sym_{\QOp}(N),
\end{equation*}
and similarly if we assume that $M$ or $N$ has a right (respectively, left) $\ROp$-module structure only (see for instance~\cite[\S 9.2]{Bar0}).
The relative composition product $M\circ_{\ROp} N$
can also be identified with the object $\Sym_{\ROp}(M,N)$
associated to $N$
by the functor $\Sym_{\ROp}(M): \E\rightarrow\E$ determined by $M$
on the category of right $\QOp$-modules $\E = \M{}_{\QOp}$.

As an illustration,
recall that an operad morphism $\phi: \POp\rightarrow\QOp$
defines adjoint extension and restriction functors
between categories of algebras over operads:
$\phi_!: {}_{\POp}\E\rightleftarrows{}_{\QOp}\E :\phi^*$.
By~\cite[Theorem 7.2.2]{Bar0},
we have an identity $\phi_! A = \Sym_{\POp}(\phi^*\QOp,A)$, for any $\POp$-algebra $A\in{}_{\POp}\E$,
where we use a restriction of structure on the right to identify the operad $\QOp$
with a $(\QOp,\POp)$-bimodule.
Symmetrically,
we have an identity $\phi^* B = \Sym_{\QOp}(\phi^*\QOp,B)$, for any $\QOp$-algebra $B\in{}_{\QOp}\E$,
where we use a restriction of structure on the left to identify the operad $\QOp$
with a $(\POp,\QOp)$-bimodule.

In~\cite[\S 7.2]{Bar0},
we also prove that extensions and restrictions of modules over operads
correspond, at the functor level,
to composites with extension and restriction functors on categories of algebras over operads.
To be explicit,
we have a functor identity $\Sym_{\ROp}(N,-) = \Sym_{\ROp}(N,\psi_!-)$,
for every right $\SOp$-module $N$,
and a functor identity $\Sym_{\SOp}(M\circ_{\ROp}\SOp,-) = \Sym_{\ROp}(M,\psi^*-)$
for every right $\ROp$-modules $M$ and every $\SOp$-algebra $B$,
where we consider the extension and restriction functors associated to an operad morphism $\psi: \ROp\rightarrow\SOp$.
Similar commutation formulas hold when we perform extensions and restrictions of modules on the left.

\subsubsection{The characterization of functors associated to modules over operads}\label{Background:FunctorOperations}
The composition of functors and the extension and retriction operations
are not the only categorical operations on functors
which can be represented at the module level.
In~\cite[\S\S 5-7]{Bar0},
we prove that the functor $\Sym_{\ROp}: M\mapsto\Sym_{\ROp}(M)$
commutes with colimits and tensor products of right $\ROp$-modules,
like the functor $\Sym: M\mapsto\Sym(M)$
on $\Sigma_*$-modules.
Similarly statements occur when we consider functors associated to left and bimodules over operads.

To retrieve the right module underlying a functor $F = \Sym_{\ROp}(M)$,
we use the simple idea that the application of~$F$
to the operad $\ROp$, viewed as an algebra over itself in the category of right modules over itself,
defines a right $\ROp$-module $F(\ROp) = \Sym_{\ROp}(M,\ROp)$
which is naturally isomorphic to $M$.
In terms of relative composition product, this identity reads $M\circ_{\ROp}\ROp = M$.

In the sequel,
we often switch from modules to functors when a construction becomes easier in the functor setting.
In any case,
we only deal with functors formed by tensor products, colimits, and operations
which have a representative at the module level
so that all our functors are properly modelled by modules over operads.

In positive characteristic,
one might consider divided power algebras,
which do not have the form of a functor $\Sym(M,E)$ associated to a $\Sigma_*$-module $M$.
In characteristic $2$,
we also have the exterior algebra functor $\Lambda(E) = \bigoplus_{r=0}^{\infty} E^{\otimes r}/(x^2\equiv 0)$
which occurs in the standard definition of the Chevalley-Eilenberg homology
and is not a functor of the form $\Sym(M,E)$ too.
In the sequel,
we do not use such constructions.
In particular,
when we deal with the Chevalley-Eilenberg complex,
we tacitely consider the most usual definition involving a symmetric algebra $\Sym(\Sigma G) = \bigoplus_{r=0}^{\infty} ((\Sigma G)^{\otimes r})_{\Sigma_r}$
on a suspension of the Lie algebra $G$,
and not an exterior algebra.

In fact,
we only apply the Chevalley-Eilenberg complex to Lie algebras
belonging to the category of connected $\Sigma_*$-modules (the definition of a connected $\Sigma_*$-module
is recalled next, in~\S\ref{Background:Kunneth})
and we deduce from observations of~\cite[\S 1.2]{FressePartitions}
that symmetric algebras $\Sym(E)$
like all functors of the form $\Sym(M,E)$
behave well in homology when the variable $E$ ranges over connected $\Sigma_*$-modules (see again~\S\ref{Background:Kunneth})
-- even when the ground ring is not a field of characteristic zero.
Therefore the generalized symmetric algebra functors associated to $\Sigma_*$-modules $\Sym(M)$
and the functors associated to modules over operads $\Sym_{\ROp}(M)$
are sufficient for our purpose.

\subsubsection{Model categories}\label{Background:ModelCategories}
We use the theory of model categories to give an abstract characterization of iterated bar constructions.

Recall that the category of dg-modules $\C$ has a standard model structure
in which a morphism is a weak-equivalence if it induces an isomorphism in homology,
a fibration if it is degreewise surjective,
a cofibration if it has the right lifting property with respect to acyclic fibrations (see for instance~\cite[\S 2.3]{Hovey}).
The category of $\Sigma_*$-modules $\M$
inherits a model structure
so that a morphism $f: M\rightarrow N$
forms a weak-equivalence, respectively a fibration, if each of its components $f: M(n)\rightarrow N(n)$
defines a weak-equivalence, respectively a fibration, in the category of dg-modules.
The category of right modules over an operad $\M{}_{\ROp}$
has a similarly defined model structure (see~\cite[\S 14]{Bar0}).
In all cases,
the cofibrations are characterized by the right lifting property with respect to acyclic fibrations,
where an acyclic fibration refers to a morphism which is both a fibration and a weak-equivalence.
By convention,
we say that a morphism of $\Sigma_*$-modules (respectively, of right $\ROp$-modules) $f: M\rightarrow N$
forms a $\C$-cofibration if its components $f: M(n)\rightarrow N(n)$
are cofibrations of dg-modules.
Any cofibration of $\Sigma_*$-modules (respectively, of right $\ROp$-modules)
forms a $\C$-cofibration (under the assumption that the operad $\ROp$ is $\C$-cofibrant in the case of right $\ROp$-modules),
but the converse implication does not hold.

The model categories $\E = \C,\M,\M{}_{\ROp}$ are all cofibrantly generated
and monoidal in the sense of~\cite[\S 11.3.3]{Bar0}
(the unit object is cofibrant and the tensor product satisfies the pushout-product axiom).
These assertions imply that the category of $\POp$-algebras in $\C$,
the category of $\POp$-algebras in $\M$ (equivalently, the category of left $\POp$-modules),
and the category of $\POp$-algebras in $\M{}_{\ROp}$ (equivalently, the category of $(\POp,\ROp)$-bimodules),
inherit a natural semi-model structure
when the operad $\POp$ is cofibrant as a $\Sigma_*$-module (see~\cite[\S 12]{Bar0}).
In all cases $\E = \C,\M,\M{}_{\ROp}$,
the forgetful functor $U: {}_{\POp}\E\rightarrow\E$
preserves cofibrations with a cofibrant domain,
but a morphism of $\POp$-algebras which forms a cofibration in the underlying category $\E$
does not form a cofibration in the category of $\POp$-algebras in general.
By convention,
we say that a morphism of $\POp$-algebras in $\E$
forms an $\E$-cofibration if it defines a cofibration in the underlying category $\E$,
a $\POp$-algebra $A\in{}_{\POp}\E$
is $\E$-cofibrant if the unit morphism $\eta: \POp(0)\rightarrow A$
forms an $\E$-cofibration.

In the next sections,
we only deal with operads such that $\POp(0) = 0$.
In this context,
a $\POp$-algebra $A$ is $\E$-cofibrant if and only if $A$ forms a cofibrant object of the underlying category $\E$.
In main theorems,
we prefer to use the latter formulation.

The model categories of $\POp$-algebras are used extensively in~\cite{Bar1},
but in this paper,
we essentially use the notion of an $\E$-cofibration
and the model structures of underlying categories of algebras over operads.

\subsubsection{Model structures and operads}\label{Background:OperadModelStructure}
The category of operads $\Op$
comes also equipped with a semi-model structure
so that a morphism $f: \POp\rightarrow\QOp$ forms a weak-equivalence, respectively a fibration,
if it defines a weak-equivalence, respectively a fibration, in the category of $\Sigma_*$-modules
(see~\cite[\S 13]{Bar0} and further bibliographical references therein).
Again,
the cofibrations of the category of operads are characterized by the right lifting property
with respect to acyclic fibration.
The forgetful functor $U: \Op\rightarrow\M$
preserves cofibrations with a cofibrant domain,
but a morphism of operads which forms a cofibration in the category of $\Sigma_*$-modules
does not form a cofibration in the category of operads in general.

By convention,
a morphism of operads which defines a cofibration in the category of $\Sigma_*$-modules
is called a $\Sigma_*$-cofibration,
an operad $\POp$ is $\Sigma_*$-cofibrant
if the unit morphism $\eta: I\rightarrow\POp$
is a $\Sigma_*$-cofibration.

Similarly,
we say that a morphism of operads $\phi: \POp\rightarrow\QOp$
is a $\C$-cofibration
if its components $f: \POp(n)\rightarrow\QOp(n)$
define cofibrations in the underlying category of dg-modules,
we say that an operad $\POp$
is $\C$-cofibrant if the unit morphism $\eta: I\rightarrow\POp$
is a $\C$-cofibration.
In the sequel,
we assume tacitely that any operad $\POp$
is at least $\C$-cofibrant
whenever we deal with model structures.
This assumption is necessary in most homotopical constructions.

In~\cite{Bar1},
we need the notion of a cofibrant operad
to define the multiplicative structure of the bar construction.
For this reason,
we still have to deal with cofibrant operads
in~\S\ref{IteratedBarModules},
where we study the iterated bar complex deduced from the construction of~\cite{Bar1}.
Otherwise,
we only deal with the model structure of $\Sigma_*$-modules.

\subsubsection{The K\"unneth isomorphism}\label{Background:Kunneth}
The composite of $\Sigma_*$-modules
has an expansion of the form
\begin{equation*}
M\circ N = \Sym(M,N) = \bigoplus_{r = 0}^{\infty} (M(r)\otimes N^{\otimes r})_{\Sigma_r}
\end{equation*}
and we have a natural K\"unneth morphism
\begin{equation*}
H_*(M)\circ H_*(N)\rightarrow H_*(M\circ N),
\end{equation*}
for all $\Sigma_*$-modules $M,N\in\M$.
In the paper,
we also use K\"unneth morphisms $\POp(H_*(N))\rightarrow H_*(\POp(N))$
associated to composites of the form $\POp(N) = \POp\circ N$,
where $\POp$ is a graded operad (equipped with a trivial differential).

In~\cite[\S\S 1.3.7-1.3.9]{FressePartitions},
we observe that the composite $M\circ N$
has a reduced expansion in which no coinvariant occurs
when $N(0) = 0$ (we say that the $\Sigma_*$-module $N$ is connected).
As a byproduct,
if $M,N$ are connected $\C$-cofibrant and the homology modules $H_*(M),H_*(N)$
consist of projective $\kk$-modules,
then the K\"unneth morphism
$H_*(M)\circ H_*(N)\rightarrow H_*(M\circ N)$
is an isomorphism.

In the sequel,
we apply most constructions to connected $\Sigma_*$-modules
for which such good properties hold.

\part*{The construction of iterated bar complexes}\label{IteratedComplexDefinition}

In this part,
we explain the construction of iterated bar complexes
\begin{equation*}
B^n: {}_{\ROp}\E\rightarrow\E
\end{equation*}
as functors $B^n(-) = \Sym_{\ROp}(B^n_{\ROp},-)$
associated to modules $B^n_{\ROp}$ over an operad $\ROp$.

First, in~\S\ref{IteratedBarModules},
we review the definition of the bar complex of $A_\infty$-algebras, of $E_\infty$-algebras,
and we prove that the construction of~\cite{Bar1}, in the context of algebras over an $E_\infty$-operad $\EOp$,
gives an iterated bar complex of the required form $B^n(-) = \Sym_{\EOp}(B^n_{\EOp},-)$.
In~\S\ref{QuasiFreeLifting},
we give a simple construction of the iterated bar module $B^n_{\EOp}$
as a lifting of quasi-free modules over operads.

Then
we aim at proving that the $n$-fold bar module $B^n_{\EOp}$
is an extension of a module $B^n_{\EOp_n}$
defined over a filtration layer $\EOp_n\subset\EOp$,
for certain $E_\infty$-operads $\EOp$
equipped with a filtration
\begin{equation}
\EOp_1\subset\dots\subset\EOp_n\subset\dots\subset\colim_n\EOp_n = \EOp
\end{equation}
such that $\EOp_n$ is an $E_n$-operad.
For this purpose,
we use that usual $E_\infty$-operads are equipped with a particular cell structure,
reviewed in~\S\ref{CompleteGraphOperad},
which refines their filtration by $E_n$-operads.
We observe that the iterated bar modules $B^n_{\COp}$ of the commutative operad $\COp$
are equipped with a cell structure of the same shape (\S\ref{KStructure})
and so are the iterated bar modules $B^n_{\EOp}$ associated to any $E_\infty$-operad $\EOp$ (\S\ref{EnDefinition}).
We draw our conclusion from the latter assertion.

The module $B^n_{\EOp_n}$
determines an extension of the iterated bar complex $B^n: {}_{\EOp}\E\rightarrow\E$
to the category of $\EOp_n$-algebras ${}_{\EOp_n}\E\supset{}_{\EOp}\E$.

\subsubsection*{Convention}
From now on, we assume tacitely that any operad $\POp$ satisfies $\POp(0) = 0$.

This assumption $\POp(0) = 0$ (we also say that the operad $\POp$ is non-unitary)
amounts to considering algebras over operads
without $0$-ary operation.
This setting simplifies the definition of iterated bar complexes and is required
in constructions and arguments of~\S\ref{EnHomology}.

Unital algebras are more usually considered in the literature,
but for the standard categories of algebras (associative algebras, commutative algebras, \dots)
we have an equivalence of categories between algebras without unit
and algebras with unit and augmentation.
In one direction, to any algebra $A$,
we can associate the algebra $A_+ = \kk 1\oplus A$,
where a unit is formally added.
This algebra inherits an augmentation $\epsilon: A_+\rightarrow\kk$
defined by the projection onto the summand $\kk 1$.
In the other direction,
we form the augmentation ideal $\bar{A} = \ker(A\rightarrow\kk)$
to associate an algebra without unit $\bar{A}$
to any augmented algebra $A$.
It is easy to check that these constructions
are inverse to each other.

Similar observations hold for coalgebras and Hopf algebras,
but in the context of Hopf algebras,
the distribution relation between products and coproducts becomes:
\begin{equation*}
\Delta(a\cdot b) = \Delta(a)\cdot\Delta(b) - \bigl\{a\otimes b + b\otimes a\bigr\}.
\end{equation*}
The standard distribution relation of Hopf algebras can be retrieved from this formula by adding terms $x\otimes 1 + 1\otimes x$
to each application of the diagonal $\Delta(x)$.

Throughout this article,
we deal with constructions of the literature which apply to augmented algebras.
Thus we just consider the equivalent construction
on the augmentation ideal of the algebra.
Actually certain constructions studied in the paper, like the bar complex, are naturally defined on the augmentation ideal,
not on the algebra itself,
and this observation gives the main motivation for the point of view which we adopt throughout the article.

\section{Iterated bar modules}\label{IteratedBarModules}
In this section,
we study the structure of the iterated bar complex~$B^n(A)$
as it arises from the construction of~\cite{Bar1}
for algebras over an $E_\infty$-operad:
we check that the iterated bar construction $B^n: {}_{\ROp}\E\rightarrow\E$,
where $\ROp$ is either the commutative operad $\ROp = \COp$
or an $E_\infty$-operad $\ROp = \EOp$,
is an instance of a functor associated to a right $\ROp$-module,
for which we adopt the notation $B^n_{\ROp}$;
we prove that the iterated bar module $B^n_{\EOp}$
associated to an $E_\infty$-operad $\EOp$
defines a cofibrant replacement, in the category of right $\EOp$-modules,
of the iterated bar module $B^n_{\COp}$
associated to the commutative operad $\COp$;
we use this result to obtain a simple characterization of the iterated bar complex $B^n: {}_{\EOp}\E\rightarrow\E$
for algebras over an $E_\infty$-operad $\EOp$.

First of all,
we review the definition of an $A_\infty$-operad, of an $E_\infty$-operad,
and we recall the definition of the bar complex $B: {}_{\ROp}\E\rightarrow\E$
for algebras over such operads.

\subsubsection{The associative and commutative operads}\label{IteratedBarModules:AssociativeCommutativeOperads}
Throughout the paper,
we use the notation $\AOp$ for the operad of associative algebras without unit
and the notation $\COp$ for the operad of associative and commutative algebras without unit (for short, commutative algebras).

The associative operad $\AOp$ is defined by:
\begin{equation*}
\AOp(r) = \begin{cases} 0, & \text{if $r = 0$}, \\
\kk[\Sigma_r], & \text{otherwise},
\end{cases}
\end{equation*}
where $\kk[\Sigma_r]$ represents the free $\kk$-module spanned by $\Sigma_r$.
The element of $\AOp$
represented by the identity permutation $\id\in\Sigma_2$
is also denoted by $\mu\in\AOp(2)$.
The operation $\mu: A^{\otimes 2}\rightarrow A$
determined by $\mu\in\AOp(2)$
represents the product of associative algebras.

The commutative operad $\COp$ is defined by:
\begin{equation*}
\COp(r) = \begin{cases} 0, & \text{if $r = 0$}, \\
\kk, & \text{otherwise},
\end{cases}
\end{equation*}
where the one dimensional $\kk$-module $\kk$ is equipped with a trivial $\Sigma_r$-action.
The generating element $\mu\in\COp(2)$
determines an operation $\mu: A^{\otimes 2}\rightarrow A$
which represents the product of commutative algebras.
For a non-empty finite set $\eset = \{e_1,\dots,e_r\}$,
we also have $\COp(\eset) = \kk$.
The expression of the commutative word $e_1\cdots e_r$
can conveniently be used to denote the generator of~$\COp(\eset)$
when necessary.

The augmentations $\alpha: \kk[\Sigma_r]\rightarrow\kk$, $r\in\NN$,
define an operad morphism $\alpha: \AOp\rightarrow\COp$
such that the restriction functor $\alpha^*: {}_{\COp}\E\rightarrow{}_{\AOp}\E$
represents the usual embedding from the category of commutative algebras
into the category of associative algebras.

Recall that an $A_\infty$-operad (respectively, of an $E_\infty$-operad)
refers to an operad equivalent to the associative operad $\AOp$ (respectively, to the commutative operad $\COp$)
in the homotopy category of operads.
In the next paragraphs,
we give a more precise definition of these structures
in a form suitable for our needs.

\subsubsection{On $A_\infty$-operads}\label{IteratedBarModules:AinfinityOperads}
For us,
an $A_\infty$-operad consists of an operad $\KOp$
augmented over the associative operad $\AOp$
so that the augmentation $\epsilon: \KOp\rightarrow\AOp$
defines an acyclic fibration in the category of operads.
The restriction functor $\epsilon^*: {}_{\AOp}\E\rightarrow{}_{\KOp}\E$,
naturally associated to the augmentation of an $A_\infty$-operad,
defines an embedding of categories from the category of associative algebras
to the category of $\KOp$-algebras.

In our work on the bar complex,
we only use a particular $A_\infty$-operad (outside the associative operad),
namely the chain operad of Stasheff's associahedra (Stasheff's operad for short)
and
the letter $\KOp$ will only refer to this $A_\infty$-operad.
The Stasheff operad is quasi-free and can be defined by a pair $\KOp = (\FOp(M),\partial)$,
where $\FOp(M)$ is a free operad and $\partial: \FOp(M)\rightarrow\FOp(M)$
is an operad derivation
that determines the differential of~$\KOp$.
The generating $\Sigma_*$-module $M$
of Stasheff's operad $\KOp$
is defined by:
\begin{equation*}
M(r) = \begin{cases} 0, & \text{if $r = 0,1$}, \\
\Sigma_r\otimes\kk\,\mu_r, & \text{otherwise},
\end{cases}
\end{equation*}
where $\mu_r$ is homogeneous of degree $r-2$.
The derivation $\partial: \Free(M)\rightarrow\Free(M)$ is determined on generating operations
by the formula
\begin{equation*}
\partial(\mu_r) = \sum_{s+t-1 = r}\Bigl\{\sum_{i=1}^{s} \pm\mu_s\circ_i\mu_t\Bigr\},
\end{equation*}
for some sign $\pm$ (see~\cite{MarklMinimal}).
This operad is cofibrant as an operad.

The augmentation $\epsilon: \KOp\rightarrow\AOp$
cancels the generating operations $\mu_r\in\KOp(r)$ such that $r>2$
and maps $\mu_2\in\KOp(2)$ to the generating operation of the associative operad $\mu\in\AOp(2)$,
the operation
which represents the product $\mu: A^{\otimes 2}\rightarrow A$
in the category of associative algebras.

The structure of an algebra over Stasheff's operad $\KOp$
is equivalent to a collection of homomorphisms $\mu_r: A^{\otimes r}\rightarrow A$, $r\in\NN$,
which give the action of the generating operations $\mu_r\in\KOp(r)$ on $A$,
so that we have the differential relation
\begin{equation*}
\delta(\mu_r)(a_1,\dots,a_r) = \sum_{s+t-1 = r}\Bigl\{\sum_{i=1}^{s}\pm\mu_s(a_1,\dots,\mu_t(a_i,\dots,a_{i+t-1}),\dots,a_r)\Bigr\},
\end{equation*}
where $\delta(\mu_r)$ refers to the differential of homomorphisms.
In this way,
we retrieve the usual definition of the structure of an $A_\infty$-algebra.

An associative algebra
is equivalent to a $\KOp$-algebra
such that $\mu_r(a_1,\dots,a_r) = 0$, for $r>2$.

\subsubsection{On $E_\infty$-operads}\label{IteratedBarModules:EinfinityOperads}
An $E_\infty$-operad is an operad $\EOp$
augmented over the commutative operad $\COp$
so that the augmentation $\epsilon: \EOp\rightarrow\COp$
defines an acyclic fibration in the category of operads.
An $E_\infty$-operad is usually assumed to be $\Sigma_*$-cofibrant.
This requirement ensures that the category of $\EOp$-algebras
is equipped with a semi-model structure.
Moreover,
different $\Sigma_*$-cofibrant $E_\infty$-algebras
have equivalent homotopy categories of algebras.

The restriction functor $\epsilon^*: {}_{\COp}\E\rightarrow{}_{\EOp}\E$,
associated to the augmentation of an $E_\infty$-operad,
defines an embedding of categories
from the category of commutative algebras to the category of $\EOp$-algebras
(as in the context of $A_\infty$-operads).
For any $E_\infty$-operad $\EOp$,
we can pick a lifting in the diagram
\begin{equation*}
\xymatrix{ \KOp\ar@{.>}[r]^{\eta}\ar@{->>}[d]^{\sim} & \EOp\ar@{->>}[d]^{\sim} \\
\AOp\ar[r]_{\alpha} & \COp }
\end{equation*}
to make $\EOp$ an object of the comma category $\Op\backslash\KOp$
of operads under $\KOp$.
The commutativity of the diagram
implies that the restriction functor $\eta^*: {}_{\EOp}\E\rightarrow{}_{\KOp}\E$
associated to the morphism $\eta: \KOp\rightarrow\EOp$
fits a commutative diagram of functors
\begin{equation*}
\xymatrix{ {}_{\KOp}\E & {}_{\EOp}\E\ar@{.>}[l]_{\eta^*} \\
{}_{\AOp}\E\ar@{^{(}->}[]!U+<0pt,2pt>;[u] &
{}_{\COp}\E\ar@{^{(}->}[]!U+<0pt,2pt>;[u]\ar[l]^{\alpha^*} }
\end{equation*}
and hence extends the usual embedding
from the category of commutative algebras to the category of associative algebras.

\subsubsection{The example of the Barratt-Eccles operad}\label{IteratedBarModules:BarrattEccles}
The Barratt-Eccles operad is a classical instance of an $E_\infty$-operad
introduced by M. Barratt and P. Eccles
in the simplicial setting~\cite{BarrattEccles}.
Throughout the paper,
we use this nice combinatorial operad to illustrate our constructions.
For our purpose,
we consider a dg-version of the Barratt-Eccles operad
which is defined by the normalized chain complexes
\begin{equation*}
\EOp(r) = N_*(E\Sigma_r),
\end{equation*}
where $E\Sigma_r$ denotes the acyclic homogeneous bar construction of the symmetric group $\Sigma_r$.
By convention,
we remove the $0$-component $N_*(E\Sigma_0)$ to assume $\EOp(0) = 0$
and to have a non-unitary analogue of this operad.
In~\cite{BergerFresse},
we observe that the Barratt-Eccles acts on the cochain complex of spaces.
Therefore this instance of $E_\infty$-operad is also well suited when we tackle topological applications
of iterated bar constructions.
In~\S\ref{CompleteGraphOperad},
we recall that the Barratt-Eccles is also equipped with a filtration of the form (*)
which refines into a particularly nice cell structure.
Therefore
the Barratt-Eccles is also particularly well suited for the connection, studied in the present article,
between iterated bar complexes and the homology of $E_n$-algebras.

The dg-module $\EOp(r)$ is spanned in degree $d$
by the $d$-simplices of permutations
\begin{equation*}
(w_0,\dots,w_d)\in\Sigma_r\times\dots\times\Sigma_r
\end{equation*}
divided out by the submodule spanned by degenerate simplices
\begin{equation*}
s_j(w_0,\dots,w_{d-1}) = (w_0,\dots,w_j,w_j,\dots,w_{d-1}),\quad\text{$j = 0,\dots,d-1$}.
\end{equation*}
The differential of a simplex $(w_0,\dots,w_d)\in\EOp(r)$
is defined by the alternate sum
\begin{equation*}
\delta(w_0,\dots,w_d) = \sum_{i=0}^{d} \pm(w_0,\dots,\widehat{w_i},\dots,w_d).
\end{equation*}
The symmetric group $\Sigma_r$ operates diagonally on $\EOp(r)$
and the composition products of $\EOp$
are yielded by an explicit substitution process on permutations.
The reader is referred to~\cite{BergerFresse}
for a detailed definition
and a comprehensive study of the Barratt-Eccles operad
in the dg-setting.

If we extend the definition of the dg-modules $\EOp(r)$ to finite sets $\eset = \{e_1,\dots,e_r\}$ as explained in~\S\ref{Background:SetIndexing},
then we obtain a dg-module $\EOp(\eset)$ spanned by non-degenerate simplices $(w_0,\dots,w_d)$
such that $w_i$ is a bijection $w_i: \{1,\dots,r\}\xrightarrow{\sim}\{e_1,\dots,e_r\}$.
In applications,
we use that such a bijection $w: \{1,\dots,r\}\rightarrow\{e_1,\dots,e_r\}$
is equivalent to an ordering of~$\eset$
and we represent this ordering by the sequence of values $(w(1),\dots,w(r))$
of the bijection~$w$.

The operad equivalence $\epsilon: \EOp\xrightarrow{\sim}\COp$
giving an $E_\infty$-operad structure to $\EOp$
is defined by the usual augmentations $\epsilon: N_*(E\Sigma_r)\xrightarrow{\sim}\kk$.
In the sequel,
we use the standard section of this augmentation $\iota: \COp(r)\rightarrow\EOp(r)$,
which identifies $\COp(r) = \kk$ with the component of $\EOp(r)$
spanned by the identity permutation $\id$ in degree $0$,
and the usual contracting chain homotopy $\nu: \EOp(r)\rightarrow\EOp(r)$
defined by $\nu(w_0,\dots,w_d) = (w_0,\dots,w_d,\id)$
for any simplex of permutations $(w_0,\dots,w_d)\in\EOp(r)$.
Naturally,
these definitions have a natural generalization giving a section $\iota: \COp(\eset)\rightarrow\EOp(\eset)$
and a contracting chain homotopy $\nu: \EOp(\eset)\rightarrow\EOp(\eset)$
for any finite set $\eset$
as long as we fix an ordering $\eset = \{e_1<\dots<e_r\}$.
In this setting,
the augmentation $\epsilon: \EOp(\eset)\rightarrow\COp(\eset)$
simply forgets the ordering of bijections $u: \{1,\dots,r\}\rightarrow\eset$
in degree $0$,
the section $\iota: \COp(\eset)\rightarrow\EOp(\eset)$ is the map
sending a commutative word $e_1\cdots e_r\in\COp(\eset)$
to the sequence $\iota(e_1\cdots e_r) = e_1\cdots e_r$
defined by the ordering of~$\eset$,
and the chain-homotopy $\nu: \EOp(\eset)\rightarrow\EOp(\eset)$
simply inserts this ordered sequence $e_1\cdots e_r$
at the final position of simplices.

The associative operad $\AOp$ can be identified with the degree $0$ part of~$\EOp$
and forms a suboperad of the Barratt-Eccles operad
such that the restriction of the augmentation $\epsilon: \EOp\xrightarrow{\sim}\COp$
to $\AOp\subset\EOp$
agrees with the morphism $\alpha:\AOp\rightarrow\COp$
defined in~\S\ref{IteratedBarModules:AssociativeCommutativeOperads}.
Thus we have a factorization
\begin{equation*}
\xymatrix{ & \EOp\ar@{.>}[dr]^{\sim} & \\ \AOp\ar@{.>}[ur]\ar[rr]_{\alpha} && \COp }
\end{equation*}
and the lifting construction of~\S\ref{IteratedBarModules:EinfinityOperads}
is not necessary for the Barratt-Eccles operad.

\subsubsection{The bar complex of algebras over Stasheff's operad}\label{IteratedBarModules:BarComplex}
The bar complex is defined naturally for $\KOp$-algebras.

Let $A$ be a $\KOp$-algebra in $\E$, where $\E = \C$ the category of dg-modules,
or $\E = \M$ the category of $\Sigma_*$-objects,
or $\E = \M{}_{\ROp}$ the category of right modules over an operad $\ROp$.
The (reduced and normalized) bar complex of~$A$
is defined by a pair $B(A) = (\Tens^c(\Sigma A),\partial)$
formed by the (non-augmented) tensor coalgebra
\begin{equation*}
\Tens^c(\Sigma A) = \bigoplus_{d=1}^{\infty} (\Sigma A)^{\otimes d},
\end{equation*}
where $\Sigma A$ is the suspension of $A$ in $\E$,
together with a homomorphism of degree $-1$
\begin{equation*}
\partial\in\Hom_{\E}(\Tens^c(\Sigma A),\Tens^c(\Sigma A)),
\end{equation*}
called abusively the bar differential,
defined pointwise by the formula
\begin{equation*}
\partial(a_1\otimes\dots\otimes a_d)
= \sum_{r=2}^{d}\Bigl\{\sum_{i=1}^{d-r+1} \pm a_1\otimes\dots
\otimes\mu_r(a_i,\dots,a_{i+r-1})\otimes\dots
\otimes a_d\Bigr\}.
\end{equation*}
The differential of the bar complex $B(A)$ is the sum $\delta+\partial$ of the natural differential of the tensor coalgebra
$\delta: \Tens^c(\Sigma A)\rightarrow\Tens^c(\Sigma A)$,
which is induced by the internal differential of~$A$,
together with the bar differential
$\partial: \Tens^c(\Sigma A)\rightarrow\Tens^c(\Sigma A)$,
which determined by the $\KOp$-action on $A$.
Note that the identity $(\delta+\partial)^2 = 0$ holds in~$B(A)$
for every algebra over Stasheff's operad.
The term ``differential'' is abusive for $\partial: \Tens^c(\Sigma A)\rightarrow\Tens^c(\Sigma A)$,
because the relation $(\delta+\partial)^2 = 0$
does not hold for the isolated term $\partial$.
Usually,
we call twisting homomorphisms the homomorphisms of dg-modules $\partial: C\rightarrow C$, like the bar differentials,
whose addition to the internal differential of~$C$
satisfies the equation of differentials $(\delta+\partial)^2 = 0$.
More recollections on this notion are given next (in~\S\ref{QuasiFreeLifting:TwistingCochains})
when we tackle the definition of an iterated bar complex by the direct construction
of a twisting homomorphism.

The definition of generalized point-tensors in~\S\ref{Background:TensorProduct}
allows us to give a sense to the formula of the bar differential
in the context of $\Sigma_*$-objects $\E = \M$ and right modules over an operad $\E = \M{}_{\ROp}$,
and not only in the context of dg-modules $\E = \C$.

\subsubsection{Coalgebra structures}\label{IteratedBarModules:CoalgebraStructure}
The tensor coalgebra $T^c(\Sigma A)$
is equipped with a diagonal $\Delta: T^c(\Sigma A)\rightarrow T^c(\Sigma A)\otimes T^c(\Sigma A)$
defined by the deconcatenation of tensors:
\begin{equation*}
\Delta(a_1\otimes\dots\otimes a_d) = \sum_{e=1}^{d} (a_1\otimes\dots\otimes a_e)\otimes(a_{e+1}\otimes\dots\otimes a_d).
\end{equation*}
This diagonal makes $T^c(\Sigma A)$ a coassociative coalgebra.
Recall that we tacitely consider non-augmented coalgebras only.
Therefore, in our definition of the bar complex,
we take a non-augmented version of the tensor coalgebra $T^c(\Sigma A)$,
where the component of tensors of order $0$
is removed.

The twisting homomorphism of the bar complex $\partial: T^c(\Sigma A)\rightarrow T^c(\Sigma A)$
satisfies the coderivation relation
$\Delta\partial = (\partial\otimes\id+\id\otimes\partial)\Delta$
with respect to the diagonal of the tensor coalgebra.
From this observation,
we conclude that the bar complex $B(A) = (T^c(\Sigma A),\partial)$
forms a coalgebra in the category of dg-modules $\E = \C$
(respectively, in the category of $\Sigma_*$-modules $\E = \M$,
in the category of right modules over an operad $\E = \M{}_{\ROp}$)
in which $B(A)$ is defined.

This structure is used in~\S\ref{EnHomology} to simplify calculations
of differentials in iterated bar complexes.

\subsubsection{The suspension morphism}\label{IteratedBarModules:Suspension}
The bar differential vanishes on the summand
\begin{equation*}
\Sigma A\subset\bigoplus_{d=1}^{\infty} (\Sigma A)^{\otimes d} = T^c(\Sigma A).
\end{equation*}
Hence,
the canonical embedding $\Sigma A\hookrightarrow T^c(\Sigma A)$
defines a natural morphism $\sigma_A: \Sigma A\rightarrow B(A)$.
This morphism is called the suspension morphism.
For short,
we can set $\sigma = \sigma_A$.

\subsubsection{The bar complex and restriction of algebra structures}\label{IteratedBarModules:BarComplexRestriction}
The bar complex defines a functor $B: {}_{\KOp}\E\rightarrow\E$
from the category of algebras over Stasheff's operad
to the underlying category.

Let $\ROp$ be an operad together with a fixed morphism $\eta: \KOp\rightarrow\ROp$,
so that $\ROp$ forms an object of the category of operads under $\KOp$.
The composite of the bar construction $B: {}_{\KOp}\E\rightarrow\E$
with the restriction functor $\eta^*: {}_{\ROp}\E\rightarrow{}_{\KOp}\E$
defines a bar construction on the category of $\ROp$-algebras.
To simplify,
we omit marking restriction functors in the notation of the bar complex $B(A)$
associated to an $\ROp$-algebra $A$.

The bar complex of an $\ROp$-algebra is given by the same construction
as in the context of a $\KOp$-algebra.
In the definition of the bar differential,
we simply take the image of the generating operations of Stasheff's operad $\mu_r\in\KOp(r)$
under the morhism $\eta: \KOp\rightarrow\ROp$
to determine the action of~$\mu_r$
on~$A$.

In the case of the associative operad $\AOp$,
the bar differential reduces to terms:
\begin{equation*}
\partial(a_1\otimes\dots\otimes a_d)
= \sum_{i=1}^{d-1} \pm a_1\otimes\dots\otimes\mu_2(a_i,a_{i+1})\otimes\dots\otimes a_d.
\end{equation*}
Thus the restriction of the functor $B: {}_{\KOp}\E\rightarrow\E$ to the category of associative algebras ${}_{\AOp}\E$
gives the usual bar complex of associative algebras.

In the case of an $E_\infty$-operad $\EOp$,
the bar construction gives a functor $B: {}_{\EOp}\E\rightarrow\E$
which extends the usual bar construction of commutative algebras
since restriction functors assemble to give a commutative diagram:
\begin{equation*}
\xymatrix{ \E & {}_{\KOp}\E\ar[l]_{B} & {}_{\EOp}\E\ar@{.>}[l]_{\eta^*} \\
& {}_{\AOp}\E\ar@{^{(}->}[]!U+<0pt,2pt>;[u] & {}_{\COp}\E\ar@{^{(}->}[]!U+<0pt,2pt>;[u]\ar[l]^{\alpha^*} }.
\end{equation*}

\subsubsection{The suspension morphism and indecomposables}\label{IteratedBarModules:Indecomposables}
The indecomposable quotient of an associative (respectively, commutative) algebra $A$
is defined by the cokernel
\begin{equation*}
A/A^2 = \coker(\mu: A\otimes A\rightarrow A),
\end{equation*}
where $\mu$ refers to the product of $A$.
In the sequel,
we also use the notation $\Indec A = A/A^2$
for this quotient.
The functor $\Indec: {}_{\ROp}\E\rightarrow\E$, for $\ROp = \AOp,\COp$,
is left adjoint to the obvious functor $\Ab: \E\rightarrow{}_{\ROp}\E$
which identifies any object $E\in\E$
with an $\ROp$-algebra in $\E$
equipped with a trivial algebra structure.

Suppose that $A$ has a trivial internal differential.
Then the suspension morphism $\sigma: \Sigma A\rightarrow B(A)$
maps $A^2 = \im(\mu: A\otimes A\rightarrow A)$
to boundaries of $B(A)$.
Hence,
we obtain that the morphism $\sigma_*: A\rightarrow H_*(B(A))$
induced by the suspension in homology
admits a factorization
\begin{equation*}
\xymatrix{ \Sigma A\ar[rr]^{\sigma_*}\ar[dr] && H_*(B(A))\\
& \Sigma\Indec A\ar@{.>}[ur]_{\overline{\sigma}_*} & }.
\end{equation*}

\subsubsection{The bar module}\label{IteratedBarModules:BarModule}
In the context of right modules over an operad $\ROp$,
the bar complex defines a functor $B: {}_{\KOp}\M{}_{\ROp}\rightarrow\M_{\ROp}$
from the category of $\KOp$-algebras in right $\ROp$-modules
to the category of right $\ROp$-modules.

Recall that an operad $\ROp$
forms an algebra over itself in the category of right modules over itself.
By restriction of structure on the left,
an operad under $\KOp$ forms an algebra over $\KOp$
in the category of right modules over itself.
We adopt the notation $B_{\ROp} = B(\ROp)$
for the bar complex of this $\KOp$-algebra.
The object $B_{\ROp}$ is the bar module associated to $\ROp$.

Note that a morphism $\psi: \ROp\rightarrow\SOp$
in the category of operads under $\KOp$
defines a morphism of $\KOp$-algebras in right $\ROp$-modules
and hence induce a natural morphism $\psi_{\sharp}: B_{\ROp}\rightarrow B_{\SOp}$
in the category of right $\ROp$-modules.
This morphism has an adjoint $\psi_{\flat}: B_{\ROp}\circ_{\ROp}\SOp\rightarrow B_{\SOp}$.

The next assertions are proved in~\cite{Bar1}:

\begin{prop}[{see~\cite[\S 1.4]{Bar1}}]\label{IteratedBarModules:BarConstructionRepresentation}
Let $\ROp$ be an operad under $\KOp$.
We have a natural isomorphism $B(A) = \Sym_{\ROp}(B_{\ROp},A)$,
for all $A\in{}_{\ROp}\E$.

The bar module $B_{\ROp}$ satisfies the relation $B_{\ROp}\simeq B_{\KOp}\circ_{\KOp}\ROp$.
More generally,
the natural morphism $\psi_{\flat}: B_{\ROp}\circ_{\ROp}\SOp\rightarrow B_{\SOp}$
associated to any morphism $\psi: \ROp\rightarrow\SOp$ in the category of operads under $\KOp$
forms an isomorphism.\qed
\end{prop}

Recall that the extension of right modules over operads represents the composition with restriction functors at the functor level.
The relation $B_{\ROp}\simeq B_{\KOp}\circ_{\KOp}\ROp$
implies that the diagram of functors
\begin{equation*}
\xymatrix{ {}_{\KOp}\E\ar[dr]_{\Sym_{\KOp}(B_{\KOp})} &&
{}_{\ROp}\E\ar[dl]^{\Sym_{\ROp}(B_{\ROp})}\ar[ll]_{\eta^*} \\
& \E & }
\end{equation*}
commutes up to a natural isomorphism.
Hence the isomorphism $B_{\ROp}\simeq B_{\KOp}\circ_{\KOp}\ROp$
reflects the definition of the bar construction on the category of $\ROp$-algebras
by the restriction of a functor $B: {}_{\KOp}\E\rightarrow\E$
on the category of $\KOp$-algebras.

\subsubsection{The representation of suspension morphisms}\label{IteratedBarModules:SuspensionRepresentation}
The suspension morphism $\sigma_A: \Sigma A\rightarrow B(A)$
can naturally be identified with a morphism of functors
\begin{equation*}
\Sym_{\ROp}(\sigma_{\ROp}): \underbrace{\Sym_{\ROp}(\Sigma\ROp)}_{= \Sigma\Id}\rightarrow\Sym_{\ROp}(B_{\ROp})
\end{equation*}
associated to a morphism of right $\ROp$-modules $\sigma = \sigma_{\ROp}: \Sigma\ROp\rightarrow B_{\ROp}$,
which is nothing but the suspension morphism of the operad $\ROp$
viewed as a $\KOp$-algebra in the category of right modules over itself.

The observation of~\S\ref{IteratedBarModules:Indecomposables}
implies that this suspension morphism $\sigma: \Sigma\ROp\rightarrow B_{\ROp}$
factors through $\Indec\ROp$
in the case $\ROp = \AOp,\COp$.
The operads $\ROp = \AOp,\COp$
are equipped with an augmentation $\epsilon: \ROp\rightarrow I$
which provides the composition unit $I$
with the structure of a right $\ROp$-module.
In both cases $\ROp = \AOp,\COp$,
we have an obvious isomorphism $\Indec\ROp\simeq I$
in the category of right $\ROp$-modules.
Consequently,
the observation of~\S\ref{IteratedBarModules:Indecomposables}
implies the existence of a factorization
\begin{equation*}
\xymatrix{ \Sigma\ROp\ar[rr]^{\sigma_*}\ar[dr] && H_*(B_{\ROp})\\
& \Sigma I\ar@{.>}[ur]_{\overline{\sigma}_*} & }
\end{equation*}
in the category of right $\ROp$-modules,
for $\ROp = \AOp,\COp$.

\subsubsection{The bar complex for algebras over the commutative operad}\label{IteratedBarModules:BarCommutativeAlgebras}
Recall that the bar complex of a commutative algebra is equipped with a commutative algebra structure
(see~\cite[Expos\'e 4]{Cartan} or~\cite[\S X.12]{MacLane}).
The product $\mu: B(A)\otimes B(A)\rightarrow B(A)$, called the shuffle product,
is defined on components of the tensor coalgebra
by sums of tensor permutations
\begin{equation*}
(\Sigma A)^{\otimes d}\otimes(\Sigma A)^{\otimes e}\xrightarrow{\sum_w w_*}(\Sigma A)^{\otimes d+e}
\end{equation*}
such that $w$ ranges over the set of $(d,e)$-shuffles in $\Sigma_{d+e}$ (see~\cite[\S X.12]{MacLane} for details).
The shuffle product is naturally associative and commutative,
but the derivation relation $\partial\mu(\alpha,\beta) = \mu(\partial\alpha,\beta) + \pm\mu(\alpha,\partial\beta)$
with respect to the bar differential $\partial: T^c(\Sigma A)\rightarrow T^c(\Sigma A)$
holds for commutative algebras only.
The definition of the shuffle product is clearly functorial.
Hence the restriction of the bar construction
to the category of commutative algebras defines a functor $B: {}_{\COp}\E\rightarrow{}_{\COp}\E$.

The definition of the shuffle product makes sense not only in the context of dg-modules $\E = \C$,
but also in the context of $\Sigma_*$-modules $\E = \M$,
and in the context of right modules over an operad $\E = \M{}_{\ROp}$.
In particular, since we consider the commutative operad $\COp$ as a commutative algebra in the category of right modules over itself,
we obtain that the bar module of the commutative operad $B_{\COp}$
forms a commutative algebra in right $\COp$-modules,
equivalently a bimodule over the commutative operad.

In~\cite[\S 2.1]{Bar1},
we observe that the identity $B(A) = \Sym_{\COp}(B_{\COp},A)$
holds in the category of functors $F: {}_{\COp}\E\rightarrow{}_{\COp}\E$.

\subsubsection{The Hopf algebra structure}\label{IteratedBarModules:HopfStructure}
In~\S\ref{IteratedBarModules:BarCommutativeAlgebras},
we recall that the shuffle product preserves the differential of the bar complex.
One proves further (see for instance~\cite[\S 0.2, \S 1.4]{Reutenauer})
that the shuffle product $\mu: T^c(\Sigma E)\otimes T^c(\Sigma E)\rightarrow T^c(\Sigma E)$
commutes with the diagonal of the tensor coalgebra $T^c(\Sigma E)$,
for every object $E\in\E$
in a symmetric monoidal category $\E$
(but, as explained in the introduction of this part,
we have to remove units from the standard commutation relation
between the diagonal and the product
since we deal with a non-augmented version of the tensor coalgebra).
As a byproduct,
the bar complex $B(A) = (T^c(\Sigma A),\partial)$
inherits a natural commutative Hopf algebra structure,
for every commutative algebra~$A$.
Again,
this assertion holds in the context of dg-modules $\E = \C$,
in the context of $\Sigma_*$-modules $\E = \M$,
and in the context of right modules over an operad $\E = \M{}_{\ROp}$.

The Hopf algebra structure is used in~\S\ref{UsualBarComplexes}
to determine the homology of $B(A)$
for generalizations of usual commutative algebras.

\medskip
The next statement is proved in~\cite[\S 2]{Bar1}
in order to extend the multiplicative structure of the bar construction
to the category of algebras over any $E_\infty$-operad:

\begin{thm}[{see~\cite[Theorem 2.A]{Bar1}}]\label{IteratedBarModules:BarMultiplicativeStructure}
Let $\EOp$ be an $E_\infty$-operad.
Suppose $\EOp$ is cofibrant as an operad.
The bar module $B_{\EOp}$
can be equipped with the structure of an $\EOp$ algebra in right $\EOp$-modules
so that the natural isomorphism of right $\COp$-modules
$\epsilon_{\flat}: B_{\EOp}\circ_{\EOp}\COp\xrightarrow{\simeq} B_{\COp}$
defines an isomorphism of $\EOp$-algebras in right $\COp$-modules,
where we use a restriction of structure on the left to make $B_{\COp}$
an $\EOp$-algebra in right $\COp$-modules.\qed
\end{thm}

This theorem implies that the functor $B(A) = \Sym_{\EOp}(B_{\EOp},A)$
lands in the category of $\EOp$-algebras.
The relation $B_{\EOp}\circ_{\EOp}\COp\simeq B_{\COp}$ in the category ${}_{\EOp}\M{}_{\COp}$
implies that the diagram of functors
\begin{equation*}
\xymatrix{ {}_{\EOp}\E\ar[rr]^{B(-) = \Sym_{\EOp}(B_{\EOp})} && {}_{\EOp}\E \\
{}_{\COp}\E\ar[rr]_{B(-) = \Sym_{\COp}(B_{\COp})}\ar@{^{(}->}[]!U+<0pt,2pt>;[u] &&
{}_{\COp}\E\ar@{^{(}->}[]!U+<0pt,2pt>;[u] }
\end{equation*}
commutes up to a natural isomorphism in the category of $\EOp$-algebras,
because extensions on the right at the module level correspond to restrictions on the source at the functor level,
restrictions on the left at the module level correspond to restrictions on the target at the functor level.

Thus theorem~\ref{IteratedBarModules:BarMultiplicativeStructure}
implies that the bar construction of commutative algebras extends to a functor from $\EOp$-algebras to $\EOp$-algebras.
As a byproduct,
we have a well defined iterated bar complex $B^n: {}_{\ROp}\E\rightarrow{}_{\ROp}\E$
defined by the $n$-fold composite of the bar construction $B: {}_{\ROp}\E\rightarrow{}_{\ROp}\E$,
for $\ROp = \COp$ and $\ROp = \EOp$,
so that the diagram
\begin{equation*}
\xymatrix{ {}_{\EOp}\E\ar[r]^{B^n} & {}_{\EOp}\E \\
{}_{\COp}\E\ar[r]_{B^n}\ar@{^{(}->}[]!U+<0pt,2pt>;[u] &
{}_{\COp}\E\ar@{^{(}->}[]!U+<0pt,2pt>;[u] }
\end{equation*}
commutes.
We use the theory of modules over operads to determine the structure of this iterated bar construction.
Observe first:

\begin{prop}\label{IteratedBarModules:IteratedBarConstructionRepresentation}
Let $\ROp = \COp$ or $\ROp = \EOp$.

The $n$-fold bar complex $B^n: {}_{\ROp}\E\rightarrow{}_{\ROp}\E$,
defined by the $n$-fold composite
\begin{equation*}
{}_{\ROp}\E\xrightarrow{B}{}_{\ROp}\E\xrightarrow{B}\cdots\xrightarrow{B}{}_{\ROp}\E,
\end{equation*}
is isomorphic to the functor $\Sym(B^n_{\ROp}): {}_{\ROp}\E\rightarrow{}_{\ROp}\E$
associated to the composite module
\begin{equation*}
B^n_{\ROp} = B_{\ROp}\circ_{\ROp}\cdots\circ_{\ROp} B_{\ROp}.
\end{equation*}
\end{prop}

\begin{proof}
Immediate consequence of~\cite[Proposition 9.2.5]{Bar0} (see recollections in~\S\ref{Background:OperadFunctors}).
\end{proof}

We have further:

\begin{prop}\label{IteratedBarModules:CofibrantModule}
The iterated bar module $B^n_{\ROp}$ is cofibrant as a right $\ROp$-module, for $\ROp = \COp$ and $\ROp = \EOp$.
Moreover we have the relation $B^n_{\EOp}\circ_{\EOp}\COp\simeq B^n_{\COp}$.
\end{prop}

\begin{proof}
The bar module $B_{\ROp}$ is cofibrant as a right $\ROp$-module
by~\cite[Proposition 1.4.6]{Bar1}.
The functor $\Sym_{\ROp}(M): {}_{\ROp}\E\rightarrow\E$
associated to a cofibrant right $\ROp$-module $M$
maps the $\ROp$-algebras which are cofibrant in $\E$ to cofibrant objects of $\E$
by~\cite[Lemma 15.1.1]{Bar0}.
In the case $M = B_{\ROp}$ and $\E = \M_{\ROp}$,
we obtain that the relative composition product $B_{\ROp}\circ_{\ROp} N$
is cofibrant as a right $\ROp$-module if $N$ is so
(recall that $\Sym_{\ROp}(M,N) = M\circ_{\ROp} N$ for $\E = \M{}_{\ROp}$).
By induction,
we conclude that $B^n_{\ROp}$ forms a cofibrant right $\ROp$-module
as asserted in the proposition.

For a commutative algebra in right $\COp$-modules $N$,
we have relations
\begin{equation*}
B_{\EOp}\circ_{\EOp}(\epsilon^* N)\simeq(B_{\EOp}\circ_{\EOp}\COp)\circ_{\COp} N\simeq\epsilon^*(B_{\COp}\circ_{\COp} N)
\end{equation*}
(we prefer to mark the restriction of structure $\epsilon^*: {}_{\EOp}\M{}_{\COp}\rightarrow{}_{\EOp}\M{}_{\COp}$ in this formula).
The relation $B^n_{\EOp}\circ_{\EOp}\COp\simeq B^n_{\COp}$
follows by induction.
\end{proof}

\subsubsection{Iterated suspension morphisms}\label{IteratedBarModules:IteratedSuspension}
The definition of the suspension morphism $\sigma_A: \Sigma A\rightarrow B(A)$
can be applied to iterated bar complexes.
In this way,
we obtain a natural transformation
\begin{equation*}
\sigma_A: \Sigma B^{n-1}(A)\rightarrow B^n(A),
\end{equation*}
for every $n\geq 1$,
and for any $\ROp$-algebra $A$,
where $\ROp = \COp,\EOp$.
(By convention, in the case $n = 1$, we set $B^0 = \Id$.)

The application of this construction to the operad itself $\ROp = \COp,\EOp$
gives a morphism of right $\ROp$-modules
\begin{equation*}
\sigma_{\ROp}: \Sigma B^{n-1}_{\ROp}\rightarrow B^n_{\ROp}
\end{equation*}
such that $\Sym_{\ROp}(\sigma_{\ROp},A) = \sigma_A$, for every $A\in{}_{\ROp}\E$.
The diagram
\begin{equation*}
\xymatrix{ \Sigma B^{n-1}_{\EOp}\ar[d]\ar[r]^{\sigma_{\EOp}} & B^n_{\EOp}\ar[d] \\
\Sigma B^{n-1}_{\COp}\ar[r]_{\sigma_{\COp}} & B^n_{\COp} }
\end{equation*}
commutes by functoriality of the suspension.
We have equivalently $\sigma_{\EOp}\circ_{\EOp}\COp = \sigma_{\COp}$.

\subsubsection{The iterated bar module associated to non-cofibrant $E_\infty$-operads}\label{IteratedBarModules:BarNonCofibrantEinfinityOperad}
We can extend the definition of the iterated bar module $B^n_{\EOp}$
to any $E_\infty$-operad,
not necessarily cofibrant as an operad.
Indeed,
an $E_\infty$-operad $\EOp$ has a cofibrant replacement $\phi: \QOp\xrightarrow{\sim}\EOp$
which forms a cofibrant $E_\infty$-operad
and hence has an associated iterated bar module $B^n_{\QOp}$.
Define the iterated bar module of $\EOp$
by the extension of structure $B^n_{\EOp} = B^n_{\QOp}\circ_{\QOp}\EOp$.
By transitivity of relative composition products,
we have still $B^n_{\EOp}\circ_{\EOp}\COp\simeq B^n_{\COp}$.
Moreover,
the object $B^n_{\EOp}$ is cofibrant in the category of right $\EOp$-modules
since the functor of extension of structure $-\circ_{\QOp}\EOp$
is the left adjoint of a Quillen equivalence (see~\cite[Theorem 16.B]{Bar0}).

The functor $\Sym_{\EOp}(B^n_{\EOp}) = \Sym_{\EOp}(B^n_{\QOp}\circ_{\QOp}\EOp)$
is identified with the composite
\begin{equation*}
{}_{\EOp}\E\xrightarrow{\phi^*}{}_{\QOp}\E\xrightarrow{B^n}{}_{\QOp}\E,
\end{equation*}
where $\phi^*: {}_{\EOp}\E\rightarrow{}_{\QOp}\E$
is the restriction functor associated to $\phi: \QOp\xrightarrow{\sim}\EOp$.

\medskip
We deduce a simple homotopical characterization of the iterated bar module $B^n_{\EOp}$
from proposition~\ref{IteratedBarModules:CofibrantModule}.
Recall that:

\begin{fact}[{see~\cite[Theorem 16.B]{Bar0}}]\label{IteratedBarModules:ExtensionRestrictionEquivalence}
Let $\EOp$ be any $E_\infty$-operad, together with an augmentation $\epsilon: \EOp\xrightarrow{\sim}\COp$.
The extension and restriction functors
\begin{equation*}
\epsilon_!: \M_{\EOp}\rightleftarrows\M_{\COp} :\epsilon^*
\end{equation*}
define Quillen adjoint equivalences of model categories.
\end{fact}

Consequently:

\begin{prop}\label{IteratedBarModules:IteratedBarCofibrantReplacement}
Let $M^n_{\EOp}$ be any cofibrant right $\EOp$-module
equipped with a weak-equivalence
\begin{equation*}
f_{\flat}: M^n_{\EOp}\circ_{\EOp}\COp\xrightarrow{\sim} B^n_{\COp}
\end{equation*}
in the category of right $\COp$-modules.
The morphism of right $\EOp$-modules
\begin{equation*}
f_{\sharp}: M^n_{\EOp}\rightarrow B^n_{\COp}
\end{equation*}
adjoint to $f_{\flat}$
defines a weak-equivalence in the category of right $\EOp$-modules.\qed
\end{prop}

This proposition applies to the iterated bar module $B^n_{\EOp}$
according to the assertion of proposition~\ref{IteratedBarModules:CofibrantModule}.
As a byproduct,
every module $M^n_{\EOp}$
which satisfies the requirement of proposition~\ref{IteratedBarModules:IteratedBarCofibrantReplacement}
is connected to the iterated bar module $B^n_{\EOp}$
by weak-equivalences
\begin{equation*}
M^n_{\EOp}\xrightarrow{\sim} B^n_{\COp}\xleftarrow{\sim} B^n_{\EOp}
\end{equation*}
in the category of right $\EOp$-modules.
Thus the homotopy type of the iterated bar module $B^n_{\EOp}$
is fully characterized by the result of proposition~\ref{IteratedBarModules:CofibrantModule}
and proposition~\ref{IteratedBarModules:IteratedBarCofibrantReplacement}.

In this chain of weak-equivalences $M^n_{\EOp}\xrightarrow{\sim} B^n_{\COp}\xleftarrow{\sim} B^n_{\EOp}$
the right $\EOp$-module $M^n_{\EOp}$ is cofibrant by assumption,
and the iterated bar module $B^n_{\EOp}$
is cofibrant as well by proposition~\ref{IteratedBarModules:CofibrantModule},
but the right $\COp$-module $B^n_{\COp}$ is not cofibrant as a right $\EOp$-module.
Nevertheless,
as usual in a model category,
we can replace our chain of weak-equivalences
by a chain of weak-equivalences
\begin{equation*}
M^n_{\EOp}\xleftarrow{\sim}\,\cdot\,\xrightarrow{\sim}\,\cdots\,\xrightarrow{\sim} B^n_{\EOp}
\end{equation*}
so that all intermediate objects are cofibrant right $\EOp$-modules.

By~\cite[Theorem 15.1.A]{Bar0},
the existence of such weak-equivalences at the module level
implies:

\begin{thm}\label{IteratedBarModules:IteratedBarEquivalence}
Suppose we have a cofibrant right $\EOp$-module $M^n_{\EOp}$, where $\EOp$ is any $E_\infty$-operad,
together with a weak-equivalence
\begin{equation*}
f_{\flat}: M^n_{\EOp}\circ_{\EOp}\COp\xrightarrow{\sim} B^n_{\COp}
\end{equation*}
in the category of right $\COp$-modules.

The functor $\Sym_{\EOp}(M^n_{\EOp}): {}_{\EOp}\E\rightarrow\E$
determined by $M^n_{\EOp}$
is connected to the iterated bar complex $B^n: {}_{\EOp}\E\rightarrow\E$
by a chain of natural morphisms
\begin{equation*}
\Sym_{\EOp}(M^n_{\EOp},A)\xleftarrow{\sim}\,\cdot\,\xrightarrow{\sim}\,\cdots\,\xrightarrow{\sim}\Sym_{\EOp}(B^n_{\EOp},A) = B^n(A),
\end{equation*}
which are weak-equivalences as long as the $\EOp$-algebra $A$
is cofibrant in the underlying category~$\E$.
\end{thm}

\begin{proof}
By~\cite[Theorem 15.1.A]{Bar0},
a weak-equivalence $f: M\xrightarrow{\sim} N$, where $M,N$ are cofibrant right $\EOp$-modules,
induces a weak-equivalence at the functor level
\begin{equation*}
\Sym_{\EOp}(f,A): \Sym_{\EOp}(M,A)\xrightarrow{\sim}\Sym_{\EOp}(N,A),
\end{equation*}
for all $\E$-cofibrant $\EOp$-algebras $A$.
Hence,
in our context,
we have a chain of weak-equivalences
\begin{equation*}
\Sym_{\EOp}(M^n_{\EOp},A)\xleftarrow{\sim}\,\cdot\,\xrightarrow{\sim}\,\cdots\,\xrightarrow{\sim}\Sym_{\EOp}(B^n_{\EOp},A) = B^n(A)
\end{equation*}
between $\Sym_{\EOp}(M^n_{\EOp},A)$ and the iterated bar complex $B^n(A) = \Sym_{\EOp}(B^n_{\EOp},A)$.
\end{proof}

According to this theorem,
the definition of a proper iterated bar complex $B^n: {}_{\EOp}\E\rightarrow\E$, where $\EOp$ is any $E_\infty$-operad,
reduces to the construction of a cofibrant right $\EOp$-module $M^n_{\EOp}$
together with a weak-equivalence $M^n_{\EOp}\circ_{\EOp}\COp\xrightarrow{\sim} B^n_{\COp}$.
In the next section,
we give an effective construction of such a cofibrant right $\EOp$-module $M^n_{\EOp}$
starting from the iterated bar module over the commutative operad $B^n_{\COp}$.

\section{Iterated bar modules as quasi-free modules}\label{QuasiFreeLifting}
Recall that the iterated bar module $B^n_{\COp}$
is given by the iterated bar complex of the commutative operad,
viewed as a commutative algebra
in the category of right modules over itself.
The bar construction $B(A)$ is defined by a twisted complex $B(A) = (T^c(\Sigma A),\partial)$
and so does any of its composite.
From this statement,
we deduce that the $n$-fold bar module $B^n_{\COp}$
is identified with a twisted right $\COp$-module
of the form $B^n_{\COp} = ((T^c\Sigma)^n(\COp),\partial_{\gamma})$,
where we take the $n$-fold composite of the functor $T^c(\Sigma -)$
underlying the bar complex $B(-)$.

We prove in this section that the right $\COp$-module~$(T^c\Sigma)^n(\COp)$
is isomorphic to a composite $(T^c\Sigma)^n(\COp)\simeq T^n\circ\COp$,
for some free $\Sigma_*$-module $T^n$.
We use this structure result
to lift the twisting homomorphism $\partial_{\gamma}: T^n\circ\COp\rightarrow T^n\circ\COp$
to the right $\EOp$-module $T^n\circ\EOp$,
for any $E_\infty$-operad $\EOp$.
We obtain from this construction a cofibrant right $\EOp$-module $B^n_{\EOp} = (T^n\circ\EOp,\partial_{\epsilon})$
such that $B^n_{\EOp}\circ_{\EOp}\COp = B^n_{\COp}$.
Hence,
the lifting construction produces a right $\EOp$-module
which satisfies the requirements of theorem~\ref{IteratedBarModules:IteratedBarEquivalence}
and, as such, determines a good iterated bar complex on the category of $\EOp$-algebras.

We can make our construction effective if we assume that the $E_\infty$-operad $\EOp$
is equipped with an effective contracting chain homotopy $\nu: \EOp\rightarrow\EOp$
such that $\delta(\nu) = \id-\iota\epsilon$,
for some section $\iota: \COp\rightarrow\EOp$ of the augmentation morphism $\epsilon: \EOp\rightarrow\COp$.
We need this effective construction in~\S\ref{EnHomology}
in order to obtain a homotopy interpretation of iterated bar complexes.

The composite $K\circ\ROp$ represents a free object in the category of right $\ROp$-modules.
The twisted objects $M = (K\circ\ROp,\partial)$
associated to a free right $\ROp$-module $K\circ\ROp$
are called quasi-free.
To begin this section,
we review the definition and usual properties of these quasi-free modules.

\subsubsection{Twisting cochains in the category of right $\ROp$-modules}\label{QuasiFreeLifting:TwistingCochains}
In certain constructions,
the natural differential of a dg-module $E$ is twisted by a homomorphism $\partial: E\rightarrow E$ of degree $-1$, called a twisting homomorphism,
to produce a new dg-module $M = (E,\partial)$,
which has the same underlying graded object as $E$,
but whose differential is given by the sum $\delta+\partial: E\rightarrow E$.
The bar complex $B(A) = (T^c(\Sigma A),\partial)$ gives an application of this construction.

To ensure that the homomorphism $\delta+\partial: E\rightarrow E$
satisfies the equation of a differential $(\delta+\partial)^2 = 0$,
we simply have to require that a twisting homomorphism $\partial: E\rightarrow E$
satisfies the equation $\delta(\partial) + \partial^2 = 0$
in $\Hom_{\E}(E,E)$.

The construction of twisted objects makes sense in the category of $\Sigma_*$-modules
and in the category of right modules over an operad,
for a twisting homomorphism $\partial$
in the hom-object of the concerned category $\E = \M,\M{}_{\ROp}$.
The bar module $B_{\ROp} = (T^c(\Sigma\ROp),\partial)$
is an instance of a twisted object in the category of right $\ROp$-modules.

\subsubsection{Free modules over operads and free $\Sigma_*$-modules}\label{QuasiFreeLifting:FreeModules}
Let $K$ be any $\Sigma_*$-module.
The composite $\Sigma_*$-module $K\circ\ROp$
inherits a right $\ROp$-action,
defined by the morphism $K\circ\mu: K\circ\ROp\circ\ROp\rightarrow K\circ\ROp$
induced by the operad composition product $\mu: \ROp\rightarrow\ROp$,
and forms naturally a right $\ROp$-module.
This object $K\circ\ROp$ is identified with a free right $\ROp$-module
associated to $K$,
in the sense that any morphism of $\Sigma_*$-modules $f: K\rightarrow M$, where $M\in\M{}_{\ROp}$,
has a unique extension
\begin{equation*}
\xymatrix{ K\ar[rr]^{f}\ar[dr]_{K\circ\eta} && M \\ & K\circ\ROp\ar@{.>}[ur]_{\tilde{f}} & }
\end{equation*}
such that $\tilde{f}: K\circ\ROp\rightarrow M$
is a morphism of right $\ROp$-modules.
Equivalently, the map $K\mapsto K\circ\ROp$
defines a left adjoint of the functor functor $U: \M{}_{\ROp}\rightarrow\M$.

Intuitively,
the object $K\circ\ROp$
is spanned by formal composites $\xi(p_1,\dots,p_r)$
of a generating element $\xi\in K(r)$
with operations $p_1,\dots,p_r\in\ROp$.
The extension to the free right $\ROp$-module $K\circ\ROp$
of a morphism of $\Sigma_*$-modules $f: K\rightarrow M$
is determined by the formula $\tilde{f}(\xi(p_1,\dots,p_r)) = f(\xi)(p_1,\dots,p_r)$.

Let $\C^{\NN}$
denote the category of collections $G = \{G(r)\}_{r\in\NN}$,
where $G(r)\in\C$.
The forgetful functor $U: \M\rightarrow\C^{\NN}$
has a left adjoint which maps a collection $G$ to an associated free $\Sigma_*$-module,
denoted by $\Sigma_*\otimes G$.
This $\Sigma_*$-module is represented by the external tensor products:
\begin{equation*}
(\Sigma_*\otimes G)(r) = \Sigma_r\otimes G(r).
\end{equation*}
The adjunction unit $\eta: G\rightarrow\Sigma_*\otimes G$
identifies $G(r)$ with the summand $\id\otimes G(r)$
of $\Sigma_*\otimes G$.
The functor on finite sets equivalent to $\Sigma_*\otimes G$
satisfies
\begin{equation*}
(\Sigma_*\otimes G)(\eset) = \Bij(\{1,\dots,r\},\eset)\otimes G(r),
\end{equation*}
for every set $\eset$ such that $\eset = \{e_1,\dots,e_r\}$.

In the case of a free $\Sigma_*$-module $K = \Sigma_*\otimes G$,
the composite $K\circ\ROp$
has an expansion of the form
\begin{equation*}
(\Sigma_*\otimes G)\circ\ROp = \Sym(\Sigma_*\otimes G,\ROp) = \bigoplus_{r=0}^{\infty} G(r)\otimes\ROp^{\otimes r},
\end{equation*}
where no coinvariant occurs.

\subsubsection{Quasi-free modules}\label{QuasiFreeLifting:QuasiFreeModules}
By convention,
a quasi-free right $\ROp$-module $M$ refers to a twisted object $M = (K\circ\ROp,\partial)$
formed from a free right $\ROp$-module $K\circ\ROp$.
In~\cite[\S 14.2]{Bar0},
we prove that a quasi-free $\ROp$-module $M = (K\circ\ROp,\partial)$
is cofibrant
if we have $K = \Sigma_*\otimes G$
for a collection of free graded $\kk$-modules $G = \{G(r)\}_{r\in\NN}$
equipped with a good filtration (see~\emph{loc. cit.} for details).
The filtration condition is automatically satisfied
if we assume that each $G(r)$ is non-negatively graded (in this case we can apply the arguments of~\emph{loc. cit.}
to the degreewise filtration).

The goal of this section is to prove that the iterated bar module $B^n_{\ROp}$,
is defined by a cofibrant quasi-free module of this form $B^n_{\ROp} = (T^n\circ\ROp,\partial)$,
for some free $\Sigma_*$-module $T^n$.
For this purpose,
we use that the twisting homomorphism $\partial: K\circ\ROp\rightarrow K\circ\ROp$
is determined by a homomorphism $\alpha: G\rightarrow K\circ\ROp$
in the category of collections $\E = \C^{\NN}$
whenever $K = \Sigma_*\otimes G$.
Indeed,
the adjunction relations
$\C^{\NN}\rightleftarrows\M\rightleftarrows\M{}_{\ROp}$
yield isomorphisms of dg-modules
\begin{equation*}
\Hom_{\C^{\NN}}(G,K\circ\ROp)
\simeq\Hom_{\M}(K,K\circ\ROp)
\simeq\Hom_{\M{}_{\ROp}}(K\circ\ROp,K\circ\ROp),
\end{equation*}
for $K = \Sigma_*\otimes G$.

Let $\partial_\alpha: K\circ\ROp\rightarrow K\circ\ROp$ denote the homomorphism of right $\ROp$-modules
equivalent to $\alpha: G\rightarrow K\circ\ROp$.
The composite $\partial_{\alpha}\partial_{\beta}: K\circ\ROp\rightarrow K\circ\ROp$
is necessarily associated to a homomorphism of collections
for which we adopt the notation $\alpha\dercirc\beta: G\rightarrow K\circ\ROp$.
The equation of twisting homomorphisms $\delta(\partial_{\alpha}) + \partial_{\alpha}^2 = 0$
is equivalent to the equation $\delta(\alpha) + \alpha\dercirc\alpha = 0$
in $\Hom_{\C^{\NN}}(G,K\circ\ROp)$.

The homomorphism $\alpha\dercirc\beta: G\rightarrow K\circ\ROp$
can be identified with the composite
\begin{equation*}
G\xrightarrow{\beta} K\circ\ROp\xrightarrow{\partial_\alpha} K\circ\ROp.
\end{equation*}
Intuitively,
the twisting homomorphism $\partial_{\alpha}$
associated to a homomorphism $\alpha: G\rightarrow K\circ\ROp$
is determined by the relation $\partial_{\alpha}(\xi(p_1,\dots,p_r)) = \alpha(\xi)(p_1,\dots,p_r)$
for any formal composite $\xi(p_1,\dots,p_r)\in K\circ\ROp$,
where $\xi\in G(r)$ and $p_1,\dots,p_r\in\ROp$.

\subsubsection{Lifting twisting homomorphisms of quasi-free modules}\label{QuasiFreeLifting:TwistingCochainLifting}
Suppose we have a quasi-free module $N = (K\circ\SOp,\partial_\beta)$
such that $K = \Sigma_*\otimes G$,
for a collection of free graded $\kk$-modules $G = \{G(r)\}_{r\in\NN}$.
Suppose further that each $G(r)$ is non-negatively graded.

Let $\psi: \ROp\rightarrow\SOp$ be an acyclic fibration of operads.
Suppose we have a section $\iota: \SOp(\eset)\rightarrow\ROp(\eset)$
and a contracting chain homotopy $\nu: \ROp(\eset)\rightarrow\ROp(\eset)$
such that $\psi\cdot\iota = \id$, $\psi\cdot\nu = 0$
and $\delta(\nu) = \id - \iota\cdot\psi$,
for every finite set $\eset$
equipped with an ordering $\eset = \{e_1<\dots<e_r\}$.
Note that such maps $\iota$ and $\nu$
are ensured to exist when $\ROp$ and $\SOp$
are non-negatively graded operads.

The section $\iota: \SOp(\eset)\rightarrow\ROp(\eset)$
and the contracting chain homotopy $\nu: \SOp(\eset)\rightarrow\SOp(\eset)$
have natural extensions
$\tilde{\iota}: \SOp^{\otimes r}(\eset)\rightarrow\ROp^{\otimes r}(\eset)$
and
$\tilde{\nu}: \ROp^{\otimes r}(\eset)\rightarrow\ROp^{\otimes r}(\eset)$
given by the tensor products
\begin{equation*}
\tilde{\iota} = \iota^{\otimes r}\\
\quad\text{and}
\quad\tilde{\nu} = \sum_{i=1}^{r} (-1)^{i-1} (\iota\psi)^{\otimes i-1}\otimes\nu\otimes\id^{\otimes r-i+1}
\end{equation*}
on the summands
of the tensor power:
\begin{equation*}
\SOp^{\otimes r}(\eset) = \bigoplus_{{\eset_1}\amalg\dots\amalg{\eset_r} = \eset} \SOp(\eset_1)\otimes\dots\otimes\SOp(\eset_r).
\end{equation*}
In this definition,
we use that each subset $\eset_i\subset\eset$
inherits a natural ordering from $\eset$.

Since $K\circ M = \bigoplus_{r=0}^{\infty} G(r)\otimes M^{\otimes r}$
for every $\Sigma_*$-object $M$,
the extension of $\iota$ and $\nu$ to tensor powers
gives also rise to a section $\tilde{\iota}: K\circ\SOp(\eset)\rightarrow K\circ\ROp(\eset)$
and a contracting chain homotopy $\tilde{\nu}: K\circ\ROp(\eset)\rightarrow K\circ\ROp(\eset)$
such that $K\circ\psi\cdot\tilde{\iota} = \id$, $K\circ\psi\cdot\tilde{\nu} = 0$
and $\delta(\tilde{\nu}) = \id - \tilde{\iota}\cdot K\circ\psi$.

For any $r\in\NN$,
we pick a lifting
\begin{equation*}
\xymatrix{ & K\circ\ROp(r)\ar[d]^{K\circ\psi} \\
G(r)\ar[r]_{\beta}\ar@{.>}[ur]!DL^{\alpha_0} & K\circ\SOp(r) }
\end{equation*}
by setting $\alpha_0 = \tilde{\iota}\cdot\beta$,
where we consider the section extension $\tilde{\iota}: K\circ\SOp(\rset)\rightarrow K\circ\ROp(\rset)$
defined with respect to the canonical ordering $\rset = \{1<\dots<r\}$.
Define a sequence of homomorphisms $\alpha_m: G(r)\rightarrow K\circ\SOp(r)$, $m\in\NN$,
by the inductive formula $\alpha_m = \sum_{p+q = m-1} \tilde{\nu}\cdot(\alpha_p\dercirc\alpha_q)$,
where we consider the chain homotopy extension $\tilde{\nu}: K\circ\ROp(\rset)\rightarrow K\circ\ROp(\rset)$
defined with respect to the canonical ordering $\rset = \{1<\dots<r\}$ too.
Recall that $\alpha_p\dercirc\alpha_q$
is given by the composite $\partial_{\alpha_p}\cdot\alpha_q: G\rightarrow K\circ\ROp$,
where $\partial_{\alpha_p}: K\circ\ROp\rightarrow K\circ\ROp$
is the homomorphism of right $\ROp$-modules extending $\alpha_p: G\rightarrow K\circ\ROp$.

Form the homomorphism $\alpha_* = \sum_{m=0}^{\infty}\alpha_m$.
Note that $\alpha_m$ decreases the degree in $K$ by $m$.
Thus the infinite sum $\alpha_* = \sum_{m=0}^{\infty}\alpha_m$
makes sense since $K$ is supposed to vanish in degree $*<0$.

\medskip
This lifting process returns the following result:

\begin{prop}\label{QuasiFreeLifting:QuasiFreeModuleExtension}
The homomorphism $\alpha_*: G\rightarrow K\circ\ROp$
determines a twisting homomorphism $\partial_\alpha: K\circ\ROp\rightarrow K\circ\ROp$
such that the diagram
\begin{equation*}
\xymatrix{ K\circ\ROp\ar@{.>}[r]^{\partial_\alpha}\ar[d]_{K\circ\psi} & K\circ\ROp\ar[d]^{K\circ\psi} \\
K\circ\SOp\ar[r]_{\partial_{\beta}} & K\circ\SOp }
\end{equation*}
commutes.
The quasi-free module $M = (K\circ\ROp,\partial_\alpha)$
defined by this twisting homomorphism satisfies the extension relation $M\circ_{\ROp}\SOp\simeq N$
with respect to the given quasi-free module $N = (K\circ\SOp,\partial_\beta)$.
\end{prop}

Note that the composite $K\circ\psi$ defines a morphism
of right $\ROp$-modules
\begin{equation*}
K\circ\psi: (K\circ\ROp,\partial_\alpha)\rightarrow(K\circ\SOp,\partial_\beta).
\end{equation*}
This assertion is an immediate consequence of the commutativity of the diagram.
The isomorphism $M\circ_{\ROp}\SOp\simeq N$ corresponds to $K\circ\psi$ under the adjunction relation
between extension and restriction functors.

\begin{proof}
The equation of twisting homomorphisms $\delta(\alpha_*) + \alpha_*\dercirc\alpha_* = 0$
follows from an immediate induction.
The assumption $\psi\nu = 0$ implies $K\circ\psi\cdot\tilde{\nu} = 0$,
from which we deduce $K\circ\psi\cdot\alpha_n = 0$ for $n>0$.
By definition,
we also have $K\circ\psi\cdot\alpha_0 = \beta$.
Hence we obtain $\psi_*\alpha_* = \beta$
and we deduce from this relation that the twisting homomorphisms equivalent to $\alpha$ and $\beta$
satisfy the relation $K\circ\psi\cdot\partial_{\alpha} = \partial_{\beta}\cdot K\circ\psi$.

This commutation relation implies the identity $\partial_{\alpha}\circ_{\ROp}\SOp = \partial_{\beta}$
from which we deduce:
\begin{equation*}
(K\circ\ROp,\partial_\alpha)\circ_{\ROp}\SOp
=((K\circ\ROp)\circ_{\ROp}\SOp,\partial_\alpha\circ_{\ROp}\SOp)
=(K\circ\SOp,\partial_\beta).
\end{equation*}
This assertion achieves the proof of proposition~\ref{QuasiFreeLifting:QuasiFreeModuleExtension}.
\end{proof}

\subsubsection{The iterated bar complex of commutative algebras
and the iterated bar module over the commutative operad}\label{QuasiFreeLifting:CommutativeCase}
The bar construction is a twisted dg-module by definition
and, as a consequence, so does any of its composite.
Thus,
for a commutative algebra $A$,
we have an identity:
\begin{equation*}
B^n(A) = ((T^c\Sigma)^n(A),\partial),
\end{equation*}
where the twisting homomorphism $\partial$ integrates all terms yielded by bar coderivations,
occurring at each level of the composite $B^n(A) = B\circ\dots\circ B(A)$.

In the case of the commutative operad $A = \COp$,
we obtain that the iterated bar module $B^n_{\COp}$
is identified with a twisted right $\COp$-module
of the form $B^n_{\COp} = ((T^c\Sigma)^n(\COp),\partial_{\gamma})$.
From this representation of $B^n_{\COp}$,
we obtain:

\begin{prop}\label{QuasiFreeLifting:StructureCommutativeCase}
The iterated bar module $B^n_{\COp}$ is identified with a quasi-free right $\COp$-module
of the form
$B^n_{\COp} = (T^n\circ\COp,\partial_{\gamma})$,
where
\begin{equation*}
T^n = (T^c\Sigma)^n(I)
\end{equation*}
is given by $n$ iterations, within the category of $\Sigma_*$-modules,
of the functor $T^c(\Sigma -)$ applied to the composition unit $I$.
\end{prop}

\begin{proof}
Recall that we have a distribution relation $(M\otimes N)\circ P\simeq(M\circ P)\otimes(N\circ P)$
between the tensor product and the composition product in the category of $\Sigma_*$-modules.
For an iterated tensor coalgebra,
we obtain an isomorphism of right $\COp$-modules
\begin{equation*}
(T^c\Sigma)^n(\COp)\simeq(T^c\Sigma)^n(I\circ\COp)\simeq(T^c\Sigma)^n(I)\circ\COp.
\end{equation*}
The conclusion follows immediately.
\end{proof}

Note further:

\begin{prop}\label{QuasiFreeLifting:SigmaFreeStructure}
The $\Sigma_*$-module $T^n = (T^c\Sigma)^n(I)$ is a free $\Sigma_*$-module
\begin{equation*}
T^n = \Sigma_*\otimes G^n
\end{equation*}
associated to a collection of non-negatively graded $\kk$-modules $G^n(r)$, $r\in\NN$,
defined inductively by:
\begin{multline*}
G^0(r) = \begin{cases} \kk, & \text{if $r=1$}, \\
0, & \text{otherwise}, \end{cases} \\
\text{and}\quad G^n(r) = \bigoplus_{d,r_1+\dots+r_d = r}\Sigma G^{n-1}(r_1)\otimes\dots\otimes\Sigma G^{n-1}(r_d),
\quad\text{for $n>0$}.
\end{multline*}
\end{prop}

The embedding $G^n(r) = \id\otimes G^n(r)\subset T^n(r)$ which yields the isomorphism $\Sigma_*\otimes G^n = T^n$
is defined in the proof of the proposition.

\begin{proof}
The identity $T^n = \Sigma_*\otimes G^n$ is obvious for $n = 0$.

For $n>0$,
we have a canonical morphism $\eta: G^n(r)\rightarrow T^n(r)$
formed inductively by tensor products
\begin{multline*}
G^n(r) = \bigoplus_{d,r_1+\dots+r_d = r}\Sigma G^{n-1}(r_1)\otimes\dots\otimes\Sigma G^{n-1}(r_d)\\
\begin{aligned}
& \xrightarrow{\eta^{\otimes *}}\bigoplus_{d,r_1+\dots+r_d = r}\Sigma T^{n-1}(r_1)\otimes\dots\otimes\Sigma T^{n-1}(r_d)\\
& \simeq\bigoplus_{d,r_1+\dots+r_d = r}\Sigma T^{n-1}(\eset_1)\otimes\dots\otimes\Sigma T^{n-1}(\eset_d)
\subset T^c(\Sigma T^{n-1})(e),
\end{aligned}
\end{multline*}
where we take the obvious ordered partition:
\begin{equation*}
\eset_i = \{r_1+\dots+r_{i-1}+1,\dots,r_1+\dots+r_{i-1}+r_i\},\quad\text{for $i = 1,\dots,d$}.
\end{equation*}
The identity $T^n = \Sigma_*\otimes G^n$ follows from an easy induction on $n$.
\end{proof}

From proposition~\ref{QuasiFreeLifting:StructureCommutativeCase}
and proposition~\ref{QuasiFreeLifting:SigmaFreeStructure},
we conclude:

\begin{thm}\label{QuasiFreeLifting:BarTwistingCochainLifting}
Let $\EOp$ be an $E_\infty$-operad such that, for every finite ordered set $\eset = \{e_1<\dots<e_r\}$,
we have a section $\iota: \COp(\eset)\rightarrow\EOp(\eset)$ of the augmentation $\epsilon: \EOp(\eset)\rightarrow\COp(\eset)$
and a contracting chain homotopy $\nu: \EOp(\eset)\rightarrow\EOp(\eset)$
satisfying $\epsilon\cdot\iota = \id$, $\epsilon\cdot\nu = 0$ and $\delta(\nu) = \id - \iota\cdot\epsilon$.

Then the construction of~\S\S\ref{QuasiFreeLifting:TwistingCochainLifting}-\ref{QuasiFreeLifting:QuasiFreeModuleExtension}
can be applied to the iterated bar module $B^n_{\COp} = (T^n\circ\COp,\partial_\gamma)$
to produce a cofibrant quasi-free module $B^n_{\EOp} = (T^n\circ\EOp,\partial_\epsilon)$
satisfying $B^n_{\EOp}\circ_{\EOp}\COp\simeq B^n_{\COp}$.
\end{thm}

The Barratt-Eccles operad, whose definition is reviewed in~\ref{IteratedBarModules:BarrattEccles},
fulfils the requirement of this theorem.

\medskip
We have moreover:

\begin{prop}\label{QuasiFreeLifting:CoderivationRelation}
In the case of iterated bar modules,
the twisting homomorphism $\partial_{\epsilon}: T^n\circ\EOp\rightarrow T^n\circ\EOp$
which arises from the construction
of~\S\S\ref{QuasiFreeLifting:TwistingCochainLifting}-\ref{QuasiFreeLifting:QuasiFreeModuleExtension}
satisfies the coderivation relation
$\Delta\partial_{\epsilon} = (\partial_{\epsilon}\otimes\Id+\Id\otimes\partial_{\epsilon})\Delta$
so that $B^n_{\EOp} = (T^n\circ\EOp,\partial_{\epsilon})$
forms a coalgebra.
\end{prop}

The diagonal $\Delta: T^n\circ\EOp\rightarrow T^n\circ\EOp$
comes from the deconcatenation coproduct of the first tensor coalgebra
in the composite $T^n\circ\EOp = T^c(\Sigma T^{n-1})\circ\EOp = T^c(\Sigma T^{n-1}\circ\EOp)$.

\begin{proof}
The twisting homomorphism $\partial_{\gamma}: T^n\circ\COp\rightarrow T^n\circ\COp$
associated to the commutative operad
forms a coderivation with respect to the coproduct of $T^n\circ\COp$
because the bar module $B^n_{\COp}$ forms a coalgebra,
like the bar complex of any commutative algebra.

The deconcatenation of tensors
yield a diagonal
\begin{equation*}
\Delta: G^n(t)\rightarrow\bigoplus_{r+s = t} G^n(r)\otimes G^n(s),
\end{equation*}
on the collection
\begin{equation*}
G^n(r) = \bigoplus_{d,r_1+\dots+r_d = r}\Sigma G^{n-1}(r_1)\otimes\dots\otimes\Sigma G^{n-1}(r_d)
\end{equation*}
of proposition~\ref{QuasiFreeLifting:SigmaFreeStructure}.
The deconcatenation coproduct of $T^n\circ\ROp = T^c(\Sigma T^{n-1}\circ\ROp)$, $\ROp = \COp,\EOp$,
is clearly identified with the diagonal
\begin{multline*}
G^n(t)\otimes\ROp^{\otimes t}
\rightarrow\bigoplus_{r+s = t} (G^n(r)\otimes G^n(s))\otimes\ROp^{\otimes t}
\simeq\bigoplus_{r+s = t} (G^n(r)\otimes\ROp^{\otimes r})\otimes(G^n(s)\otimes\ROp^{\otimes s})
\end{multline*}
induced by the deconcatenation coproduct of~$G^n$.

One checks easily that the homomorphisms of~\S\ref{QuasiFreeLifting:TwistingCochainLifting}
\begin{equation*}
\tilde{\iota}: G^n(r)\otimes\SOp^{\otimes r}\rightarrow G^n(r)\otimes\ROp^{\otimes r}
\quad\text{and}
\quad\tilde{\nu}: G^n(r)\otimes\ROp^{\otimes r}\rightarrow G^n(r)\otimes\ROp^{\otimes r}
\end{equation*}
satisfy
$\Delta\cdot\tilde{\iota} = \tilde{\iota}\otimes\tilde{\iota}\cdot\Delta$
and
$\Delta\cdot\tilde{\nu} = (\tilde{\nu}\otimes\id+\tilde{\iota}\psi\otimes\tilde{\nu})\cdot\Delta$.
Then
an easy induction shows that each $\epsilon_*: G^n\rightarrow T^n\circ\EOp$
forms a coderivation.
The conclusion follows immediately.
\end{proof}

The iterated bar modules $B^n_{\EOp}$ which arise from the construction
of~\S\ref{IteratedBarModules}
are connected by suspension morphisms
$\sigma_{\EOp}: \Sigma B^{n-1}_{\EOp}\rightarrow B^{n}_{\EOp}$.
So do the iterated bar modules yielded by the construction
of theorem~\ref{QuasiFreeLifting:BarTwistingCochainLifting}:

\begin{prop}\label{QuasiFreeLifting:Suspension}
The iterated bar modules of theorem~\ref{QuasiFreeLifting:BarTwistingCochainLifting}
are connected by suspension morphisms
\begin{equation*}
\sigma_{\EOp}: \Sigma B^{n-1}_{\EOp}\rightarrow B^{n}_{\EOp}
\end{equation*}
that fit commutative diagrams
\begin{equation*}
\xymatrix{ \Sigma B^{n-1}_{\EOp}\ar[d]\ar@{.>}[r]^{\sigma_{\EOp}} & B^n_{\EOp}\ar[d]\\
\Sigma B^{n-1}_{\COp}\ar[r]_{\sigma_{\COp}} & B^n_{\COp} },
\end{equation*}
for every $n>0$.
We have moreover $\sigma_{\EOp}\circ_{\EOp}\COp = \sigma_{\COp}$.
\end{prop}

This proposition is an immediate consequence of the next lemma:

\begin{lemm}\label{QuasiFreeLifting:SuspensionDifferential}
The canonical embedding
\begin{equation*}
\Sigma T^{n-1}\circ\EOp\hookrightarrow T^c(\Sigma T^{n-1}\circ\EOp)\simeq T^n\circ\EOp
\end{equation*}
commutes with the twisting homomorphism of theorem~\ref{QuasiFreeLifting:BarTwistingCochainLifting},
for every $n>0$.
\end{lemm}

\begin{proof}
Observe that the embedding of the lemma
is realized at the level of the generating collection $G^n$
of proposition~\ref{QuasiFreeLifting:SigmaFreeStructure}
and is clearly preserved by the homomorphisms $\tilde{\iota}$ and $\tilde{\nu}$
of~\S\ref{QuasiFreeLifting:TwistingCochainLifting}.
The lemma follows by an easy induction
from the definition of the twisting homomorphism $\partial_{\epsilon}$.
\end{proof}

We have moreover:

\begin{prop}\label{QuasiFreeLifting:SuspensionCofibration}
The suspension morphisms
of proposition~\ref{QuasiFreeLifting:Suspension}
\begin{equation*}
\sigma_{\EOp}: \Sigma B^{n-1}_{\EOp}\rightarrow B^{n}_{\EOp}
\end{equation*}
are cofibrations in the category of right $\EOp$-modules.
\end{prop}

\begin{proof}
Use the degreewise filtration of the $\Sigma_*$-module $T^n$
to split the embedding of quasi-free modules
\begin{equation*}
\Sigma(T^{n-1}\circ\EOp,\partial_{\epsilon})\hookrightarrow(T^n\circ\EOp,\partial_{\epsilon})
\end{equation*}
into a sequence of generating cofibrations of right $\EOp$-modules
(use the overall ideas of~\cite[\S 11.2,\S 14.2]{Bar0}).
\end{proof}

In the next sections,
we use the effective construction
of~\S\S\ref{QuasiFreeLifting:TwistingCochainLifting}-\ref{QuasiFreeLifting:QuasiFreeModuleExtension}
to prove that, for some good $E_\infty$-operads $\EOp$ equipped with a filtration of the form (*),
the twisting homomorphism $\partial_\epsilon: T^n\circ\EOp\rightarrow T^n\circ\EOp$
preserves the subobject $T^n\circ\EOp_n$
and restricts to a twisting homomorphism on this right $\EOp_n$-module $T^n\circ\EOp_n$.
For this aim,
we use the existence of nice cell decompositions $\EOp = \colim_{\kappa}\EOp_\kappa$
refining the filtration of $\EOp$.
Therefore we recall the overall definition of these cell structures
before going further into the study of iterated bar complexes.

\section{Interlude: operads shaped on complete graph posets}\label{CompleteGraphOperad}
For our analysis of iterated bar modules,
we are going to use a particular cell structure of $E_\infty$-operads, introduced in~\cite{BergerCell},
and modelled by a certain operad in posets $\K$,
the complete graph operad.
The main purpose of this section is to revisit definitions of~\cite{BergerCell}
in order to give an abstract formalization
of the complete graph cell decompositions of $E_\infty$-operads
and to extend the applications of these cell structures
to iterated bar modules.

To begin with,
we review the definition of the complete graph operad.
For our needs,
it is more convenient to define directly the functor on finite sets
underlying the complete graph operad $\K$.

\subsubsection{The complete graph posets}\label{CompleteGraphOperad:PosetStructure}
Let $\eset$ be any set with $r$ elements $\{e_1,\dots,e_r\}$.
The complete graph poset $\K(\eset)$
consists of pairs $\kappa = (\mu,\sigma)$,
such that $\mu$ is a collection of non-negative integers $\mu_{e f}\in\NN$
indexed by pairs $\{e,f\}\subset\eset$
and $\sigma$ is a bijection $\sigma: \{1,\dots,r\}\rightarrow\eset$,
which amounts to ordering the set~$\eset$.

For a pair $\{e,f\}\subset\eset$,
we define the restriction $\sigma|_{e f}$ of an ordering $\sigma$
as the ordering of~$\{e,f\}$
defined by the occurrences of~$\{e,f\}$
in the sequence $\sigma = (\sigma(1),\dots,\sigma(r))$.

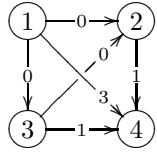
\begin{figure}[h]%
\[\vcenter{\xymatrix@M=0pt@!R=6mm@!C=6mm{ *+<3mm>[o][F]{1}\ar[r]|*+<2pt>{\scriptstyle 0}
\ar[d]|*+<2pt>{\scriptstyle 0}
\ar[dr]|(0.7)*+<2pt>{\scriptstyle 3}
& *+<3mm>[o][F]{2}
\ar[d]|*+<2pt>{\scriptstyle 1} \\
*+<3mm>[o][F]{3}\ar[r]|*+<2pt>{\scriptstyle 1}
\ar[ur]|*+<1mm>{\hole}|(0.7)*+<2pt>{\scriptstyle 0} &
*+<3mm>[o][F]{4} }}\]
\caption{An element of the complete graph operad}\label{Figure:WeightedDirectedGraph}
\end{figure}

The elements of~$\K(\eset)$ are represented by complete graphs on $r$-vertices, indexed by $\eset$,
whose edges are coherently oriented and equipped with a weight (see figure~\ref{Figure:WeightedDirectedGraph}).
The weight of the edge $\{e,f\}$
is defined by the integer $\mu_{e f}\in\NN$.
The orientation of~$\{e,f\}$
is defined by the ordering $\sigma|_{e f}\in\{(e,f),(f,e)\}$.
The coherence of the orientations amounts to the requirement
that the local orderings $\sigma|_{e f}\in\{(e,f),(f,e)\}$
assemble to a global ordering of the set $\{e_1,\dots,e_r\}$.

For elements $(\mu,\sigma),(\nu,\tau)\in\K(r)$,
we set
\begin{equation*}
(\mu,\sigma)\leq(\nu,\tau)
\end{equation*}
if we have
\begin{equation*}
(\mu_{e f}<\nu_{e f})\quad\text{or}\quad(\mu_{e f},\sigma|_{e f}) = (\nu_{e f},\tau|_{e f}),
\end{equation*}
for every pair $\{e,f\}\subset\eset$.
This relation defines clearly a poset structure on~$\K(\eset)$.

The collection of posets $\K(\eset)$
defines clearly a functor on the category of finite sets and bijections.
The morphism $u_*: \K(\eset)\rightarrow\K(\fset)$
induced by a bijection $u: \eset\rightarrow\fset$
is simply defined by reindexing the vertices
of complete graphs.

In the case $\eset = \{1,\dots,r\}$,
we simply replace the bijection $\sigma: \{1,\dots,r\}\rightarrow\eset$
by an equivalent permutation $\sigma\in\Sigma_r$
in the definition of $\K(r) = \K(\{1,\dots,r\})$.

\subsubsection{The complete graph operad}\label{CompleteGraphOperad:OperadStructure}
The collection of posets
\begin{equation*}
\K(r) = \K(\{1,\dots,r\})
\end{equation*}
are equipped with an operad structure.
The action of permutations $w_*: \K(r)\rightarrow\K(r)$
arises from the reindexing process of the previous paragraph.

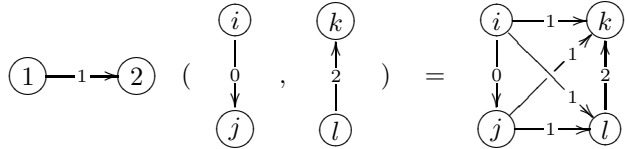
\begin{figure}[h]%
\[\vcenter{\xymatrix@M=0pt@!R=6mm@!C=6mm{ *+<3mm>[o][F]{1}\ar[r]|*+<2pt>{\scriptstyle 1} & *+<3mm>[o][F]{2} }}
\quad(\quad\vcenter{\xymatrix@M=0pt@!R=6mm@!C=6mm{ *+<3mm>[o][F]{i}\ar[d]|*+<2pt>{\scriptstyle 0} \\ *+<3mm>[o][F]{j} }}\quad,
\quad\vcenter{\xymatrix@M=0pt@!R=6mm@!C=6mm{ *+<3mm>[o][F]{k} \\ *+<3mm>[o][F]{l}\ar[u]|*+<2pt>{\scriptstyle 2} }}\quad)
\quad=\quad\vcenter{\xymatrix@M=0pt@!R=6mm@!C=6mm{ *+<3mm>[o][F]{i}\ar[r]|*+<2pt>{\scriptstyle 1}
\ar[d]|*+<2pt>{\scriptstyle 0}
\ar[dr]|(0.7)*+<2pt>{\scriptstyle 1}
& *+<3mm>[o][F]{k} \\
*+<3mm>[o][F]{j}\ar[r]|*+<2pt>{\scriptstyle 1}
\ar[ur]|*+<1mm>{\hole}|(0.7)*+<2pt>{\scriptstyle 1} &
*+<3mm>[o][F]{l}\ar[u]|*+<2pt>{\scriptstyle 2} }}\]
\caption{A composite in the complete graph operad}\label{Figure:GraphComposite}
\end{figure}

The operadic composite $\kappa(\pi_1,\dots,\pi_r)\in\K(\eset_1\amalg\dots\amalg\eset_r)$
is defined by the substitution of the vertices of~$\kappa\in\K(r)$
by the complete graphs $\pi_i\in\K(\eset_i)$, $i = 1,\dots,r$.
Explicitly,
the weight and the orientation of the edges of $\kappa(\pi_1,\dots,\pi_r)$
are determined by the following rules:
the edges of $\kappa(\pi_1,\dots,\pi_r)$
between vertices $e,f\in\eset_1\amalg\dots\amalg\eset_r$
such that $e,f\in\eset_i$, for some $i\in\{1,\dots,r\}$,
are copies of the edge $\{e,f\}$ of the graph $\pi_i$;
the edges of $\kappa(\pi_1,\dots,\pi_r)$
between vertices $e,f\in\eset_1\amalg\dots\amalg\eset_r$ such that $e\in\eset_i$ and $f\in\eset_j$, for a pair $i\not=j$,
are copies of the edge $\{i,j\}$ of the graph $\kappa$.
An example is represented in figure~\ref{Figure:GraphComposite}.

Since we only consider non-unitary operads
in this article,
we adopt the convention $\K(0) = \emptyset$, different from~\cite{BergerCell},
for the complete graph operad.

\medskip
The next observation is a simple consequence of the definition
of the composition of complete graphs:

\begin{obsv}\label{CompleteGraphOperad:CompositeBound}
Define the restriction $\mu|_{\eset_i}$ of an element $\mu\in\K(\eset_1\amalg\dots\amalg\eset_r)$
as the subgraph of $\mu$ generated by the vertices of $\eset_i$.

For a composite $\kappa(\pi_1,\dots,\pi_r)\in\K(\eset_1\amalg\dots\amalg\eset_r)$,
where $\kappa\in\K(r)$ and $\pi_1\in\K(\eset_1),\dots,\pi_r\in\K(\eset_r)$,
we have $\kappa(\pi_1,\dots,\pi_r)\leq\mu$
if and only if $\pi_1\leq\mu|_{\eset_1},\dots,\pi_r\leq\mu|_{\eset_r}$
and $\kappa(\mu|_{\eset_1},\dots,\mu|_{\eset_r})\leq\mu$.
\end{obsv}

\subsubsection{Operads shaped on complete graph posets}\label{CompleteGraphOperad:KOperads}
Define a $\K$-operad as a collection of $\K(r)$-diagrams $\{\POp_{\kappa}\}_{\kappa\in\K(r)}$
together with $\Sigma_r$-actions
\begin{equation*}
\POp_{\kappa}\xrightarrow{w_*}\POp_{w\kappa},\quad w\in\Sigma_*,
\end{equation*}
and composition products
\begin{equation*}
\POp_{\kappa}\otimes(\POp_{\pi_1}\otimes\dots\otimes\POp_{\pi_r})
\rightarrow\POp_{\kappa(\pi_1,\dots,\pi_r)}
\end{equation*}
which satisfy a natural extension of the standard axioms of operads.
Naturally,
the $\K(r)$-diagrams $\{\POp_{\kappa}\}_{\kappa\in\K(r)}$
are equivalent to $\K(\eset)$-diagrams $\{\POp_{\kappa}\}_{\kappa\in\K(\eset)}$,
associated to all finite sets $\eset$,
so that any bijection $f: \eset\rightarrow\fset$
defines a morphism of diagrams $f_*: \POp_{\kappa}\rightarrow\POp_{f_*\kappa}$.

We adopt the notation ${}_{\K}\Op$ to refer to the category of $\K$-operads.
Naturally,
a morphism of $\K$-operads consists of a collection of $\K(r)$-diagram morphisms $\phi: \POp_{\kappa}\rightarrow\QOp_{\kappa}$
which commute with $\Sigma_r$-actions and composition structures.

The colimit of the underlying $\K(r)$-diagrams of a $\K$-operad
\begin{equation*}
(\colim_{\K}\POp)(r) = \colim_{\kappa\in\K(r)}\POp_{\kappa}
\end{equation*}
inherits a natural operad structure.
Hence
we have a functor $\colim_{\K}: {}_{\K}\Op\rightarrow\Op$.
This colimit functor is left adjoint to the obvious functor $\cst: \Op\rightarrow{}_{\K}\Op$
which maps an operad $\POp$ to the constant $\KOp$-operad
such that $\POp_{\kappa} = \POp(r)$
for every $\kappa\in\K(r)$.

We say that an operad $\POp\in\Op$ is equipped with a $\K$-structure
if we have a $\K$-operad $\POp_{\kappa}$ such that $\colim_{\K}\POp_{\kappa}\xrightarrow{\simeq}\POp$
for the adjoint morphism of an embedding of the $\K$-operad $\POp_{\kappa}$
into the constant $\K$-operad defined by $\POp\in\Op$.
In~\cite{Bar3},
we use a more general definition,
where we do not necessarily assume that the components of $\K$-operad $\POp_{\kappa}$
are mapped injectively into $\POp$.
The convention $\POp_{\kappa}\subset\POp$ will simplify the presentation
of the constructions of this article.

In general,
when an operad $\POp$ is equipped with a $\K$-structure
we use the same notation for the genuine operad $\POp\in\Op$
and the underlying $\K$-operad $\POp_{\kappa}$
associated to $\POp$.

\subsubsection{The example of the Barratt-Eccles operad}\label{CompleteGraphOperad:BarrattEccles}
The Barratt-Eccles operad is an instance of an operad equipped with a nice $\K$-structure.
Recall that the dg-version of this operad consists of the normalized chain complexes $\EOp(r) = N_*(E\Sigma_r)$
so that the dg-module $\EOp(r)$ is spanned in degree $d$
by the $d$-simplices of non-degenerate permutations $(w_0,\dots,w_d)\in\Sigma_r\times\dots\times\Sigma_r$.

For $\kappa = (\mu,\sigma)\in\KOp(r)$,
we form the module $\EOp_{\kappa}\subset\EOp(r)$
spanned by the simplices of permutations $(w_0,\dots,w_d)$
such that, for every pair $\{i,j\}\subset\{1,\dots,r\}$,
the sequence $(w_0|_{i j},\dots,w_d|_{i j})$
has less than $\mu_{i j}$ variations,
or has exactly $\mu_{i j}$ variations
and satisfies $w_d|_{i j} = \sigma|_{i j}$.
This module is clearly preserved by the differential of the Barratt-Eccles operad
and hence forms a dg-submodule of~$\EOp(r)$.
Moreover, we have clearly $\kappa\leq\lambda\Rightarrow\EOp_{\kappa}\subset\EOp_{\lambda}$
in~$\EOp(r)$
and we can easily check that $\colim_{\kappa\in\K(r)}\EOp_{\kappa} = \EOp(r)$, for every $r\in\NN$.

The action of a permutation $w\in\Sigma_r$ on $\EOp(r)$
maps the subcomplex $\EOp_{\kappa}\subset\EOp(r)$
into $\EOp_{w\kappa}\subset\EOp(r)$
and the composition product of $\EOp$
restricts to morphisms
\begin{equation*}
\EOp_{\kappa}\otimes(\EOp_{\pi_1}\otimes\dots\otimes\EOp_{\pi_r})
\xrightarrow{\mu}\EOp_{\kappa(\pi_1,\dots,\pi_r)}
\end{equation*}
for every $\kappa\in\K(r)$, $\pi_1,\dots,\pi_r\in\K$.
Thus the collection of diagrams $\{\EOp_{\kappa}\}_{\kappa\in\K(r)}$,
inherits the structure of a $\K$-operad
so that the inclusions $i: \EOp_{\kappa}\hookrightarrow\EOp$
yield an isomorphism of operads $i: \colim_{\K}\EOp_{\kappa}\xrightarrow{\simeq}\EOp$.

Recall that we simply replace permutations $w_i\in\Sigma_r$ by bijections $w_i: \{1,\dots,r\}\rightarrow\eset$
in the definition of the Barratt-Eccles operad
to form the dg-module $\EOp(\eset)$ associated to a finite set $\eset$.
The extension of the definition to any finite indexing set works same for the dg-modules~$\EOp_{\kappa}$
when we assume $\kappa\in\K(\eset)$.

\subsubsection{Complete graph posets and cell decompositions of $E_\infty$-operads}\label{CompleteGraphOperad:CellDecompositions}
The complete graph operad $\K$
has a nested sequence of suboperads
\begin{equation*}
\K_1\subset\cdots\subset\K_n\subset\cdots\subset\colim_n\K_n = \K
\end{equation*}
defined by bounding the weight of edges in complete graphs.
Explicitly,
the subposet $\K_n(\eset)\subset\K(\eset)$ consists of complete graphs $\kappa = (\mu,\sigma)\in\K(\eset)$
such that $\mu_{e f}<n$,
for every pair $\{e,f\}\subset\eset$.

For any operad $\POp$ equipped with a $\K$-structure,
we have a sequence of operads
\begin{equation*}
\POp_1\rightarrow\cdots\rightarrow\POp_n\rightarrow\cdots\rightarrow\colim_n\POp_n\xrightarrow{\simeq}\POp
\end{equation*}
such that $\POp_n = \colim_{\K_n}\POp_{\kappa}$.

The main theorem of~\cite{BergerCell}
implies
that the sequence of dg-operads $\POp_n = \colim_{\K_n}\POp_{\kappa}$ is homotopy equivalent to the nested sequence
of the chain operads of little $n$-cubes
\begin{equation*}
C_*(\DOp_1)\rightarrow\cdots\rightarrow C_*(\DOp_n)\rightarrow\cdots\rightarrow\colim_n C_*(\DOp_n) = C_*(\DOp_{\infty}),
\end{equation*}
when we assume:
\begin{enumerate}\renewcommand{\theenumi}{K\arabic{enumi}}\renewcommand{\labelenumi}{(\theenumi)}
\item\label{ReedyKCells}
the collection $\POp_{\kappa}$
forms a cofibrant $\KOp(r)$-diagram in dg-modules
(with respect to the standard model structure of diagrams in a cofibrantly generated model category),
for every $r\in\NN$;
\item\label{ContractibleKCells}
we have a pointwise equivalence of $\K$-operads $\epsilon: \POp_{\kappa}\xrightarrow{\sim}\COp_{\kappa}$,
where $\COp_{\kappa}$ is the constant $\K$-operad defined by the commutative operad $\COp$;
\end{enumerate}
Moreover,
each morphism $\colim_{\K_{n-1}}\POp_{\kappa}\rightarrow\colim_{\K_{n}}\POp_{\kappa}$
is a cofibration of dg-modules
whenever the condition~(\ref{ReedyKCells} holds.
In the sequel,
we say that an operad $\POp$ is a $\K$-cellular $E_\infty$-operad
when $\POp$ is equipped with a $\K$-structure satisfying~(\ref{ReedyKCells}-\ref{ContractibleKCells}).

\subsubsection{The example of the Barratt-Eccles operad (continued)}\label{CompleteGraphOperad:BarrattEcclesCellStructure}
The $\K$-structure $\EOp_{\kappa}$ defined in~\S\ref{CompleteGraphOperad:BarrattEccles}
for the Barratt-Eccles operad satisfies the condition~(\ref{ReedyKCells}-\ref{ContractibleKCells}).
Hence
so that we have a nested sequence of operads $\EOp_n = \colim_{\K_n}\EOp_{\kappa}$, formed from the Barratt-Eccles operad $\EOp$,
and weakly-equivalent to the nested of the chain operads of little $n$-cubes (see~\cite{BergerCell}).

We refer to~\cite{Bar3}
for details about the cofibration condition~(\ref{ReedyKCells})
in the dg-setting.
We just recall the proof of the acyclicity condition~(\ref{ContractibleKCells}).
We arrange the definition of the standard contracting chain homotopy $\nu: \EOp(\eset)\rightarrow\EOp(\eset)$
to prove that each $\EOp_{\kappa}$
is contractible.
Explicitly,
for any simplex $(w_0,\dots,w_d)\in\EOp_{\kappa}$, where $\kappa = (\mu,\sigma)$,
we set:
\begin{equation*}
\nu(w_0,\dots,w_d) = (w_0,\dots,w_d,\sigma).
\end{equation*}
Note that $(w_0,\dots,w_d)\in\EOp_{(\mu,\sigma)}\Rightarrow(w_0,\dots,w_d,\sigma)\in\EOp_{(\mu,\sigma)}$.
Thus we have we well-defined homomorphism $\nu: \EOp_{(\mu,\sigma)}\rightarrow\EOp_{(\mu,\sigma)}$.

Recall that the standard section $\eta: \COp(r)\rightarrow\EOp(r)$
of the augmentation $\epsilon: \EOp(r)\rightarrow\COp(r)$
identifies $\COp(r)$
with the summand of $\EOp(r)$ spanned by the identity permutation $\id\in\Sigma_r$.
Since we have $(\id)\in\EOp_{\kappa}$
for every $\kappa\in\K(r)$,
we immediately obtain a map $\iota: \COp_{\kappa}\rightarrow\EOp_{\kappa}$
such that $\epsilon\iota = \id$.
On the other hand,
we easily check that the modified chain contraction $\nu$
satisfies $\delta(\nu) = \id-\iota\epsilon$
on each $\EOp_{\kappa}$.
Hence we conclude that the augmentation $\epsilon: \EOp\rightarrow\COp$
gives rise to a pointwise equivalence of $\K$-operads $\epsilon: \EOp_{\kappa}\xrightarrow{\sim}\COp_{\kappa}$.

\subsubsection{Modules shaped on complete graph posets}\label{CompleteGraphOperad:KModules}
The definition of a $\K$-structure has an obvious generalization in the context of modules
over an operad:
a right $\K$-module over a $\K$-operad $\ROp$
consists of a collection of $\K(r)$-diagrams $\{M_{\kappa}\}_{\kappa\in\K(r)}$
together with $\Sigma_r$-actions
\begin{equation*}
M_{\kappa}\xrightarrow{w_*} M_{w\kappa},\quad w\in\Sigma_*,
\end{equation*}
and composition products
\begin{equation*}
M_{\kappa}\otimes(\ROp_{\pi_1}\otimes\dots\otimes\ROp_{\pi_r})\xrightarrow{\mu} M_{\kappa(\pi_1,\dots,\pi_r)}
\end{equation*}
which satisfy a natural extension of the standard axioms of modules over operads.

For a right $\K$-module $M$ over a $\K$-operad $\ROp$,
the colimit $M = \colim_{\K} M_{\kappa}$ forms naturally a right module over the operad $\ROp = \colim_{\K}\ROp_{\kappa}$.
In the other direction,
the constant diagrams $M_{\kappa} = M(r)$, $\kappa\in\K(r)$, associated to a right module over $M$
forms a right $\K$-module over $\ROp$
and we have an adjunction bewteen this constant functor and the colimit construction.

Say that a right $\ROp$-module $M$ is equipped with a $\K$-structure
if we have a right $\K$-module $M_{\kappa}$ over the right $\K$-operad $\ROp$
such that $\colim_{\K}\POp_{\kappa}\xrightarrow{\simeq}\POp$
for the adjoint morphism of an embedding of the right $\K$-moduel $M_{\kappa}$
into the constant object defined by $M$.

\section{The structure of iterated bar modules}\label{KStructure}
In the next section,
we prove that the twisting homomorphism of the $n$-fold bar module $\partial_{\epsilon}: T^n\circ\EOp\rightarrow T^n\circ\EOp$
factors through $T^n\circ\EOp_n\subset T^n\circ\EOp$
when $\EOp$ be a $\K$-cell $E_\infty$-operad,
and we conclude from this verification that the $n$-fold bar module $B^n_{\EOp} = (T^n_{\EOp},\partial_{\epsilon})$
is defined over~$\EOp_n$.

To reach this objective,
we study the iterated bar module of the commutative operad first:
we prove that $B^n_{\COp}$ is equipped with a $\K$-structure
over the constant $\K$-operad $\COp$.
For this purpose,
we essentially have to check the definition of the iterated bar complex~$B^n_{\COp} = B^n(\COp)$.

\subsubsection{The bar complex in $\Sigma_*$-modules}\label{KStructure:BarComplexReview}
Let $M$ be any $\Sigma_*$-module.
The definition of the tensor product of $\Sigma_*$-modules in~\S\ref{Background:TensorProduct}
gives an expansion of the form
\begin{equation*}
T^c(\Sigma M)(\eset)
= \bigoplus_{\substack{\eset_1\amalg\dots\amalg\eset_r = \eset\\r\geq 1}}\Sigma M(\eset_1)\otimes\dots\otimes\Sigma M(\eset_r)
\end{equation*}
for the tensor coalgebra~$T^c(\Sigma M)$.
The sum ranges over all integers $r\in\NN^*$
and all partitions $\eset_1\amalg\dots\amalg\eset_r = \eset$.
If we assume that $M$ is connected ($M(0) = 0$),
then the sum ranges over partitions $\eset_1\amalg\dots\amalg\eset_r = \eset$
such that $\eset_i\not=\emptyset$, for $i = 1,\dots,r$.

Suppose $A$ is a commutative algebra in the category of $\Sigma_*$-modules.
The bar differential $\partial: T^c(\Sigma A)\rightarrow T^c(\Sigma A)$
and the shuffle product $\mu: T^c(\Sigma A)\rightarrow T^c(\Sigma A)$
have a transparent representation
in terms of this expansion.

The product of~$A$
is defined by a collection of morphisms $\mu: A(\eset)\otimes A(\fset)\rightarrow A(\eset\amalg\fset)$.
The bar differential has a component
\begin{multline*}
\partial: \Sigma A(\eset_1)\otimes\dots\otimes\Sigma A(\eset_{i})\otimes\Sigma A(\eset_{i+1})\otimes\dots\otimes\Sigma A(\eset_r)\\
\rightarrow\Sigma A(\eset_1)\otimes\dots\otimes\Sigma A(\eset_{i}\amalg\eset_{i+1})\otimes\dots\otimes\Sigma A(\eset_r)
\end{multline*}
induced by the product $\mu: A(\eset_{i})\otimes A(\eset_{i+1})\rightarrow A(\eset_{i}\amalg\eset_{i+1})$,
together with the usual sign,
for each merging $\eset_{i},\eset_{i+1}\mapsto\eset_{i}\amalg\eset_{i+1}$.

The shuffle product $\mu: T^c(\Sigma A)(\eset)\otimes T^c(\Sigma A)(\fset)\rightarrow T^c(\Sigma A)(\eset\amalg\fset)$
has a component
\begin{multline*}
\mu: (\Sigma A(\eset_1)\otimes\dots\otimes\Sigma A(\eset_r))\otimes(\Sigma A(\fset_1)\otimes\dots\otimes\Sigma A(\fset_s))\\
\rightarrow\Sigma A(\gset_1)\otimes\dots\otimes\Sigma A(\gset_{r+s})
\end{multline*}
for each shuffle $\gset_1\amalg\dots\amalg\gset_{r+s} = \eset\amalg\fset$
of partitions $\eset_1\amalg\dots\amalg\eset_r = \eset$ and $\fset_1\amalg\dots\amalg\fset_s = \fset$.

\subsubsection{The $\K$-structure of the iterated bar module: the generating $\Sigma_*$-module}\label{KStructure:InductiveDefinition}
The purpose of the construction of this paragraph and of the next one is to prove that a composite $T^n\circ\POp$,
where $\POp$ is a $\K$-operad,
inherits a $\K$-structure.

To begin with,
we define a $\K(\eset)$-diagram $(T^n)_{\kappa}\subset T^n(\eset)$
such that $T^n(\eset) = \colim_{\kappa\in\K(\eset)} (T^n)_{\kappa}$,
for every finite set $\eset$.

Let $\kappa = (\mu,\sigma)\in\K(\eset)$.
Recall that $T^n$ is defined as the iterated tensor coalgebra $T^n = (T^c\Sigma)^n(I) = T^c(\Sigma T^{n-1})$.
The submodule $(T^n)_{\kappa}\subset T^n(\eset)$
is defined by induction on $n$.
For this purpose,
we use the expansion of the tensor coalgebra
\begin{equation*}
T^n(\eset) = T^c(\Sigma T^{n-1})(\eset)
= \bigoplus_{\substack{\eset_1\amalg\dots\amalg\eset_r = \eset\\r\geq 1}} T^{n-1}(\eset_1)\otimes\dots\otimes T^{n-1}(\eset_r)
\end{equation*}
(we omit marking suspensions to simplify the writing).
A tensor $\xi = x_1\otimes\dots\otimes x_r\in T^{n-1}(\eset_1)\otimes\dots\otimes T^{n-1}(\eset_r)$
belongs to the submodule $(T^n)_{\kappa}\subset T^n$
when:
\begin{enumerate}
\item
each factor $x_i\in T^{n-1}$ satisfies $x_i\in(T^n)_{\kappa|_{\eset_i}}$,
where $\kappa|_{\eset_i}$
is the restriction of the complete graph $\kappa$
to the subset $\eset_i\subset\eset$;
\item
the indices $e,f\in\eset$ such that $\mu_{e f}<n-1$
belong to a same component $\eset_i$ of the partition $\eset = \eset_1\amalg\dots\amalg\eset_r$;
\item
for indices $e,f\in\eset$ in separate components $e\in\eset_i,f\in\eset_j$, $i\not=j$,
and such that $\mu_{e f}=n-1$,
the ordering $\sigma|_{e f}$ of the edge $\{e,f\}$
agrees with the order of~$\{\eset_i,\eset_j\}$ in the decomposition $\eset = \eset_1\amalg\dots\amalg\eset_n$.
\end{enumerate}
We have $\kappa<\pi\Rightarrow(T^n)_{\kappa}\subset(T^n)_{\pi}$.

Note that every element $\xi\in T^n(\eset)$ belongs to a component $(T^n)_{(\mu,\sigma)}$
such that $\max(\mu_{e f})<n$, $\forall e,f$.
Moreover,
we have an equality $T^n(\eset) = (T^n)_{\kappa}$, for $\kappa$ sufficiently large,
from which we deduce the identity $\colim_{\kappa\in\K(\eset)} (T^n)_{\kappa} = T^n(\eset)$.
The next observation follows from similar easy verifications:

\begin{obsv}\label{KStructure:GeneratingObject}
Every element $\xi$ of the generating collection $G^n(r)\subset T^n(r)$
belongs to a component $(T^n)_{(\mu,\sigma)}$
such that $\max(\mu_{e f})<n$ and $\sigma = \id$.

If we have $\xi\in G^n\cap(T^n)_{\kappa}$ for some $\kappa\in\K(r)$,
then there is an element of the form $(\mu,\id)$
such that $(\mu,\id)\leq\kappa$
and $\xi\in G^n\cap(T^n)_{(\mu,\id)}$.
\end{obsv}

The verification of this observation is straightforward from the inductive definition of $G^n(r)\subset T^n(r)$
in proposition~\ref{QuasiFreeLifting:SigmaFreeStructure}
and from the inductive definition of the $\K(r)$-diagrams $(T^n)_{\kappa}$.

\subsubsection{The $\K$-structure of the iterated bar module: the underlying free module}\label{KStructure:CompositeModuleDefinition}
Suppose now that $\POp$ is a $\K$-operad.
For each $m\in\NN$, we have a $\K(m)$-diagram formed by the dg-modules $(T^n\circ\POp)_{\kappa}$
spanned by formal composites $\xi(p_1,\dots,p_r)\in T^n\circ\POp(m)$,
such that
\begin{equation*}
\theta(\pi_1,\dots,\pi_r)\leq\kappa,
\end{equation*}
where $\xi\in(T^n)_{\theta}$ and $p_1\in\POp_{\pi_1},\dots,p_r\in\POp_{\pi_r}$.
One can easily check that the objects $(T^n\circ\POp)_{\kappa}$
inherits the structure of a right $\K$-module over $\POp$
so that $\colim_{\K}(T^n\circ\POp)_{\kappa} = T^n\circ\POp$.
Hence,
we obtain that the free right $\POp$-module $T^n\circ\POp$
admits a $\K$-structure.

\medskip
For the moment,
we only need to apply the definition of the object $(T^n\circ\POp)_{\kappa}$
to the commutative operad $\POp = \COp$,
viewed as a constant $\K$-operad.
In this particular case,
we have:

\begin{prop}\label{KStructure:CommutativeCase}
The module $(T^n\circ\COp)_{(\mu,\sigma)}$
is spanned by composites
$\xi(p_1,\dots,p_r)\in(T^n\circ\COp)(\eset)$, $\xi\in(T^n)_{\kappa}$, $p_1\in\COp(\eset_1),\dots,p_r\in\COp(\eset_r)$,
such that
\begin{equation*}
\kappa((0,\sigma|_{\eset_1}),\dots,(0,\sigma|_{\eset_r}))\leq(\mu,\sigma),
\end{equation*}
where $0$ represents the null collection $0_{e f}\equiv 0$
and $\sigma|_{\eset_i}$ refers to the restriction of the ordering $\sigma$
to the subset $\eset_i$.
\end{prop}

\begin{proof}
By definition,
the object $(T^n\circ\COp)_{(\mu,\sigma)}\subset T^n\circ\COp$
is spanned by formal composites $\xi(p_1,\dots,p_r)$, $\xi\in(T^n)_{\kappa}$, $p_1\in\COp_{\pi_1},\dots,p_r\in\COp_{\pi_r}$
such that $\kappa(\pi_1,\dots,\pi_r)\leq(\mu,\sigma)$.
By observation~\ref{CompleteGraphOperad:CompositeBound},
we have $\kappa(\pi_1,\dots,\pi_r)\leq(\mu,\sigma)$
if and only if $\pi_i\leq(\mu|_{\eset_i},\sigma|_{\eset_i})$, for every $i = 1,\dots,r$,
and
\begin{equation*}
\kappa((\mu|_{\eset_1},\sigma|_{\eset_1}),\dots,(\mu|_{\eset_r},\sigma|_{\eset_r}))\leq(\mu,\sigma).
\end{equation*}
Since we have clearly
\begin{multline*}
(0,\sigma|_{\eset_i})\leq(\mu|_{\eset_i},\sigma|_{\eset_i}),
\quad\COp_{(0,\sigma|_{\eset_i})} = \COp_{(\mu|_{\eset_i},\sigma|_{\eset_i})},
\quad\text{for every $i$},\\
\text{and}\quad\kappa((0,\sigma|_{\eset_1}),\dots,(0,\sigma|_{\eset_r}))
\leq\kappa((\mu|_{\eset_1},\sigma|_{\eset_1}),\dots,(\mu|_{\eset_r},\sigma|_{\eset_r}))\leq(\mu,\sigma),
\end{multline*}
we obtain that $\xi(p_1,\dots,p_r)$
is a composite of the form of the proposition.
\end{proof}

The crucial property which allows us to factor the twisting homomorphism of iterated bar modules to $\EOp_n$-operads
is given by the next proposition:

\begin{prop}\label{KStructure:TwistingCochain}
The twisting homomorphism $\partial_{\gamma}: T^n\circ\COp\rightarrow T^n\circ\COp$
of the iterated bar module $B^n_{\COp}$
satisfies $\partial_{\gamma}((T^n)_{\kappa})\subset(T^n\circ\COp)_{\kappa}$,
for every $\kappa\in\K$.
\end{prop}

\begin{proof}
Recall that $\partial_{\gamma}$ is determined by the bar coderivations $\partial: T^c(\Sigma A)\rightarrow T^c(\Sigma A)$
which occur at each level of the composite
$T^n\circ\COp = (T^c\Sigma\dots T^c\Sigma)(\COp)$.
We prove the implication $\xi\in(T^n)_{\kappa}\Rightarrow\partial_{\gamma}(\xi)\in(T^n\circ\COp)_{\kappa}$
by induction on $n\geq 1$.

For $n = 1$,
the claim is checked by a quick inspection of the formula of~\S\ref{KStructure:BarComplexReview}
applied to $A = \COp$.

Suppose that the assertion of the lemma holds for $n-1$,
where we assume $n>1$.
Let $\xi = x_1\otimes\dots\otimes x_r\in T^{n-1}(\eset_1)\otimes\dots\otimes T^{n-1}(\eset_r)$
be a tensor such that $\xi\in(T^n)_{\kappa}$, for some $\kappa\in\K(\eset)$.
By definition, we have:
\begin{align*}
\partial_{\gamma}(\xi) & = \sum_{i=1}^{r-1}\pm x_1\otimes\dots\otimes\mu(x_i,x_{i+1})\otimes\dots\otimes x_r\\
& + \sum_{i=1}^{r}\pm x_1\otimes\dots\otimes\partial_{\gamma}(x_i)\otimes\dots\otimes x_r,
\end{align*}
where $\mu(x_i,x_{i+1})$ refers to the shuffle product of $x_i$ and $x_{i+1}$ in $T^{n-1} = T^c(\Sigma T^{n-2})$
and $\partial_{\gamma}(x_i)$ is determined by the bar differential $\partial_{\gamma}: T^{n-1}\circ\COp\rightarrow T^{n-1}\circ\COp$
coming from lower iterations of the bar complex.

By induction,
we have $\partial_{\gamma}(x_i)\in(T^{n-1}\circ\COp)_{\kappa|_{\eset_i}}$,
from which we deduce readily $x_1\otimes\dots\otimes\partial_{\gamma}(x_i)\otimes\dots\otimes x_r\in(T^n\circ\COp)_{\kappa}$.

For indices $e,f\in\eset$ such that $e\in\eset_i$ and $f\in\eset_{i+1}$,
we have necessarily $\mu_{e f}\geq n-1$
according to the definition of~\S\ref{KStructure:InductiveDefinition}.
This condition ensures that the shuffle product $\mu(x_i,x_{i+1})\in T^c(\Sigma T^{n-2})(\eset_i\amalg\eset_{i+1})$
belongs to the cell
\begin{equation*}
(T^{n-1})_{\kappa|\eset_i\amalg\eset_{i+1}}\subset T^{n-1}(\eset_i\amalg\eset_{i+1}) = T^c(\Sigma T^{n-2})(\eset_i\amalg\eset_{i+1}),
\end{equation*}
because:
\begin{itemize}
\item
the product $\mu(x_i,x_{i+1})$
consists by definition of shuffles of tensors
\begin{align*}
x_i & = y_1\otimes\dots\otimes y_p\in T^{n-2}(\gset_1)\otimes\dots\otimes T^{n-2}(\gset_p)\\
\text{and}\quad x_{i+1} & = z_1\otimes\dots\otimes z_q\in T^{n-2}(\hset_1)\otimes\dots\otimes T^{n-2}(\hset_q),
\end{align*}
where $\eset_i = \gset_1\amalg\dots\amalg\gset_p$ and $\eset_{i+1} = \hset_1\amalg\dots\amalg\hset_q$;
\item
only pairs of the form $(y_j,z_k)$ can be permuted in $\mu(x_i,x_{i+1})$,
but the condition $\mu_{e f}\geq n-1$ for elements $e\in\gset_j\subset\eset_i$ and $f\in\hset_k\subset\eset_{i+1}$
ensures that the order between these pairs does not matter in the relation
$\mu(x_i,x_{i+1})\in(T^{n-1}\circ\COp)_{\kappa|\eset_i\cup\eset_{i+1}}$.
\end{itemize}

From this observation,
we conclude readily, by an easy inspection of conditions of~\S\ref{KStructure:InductiveDefinition},
that the tensor product $x_1\otimes\dots\otimes\mu(x_i,x_{i+1})\otimes\dots\otimes x_r$
belongs to $(T^n\circ\COp)_{\kappa}$.
Therefore,
we obtain finally $\xi\in(T^n)_{\kappa}\Rightarrow\partial_{\gamma}(\xi)\in(T^n\circ\COp)_{\kappa}$,
for every $n\geq 1$.
\end{proof}

\section{The restriction of iterated bar modules to $E_n$-operads}\label{EnDefinition}
The goal of this section
is to prove that the twisting homomorphism of the iterated bar module $B^n_{\EOp} = (T^n\circ\EOp,\partial_{\epsilon})$
factors through $T^n\circ\EOp_n\subset T^n\circ\EOp$ when $\EOp$ is a $\K$-cellular $E_\infty$-operads.
From this observation,
we conclude that the $n$-fold bar complex $B^n: {}_{\EOp}\E\rightarrow\E$
admits an extension to the category of $\EOp_n$-algebras ${}_{\EOp_n}\E$.

To begin with,
we review briefly the structure of a $\K$-operad $\EOp$,
because we put additional assumptions on $\EOp$
in order to simplify our construction.
The Barratt-Eccles operad still gives an example of a $\K$-cellular $E_\infty$-operad
which fulfils all our requirements.

\subsubsection{Assumptions on $\K$-cellular $E_\infty$-operads}\label{EnDefinition:KCellOperads}
Recall that a $\K$-cellular $E_\infty$-operad consists of an $E_\infty$-operad
equipped with a $\K$-structure
defined by acyclic dg-modules $\EOp_{\kappa}\subset\EOp(r)$, $\kappa\in\K(r)$, $r\in\NN$,
so that the colimit
\begin{multline*}
\EOp_n = \colim_{\K_n}\EOp_\kappa,\\
\text{where}\quad
\K_n(r) = \bigl\{(\mu,\sigma)\in\K(r)\ \text{such that}\ \max(\mu_{i j})<n\bigr\}\subset\K(r),
\end{multline*}
forms an $E_n$-operad.

In the construction of~\S\ref{QuasiFreeLifting},
we assume that $\EOp$ is any $E_\infty$-operad equipped with a section $\iota: \COp(\eset)\rightarrow\EOp(\eset)$
of the augmentation $\epsilon: \EOp(\eset)\rightarrow\COp(\eset)$
and with a contracting chain homotopy $\nu: \EOp(\eset)\rightarrow\EOp(\eset)$
such that $\epsilon\iota = \id$, $\epsilon\nu = 0$ and $\delta(\nu) = \id - \iota\epsilon$,
for every finite set $\eset$ equipped with an ordering $\eset = \{e_1,\dots,e_r\}$.
From now on,
we assume that $\EOp$ is a $\K$-cellular $E_\infty$-operad
together with a section $\iota$ such that $\iota(\COp(r))\subset\EOp_{(0,\id)}$
and a chain-homotopy $\nu: \EOp\rightarrow\EOp$ satisfying $\nu(\EOp_{\kappa})\subset\EOp_{\kappa}$,
for every $\kappa\in\K$ of the form $\kappa = (\mu,\id)$,
where we apply the definition of $\iota$ and $\nu$
to the finite set $\rset = \{1<\dots<r\}$.

The section $\iota$ and the contracting chain-homotopy defined in~\S\ref{IteratedBarModules:BarrattEccles}
for the Barratt-Eccles operad satisfy these conditions (check definitions).

\medskip
Under the assumptions of~\S\ref{EnDefinition:KCellOperads},
we obtain:

\begin{lemm}\label{EnDefinition:LiftingKStructure}
The lifting of~\S\ref{QuasiFreeLifting:BarTwistingCochainLifting}
\begin{equation*}
\xymatrix{ & T^n\circ\EOp\ar[d] \\ G^n\ar@{.>}[ur]^{\epsilon_0}\ar[r]_{\gamma} & T^n\circ\COp }
\end{equation*}
satisfies $\epsilon_0(G^n\cap(T^n)_{\kappa})\subset(T^n\circ\EOp)_{\kappa}$,
for every $\kappa\in\K$ of the form $\kappa = (\mu,\id)$.
\end{lemm}

\begin{proof}
Recall that $\epsilon_0$ is defined by a composite $\epsilon_0 = \tilde{\iota}\cdot\gamma$.

By proposition~\ref{KStructure:TwistingCochain},
we have $\partial_{\gamma}((T^n)_{\kappa})\subset(T^n\circ\COp)_{\kappa}$
for every $\kappa\in\K$.
Thus,
for an element of the form $\kappa = (\mu,\id)$,
the expansion of $\partial_{\gamma}(\xi)$, $\xi\in G^n\cap(T^n)_{\kappa}$,
consists by proposition~\ref{KStructure:CommutativeCase}
of composites $x(p_1,\dots,p_r)\in T^n\circ\COp$,
where $x\in(T^n)_{\theta}$ and $p_1\in\COp(\eset_1),\dots,p_r\in\COp(\eset_r)$,
for some $\theta$ such that $\theta((0,\id),\dots,(0,\id))\leq\kappa$.
By $\Sigma_r$-invariance of composites,
we can moreover assume $x\in G^n\cap(T^n)_{\theta}$.

We have by definition $\tilde{\iota}(x(p_1,\dots,p_r)) = x(\iota(p_1),\dots,\iota(p_r))$
and our assumption on $\iota$ implies $\iota(p_i)\in\EOp_{(0,\id)}$.
Therefore we obtain
\begin{equation*}
\tilde{\iota}(x(p_1,\dots,p_r))
\in(T^n\circ\COp)_{\theta((0,\id),\dots,(0,\id))}
\subset(T^n\circ\COp)_{\kappa}
\end{equation*}
and the lemma follows.
\end{proof}

\begin{lemm}\label{EnDefinition:TwistingCochainKStructure}
The homomorphism $\epsilon_*: G^n\rightarrow T^n\circ\EOp$ of~\S\ref{QuasiFreeLifting:BarTwistingCochainLifting}
satisfies $\epsilon_*(G^n\cap(T^n)_{\kappa})\subset(T^n\circ\EOp)_{\kappa}$,
for every element $\kappa\in\K$ of the form $\kappa = (\mu,\id)$.

The associated twisting homomorphism $\partial_{\epsilon}: T^n\circ\EOp\rightarrow T^n\circ\EOp$
satisfies $\partial_{\epsilon}((T^n\circ\EOp)_{\kappa})\subset(T^n\circ\EOp)_{\kappa}$,
for every $\kappa\in\K$.
\end{lemm}

\begin{proof}
We have by definition $\epsilon_* = \sum_{m=0}^{\infty} \epsilon_m$,
where $\epsilon_m = \sum_{p+q=m-1} \tilde{\nu}(\epsilon_p\dercirc\epsilon_q)$, for $m>0$.
We check that the assertions of the lemma are satisfied by each term $\epsilon_m$, $m\in\NN$,
of $\epsilon_*$.

We have by definition $\epsilon_p\dercirc\epsilon_q(\xi) = \partial_{\epsilon_p}\cdot\partial_{\epsilon_q}(\xi)$,
where $\partial_{\epsilon_p}: T^n\circ\EOp\rightarrow T^n\circ\EOp$
is the homomorphism associated to $\epsilon_p$.
We assume by induction that $\partial_{\epsilon_p}((T^n\circ\EOp)_{\kappa})\subset(T^n\circ\EOp)_{\kappa}$
for every element $\kappa\in\K$ and for every $p<m$.

Suppose $\kappa$ is an element of the form $\kappa = (\mu,\id)$.
For $\xi\in G^n\cap(T^n)_{\kappa}$,
the expansion of $\partial_{\epsilon_p}\cdot\partial_{\epsilon_q}(\xi)$,
consists by observation~\ref{CompleteGraphOperad:CompositeBound}
of composites $x(p_1,\dots,p_r)\in T^n\circ\EOp$,
where $x\in(T^n)_{\theta}$ and $p_1\in\EOp_{\kappa|_{\eset_1}},\dots,p_r\in\EOp_{\kappa|_{\eset_r}}$,
for some $\theta\in\K(r)$ such that $\theta(\kappa|_{\eset_1},\dots,\kappa|_{\eset_r})\leq\kappa$.
By $\Sigma_r$-invariance of composites,
we can moreover assume $x\in G^n\cap(T^n)_{\theta}$.

We have by definition
\begin{equation*}
\tilde{\nu}(x(p_1,\dots,p_r))
= \sum_{i=1}^{r}\pm x(\iota\epsilon(p_1),\dots,\iota\epsilon(p_{i-1}),\nu(p_i),p_{i+1},\dots,p_r).
\end{equation*}
Since $\iota(\COp(r))\subset\EOp_{(0,\id|_{\eset_j})}$,
we have $\iota\epsilon(p_j)\in\EOp_{(0,\id)}\subset\EOp_{\kappa|_{\eset_j}}$,
for $j = 1,\dots,i-1$.
By assumption, we have moreover $p_i\in\EOp_{\kappa|_{\eset_i}}\Rightarrow\nu(p_i)\in\EOp_{\kappa|_{\eset_i}}$.
Hence we obtain
\begin{equation*}
\tilde{\nu}(x(p_1,\dots,p_r))\in(T^n\circ\EOp)_{\theta(\kappa|_{\eset_1},\dots,\kappa|_{\eset_r})}\subset(T^n\circ\EOp)_{\kappa},
\end{equation*}
from which we conclude $\xi\in G^n\cap(T^n)_{\kappa}\Rightarrow\epsilon_m(\xi)\in(T^n\circ\EOp)_{\kappa}$.

The homomorphism $\partial_{\epsilon_m}: T^n\circ\EOp\rightarrow T^n\circ\EOp$
is defined on composites $x(p_1,\dots,p_r)\in T^n\circ\EOp$ such that $x\in G^n(r)$
by $\partial_{\epsilon_m}(x(p_1,\dots,p_r)) = \epsilon_m(x)(p_1,\dots,p_r)$.

Suppose $x\in G^n\cap(T^n)_{\theta}$
and $p_1\in\EOp_{\pi_1},\dots,p_r\in\EOp_{\pi_r}$.
By observation~\ref{KStructure:GeneratingObject},
we can assume that $\theta$ has the form $\theta = (\mu,\id)$.
We have then $\epsilon_m(x)\in(T^n\circ\EOp)_{\theta}$
and
\begin{equation*}
\epsilon_m(x)(p_1,\dots,p_r))\in(T^n\circ\EOp)_{\theta(\pi_1,\dots,\pi_r)}.
\end{equation*}
From this assertion,
we conclude that $\partial_{\epsilon_m}((T^n\circ\EOp)_{\kappa})\subset(T^n\circ\EOp)_{\kappa}$
for every $\kappa\in\K$.

This verification achieves the proof of the lemma.
\end{proof}

\begin{thm}\label{EnDefinition:TwistingCochainRestriction}
In the setting of~\S\ref{EnDefinition:KCellOperads},
the twisting homomorphism $\partial_{\epsilon}: T^n\circ\EOp\rightarrow T^n\circ\EOp$
which arises from the definition of the iterated bar module $B^n_{\EOp} = (T^n\circ\EOp,\partial_{\epsilon})$
in~\S\S\ref{QuasiFreeLifting:TwistingCochainLifting}-\ref{QuasiFreeLifting:QuasiFreeModuleExtension}
satisfies $\partial_{\epsilon}(T^n)\subset T^n\circ\EOp_n$
and admits a restriction to $T^n\circ\EOp_n\subset T^n\circ\EOp$.

Thus we have a quasi-free right $\EOp_n$-module $B^n_{\EOp_n} = (T^n\circ\EOp_n,\partial_{\epsilon})$
defined by the restriction of~$\partial_{\epsilon}$
to $T^n\circ\EOp_n$
and this quasi-free module satisfies the relation $B^n_{\EOp_n}\circ_{\EOp_n}\EOp\simeq B^n_{\EOp}$.
\end{thm}

\begin{proof}
Recall (see~\S\ref{KStructure:InductiveDefinition})
that any element $\xi\in T^n$
belongs to a submodule $(T^n)_{\kappa}\subset T^n$ such that $\kappa\in\K_n$
and $\EOp_n = \colim_{\kappa\in\K_n}\EOp_{\kappa}$.
Therefore lemma~\ref{EnDefinition:TwistingCochainKStructure}
implies $\partial_{\epsilon}(\xi)\in T^n\circ\EOp_n$,
for every $\xi\in T^n$,
and $\partial_{\epsilon}(\xi(p_1,\dots,p_r))\in T^n\circ\EOp_n$
for every composite $\xi(p_1,\dots,p_r)$ such that $p_1,\dots,p_r\in\EOp_n$.

The relation $B^n_{\EOp_n}\circ_{\EOp_n}\EOp\simeq B^n_{\EOp}$
is an immediate consequence of the identity $(T^n\circ\EOp_n)\circ_{\EOp_n}\EOp\simeq T^n\circ\EOp$
for a free module.
\end{proof}

The quasi-free right $\EOp_n$-module $B^n_{\EOp_n} = (T^n\circ\EOp_n,\partial_{\epsilon})$
determines a functor $\Sym_{\EOp_n}(B^n_{\EOp_n}): {}_{\EOp_n}\E\rightarrow\E$
of the form
\begin{equation*}
\Sym_{\EOp_n}(B^n_{\EOp_n},A) = (\Sym_{\EOp_n}(T^n\circ\EOp_n,A),\partial_{\epsilon}) = ((T^c\Sigma)^n(A),\partial_{\epsilon}),
\end{equation*}
for every $\EOp_n$-algebra $A$.
Thus we obtain that $\Sym_{\EOp_n}(B^n_{\EOp_n},A)$
is a twisted dg-module defined by the same underlying functor as the $n$-fold tensor coalgebra $(T^c\Sigma)^n(A)$
together with the twisting homomorphism yielded by the construction
of theorem~\ref{EnDefinition:TwistingCochainRestriction}
at the module level.
The relation $B^n_{\EOp_n}\circ_{\EOp_n}\EOp\simeq B^n_{\EOp}$
implies that the diagram
\begin{equation*}
\xymatrix{ {}_{\EOp_n}\E\ar@{.>}[dr]_{\Sym_{\EOp_n}(B^n_{\EOp_n})} &&
{}_{\EOp}\E\ar@{_{(}->}[]!L-<4pt,0pt>;[ll]\ar[dl]^{\Sym_{\EOp}(B^n_{\EOp}) = B^n(-)} \\
& \E & }
\end{equation*}
commutes.
To summarize:

\begin{thm}\label{EnDefinition:Result}
If the operad $\EOp$ fulfils the requirements of~\S\ref{EnDefinition:KCellOperads},
then the $n$-fold bar complex $B^n: {}_{\EOp}\E\rightarrow\E$
admits an extension to the category of $\EOp_n$-algebras
\begin{equation*}
\xymatrix{ {}_{\EOp_1}\E & \ar@{_{(}->}[]!L-<4pt,0pt>;[l] & \ar@{.}[l] &
{}_{\EOp_n}\E\ar@{_{(}->}[]!L-<4pt,0pt>;[l]\ar@{.>}@/_/[]!D;[drrr]_{B^n} & \ar@{_{(}->}[]!L-<4pt,0pt>;[l] &  \ar@{.}[l] &
\ar@{_{(}->}[]!L-<4pt,0pt>;[l]{}_{\EOp}\E\ar[d]^{B^n} \\
&&&&&& \E }.
\end{equation*}
This extension is defined by a twisted object of the form $B^n(A) = ((T^c\Sigma)^n(A),\partial_{\epsilon})$
for a twisting homomorphism $\partial_{\epsilon}$,
deduced from the result of theorem~\ref{EnDefinition:TwistingCochainRestriction},
lifting the twisting homomorphism of the $n$-fold bar complex of commutative algebras.
\qed
\end{thm}

\part*{Iterated bar complexes and homology theories}\label{HomologyTheories}

The goal of this part is to prove that the $n$-fold desuspension of the $n$-fold bar complex $\Sigma^{-n} B^n(A)$
determines the $\EOp_n$-homology $H^{\EOp_n}_*(A)$,
for every $n\in\NN$,
including $n = \infty$.
The infinite bar complex $\Sigma^{-\infty} B^{\infty}(A)$
is just defined by the colimit of the complexes $\Sigma^{-n} B^n(A)$
over the suspension morphisms $\sigma: \Sigma^{1-n} B^{n-1}(A)\rightarrow\Sigma^{-n} B^{n}(A)$.

The preliminary section (\S\ref{UsualBarComplexes})
is devoted to the computation of the homology of the iterated bar complex of usual commutative algebras,
like trivial algebras and symmetric algebras.
The next section (\S\ref{OperadicHomology})
is devoted to a short review of the definition
of the homology theory $H^{\ROp}_*(-)$
associated to operad $\ROp$.
In the core section (\S\ref{EnHomology}),
we determine the homology of the bar module $\Sigma^{-n} B^n_{\EOp_n}$
and we use the result to prove that $\Sigma^{-n} B^n(A)$
determines the $\EOp_n$-homology $H^{\EOp_n}(A)$ in the case $n<\infty$.
In the last section of this part (\S\ref{EinfinityHomology}),
we prove that $\Sigma^{-n} B^{n}(A)$
determines the $\EOp_n$-homology $H^{\EOp_n}(A)$ in the case $n=\infty$ too.

To obtain our result,
we use that the homology $H^{\ROp}_*(-)$
can be represented by a generalized $\Tor$-functor $\Tor^{\ROp}_*(I,-)$.
The idea is to prove the relation
\begin{equation*}
H_*(\Sigma^{-n} B^n(A)) = \Tor^{\EOp_n}_*(I,A) = H^{\EOp_n}_*(A)
\end{equation*}
by checking that the module $\Sigma^{-n} B^n_{\EOp_n}$, which represents the $n$-fold bar complex $\Sigma^{-n} B^n(A)$,
forms a cofibrant replacement of the composition unit $I$
in the category of right $\EOp_n$-modules.

\subsubsection*{Koszul complexes of operads}
In our arguments (proposition~\ref{UsualBarComplexes:HarrisonKoszul} and proposition~\ref{EnHomology:DtwoTerm}),
we make appear operadic Koszul complexes $K(I,\POp,\POp)$
associated to the commutative and Lie operads $\POp = \COp,\LOp$.
Recall simply that the Koszul complex $K(I,\POp,\POp)$
associated to a Koszul operad $\POp$
is an acyclic complex of right $\POp$-modules.
Let $K_{\POp}(-) = \Sym_{\POp}(K(I,\POp,\POp),-)$ be the functor on $\POp$-algebras determined by this right $\POp$-module $K(I,\POp,\POp)$.
In the case $\POp = \COp$,
the functor $K_{\COp}: {}_{\COp}\E\rightarrow\E$ is identified with the Harrison complex of commutative dg-algebras
(see~\S\ref{UsualBarComplexes:Harrison} for short recollections and~\cite[\S 6]{FressePartitions} for more detailed explanations).
In the case $\POp = \LOp$,
the functor $K_{\LOp}: {}_{\LOp}\E\rightarrow\E$ is identified with the usual Chevalley-Eilenberg complex of Lie dg-algebras
(see~\cite[\S 6]{FressePartitions}).
In this paper,
we use the extension of these complexes to $\POp$-algebras in $\Sigma_*$-modules
and, in order to identify the object $K(I,\POp,\POp)$,
we apply the relation
\begin{equation*}
K(I,\POp,\POp) = \Sym_{\POp}(K(I,\POp,\POp),\POp) = K_{\POp}(\POp)
\end{equation*}
where the operad $\POp$ is viewed as an algebra over itself in the category of $\Sigma_*$-modules
(see explanations of~\S\ref{Background:FunctorOperations}).

Recall that we only deal with generalizations of modules of symmetric tensors.
For us,
the Chevalley-Eilenberg complex is defined by a symmetric algebra $\Sym(\Sigma G)$
on a suspension of Lie algebras $G$,
as in the context of dg-modules over a field of characteristic zero,
but we apply this construction to $\Sigma_*$-modules over any ring (see again~\S\ref{Background:FunctorOperations} for more explanations).
In~\S\ref{UsualBarComplexes:Harrison},
we apply a similar convention for the definition of the Harrison complex.

\subsubsection*{The operadic suspension of $\Sigma_*$-modules}
In this part,
we use a functor $\Lambda: \M\rightarrow\M$
such that:
\begin{equation*}
\Sym(\Lambda M,\Sigma E)\simeq\Sigma\Sym(M,E),
\end{equation*}
for every $\Sigma_*$-module $M\in\M$.
We call the $\Sigma_*$-module $\Lambda M$ the operadic suspension of~$M$.
This $\Sigma_*$-module $\Lambda M$
is defined in arity $r$ by the tensor product:
\begin{equation*}
\Lambda M(r) = \Sigma^{1-r} M(r)\otimes\sgn_r,
\end{equation*}
where $\sgn_r$ refers to the signature representation of $\Sigma_r$ (see~\cite{GetzlerJones}).

The operadic suspension of an operad $\Lambda\POp$ inherits a natural operad structure
and the suspension $\Sigma: E\mapsto\Sigma E$
induces an isomorphism from the category of $\POp$-algebras
to the category of $\Lambda\POp$-algebras (see~\emph{loc. cit.} for details).
Note that the operadic suspension $\Lambda\POp$
has nothing to do with the suspension of the semi-model category of operads.

We have by definition $\Sigma\POp(E) = \Lambda\POp(\Sigma E)$ for a free $\POp$-algebra $A = \POp(E)$.
In the paper,
we often use the equivalent relation $\POp(\Sigma E) = \Sigma\Lambda^{-1}\POp(E)$
which makes appear the operadic desuspension of~$\POp$.

\section{Prelude: iterated bar complexes of usual commutative algebras}\label{UsualBarComplexes}

The homology of the iterated bar complexes of usual commutative algebras (exterior, polynomial, divided power, \dots)
is determined in~\cite{Cartan} in the context of dg-modules.
The purpose of this section is to review these classical homology calculations
in the context of $\Sigma_*$-modules.
To be specific,
we determine the homology $H_*(B^n(A))$ of a free commutative algebra $A = \COp(M)$
and of a trivial algebra $A = M$.

For our needs,
we consider, throughout this section, a $\Sigma_*$-module $M$
such that:
\begin{enumerate}\renewcommand{\theenumi}{M\arabic{enumi}}\renewcommand{\labelenumi}{(\theenumi)}
\item\label{UsualBarComplexes:ConnectednessAssumption}
we have $M(0) = 0$ -- thus $M$ is connected as a $\Sigma_*$-module;
\item\label{UsualBarComplexes:DifferentialAssumption}
the differential of $M$ is zero -- in other words $M$ defines a $\Sigma_*$-objects in the category of graded modules;
\item\label{UsualBarComplexes:ModuleAssumption}
and each component $M(r)$ is projective as a $\kk$-module -- but we do not assume that $M$ is projective as a $\Sigma_*$-module.
\end{enumerate}

\subsubsection{Free commutative algebras in $\Sigma_*$-modules}\label{UsualBarComplexes:FreeCommutativeAlgebras}
The free (non-unital) commutative algebra $\COp(M)$
is identified with the (non-unital) symmetric algebra:
\begin{equation*}
\COp(M) = \bigoplus_{r=1}^{\infty} (\COp(r)\otimes M^{\otimes r})_{\Sigma_r}
= \bigoplus_{r=1}^{\infty} (M^{\otimes r})_{\Sigma_r}.
\end{equation*}
The canonical morphism $\eta: M\rightarrow\COp(M)$
identifies $M$ with the summand $M\subset \bigoplus_{r=1}^{\infty} (M^{\otimes r})_{\Sigma_r}$
of $\COp(M)$.

The element of $\COp(M)$ represented by the tensor $x_1\otimes\dots\otimes x_r\in M^{\otimes r}$
is usually denoted by $x_1\cdot\ldots\cdot x_r$
since this tensor represent the product of $x_1,\dots,x_r\in M$
in $\COp(M)$.

\medskip
For a free commutative algebra $\COp(M)$,
the composite of the canonical morphism $\eta: M\rightarrow\COp(M)$
with the suspension $\sigma: \Sigma\COp(M)\rightarrow B(\COp(M))$
defines a natural morphism of $\Sigma_*$-modules
$\sigma: \Sigma M\rightarrow B(\COp(M))$.
Form the morphism of commutative algebras
\begin{equation*}
\nabla: \COp(\Sigma M)\rightarrow B(\COp(M)).
\end{equation*}
such that $\nabla|_{\Sigma M} = \sigma$.

\begin{prop}\label{UsualBarComplexes:FreeCommutativeAlgebra}
The morphism $\nabla: \COp(\Sigma M)\rightarrow B(\COp(M))$
defines a weak-equivalence of commutative algebras in $\Sigma_*$-modules
whenever the $\Sigma_*$-module $M$
satisfies the requirements~(\ref{UsualBarComplexes:ConnectednessAssumption}-\ref{UsualBarComplexes:ModuleAssumption}).
\end{prop}

The result becomes much more complicated
when $M(0)\not=0$:
the module~$B(\COp(M))$ carries higher homological operations, like divided powers (see~\cite{Cartan}),
and the homology of~$B(\COp(M))$ does not reduce to the free commutative algebra~$\COp(\Sigma M)$.

\begin{proof}
We adapt classical arguments (given for instance in~\cite{Cartan}).

Recall that a non-unital commutative algebra $A$
is equivalent to a unital augmented commutative algebra such that $A_+ = \unit\oplus A$.
In this proof (and in this proof only),
we use a free commutative algebra with unit $\COp_+(M) = \unit\oplus\COp(M)$
and a unital version of the bar complex $B_+(A)$
such that $B_+(A) = \unit\oplus B(A)$.

We consider the acyclic bar complex $B(A_+,A,\unit)$
formed by the tensor products
\begin{equation*}
B(A_+,A,\unit) = A_+\otimes B_+(A) = \bigoplus_{d=0}^{\infty} A_+\otimes\Sigma A^{\otimes d}
\end{equation*}
together with the differential such that
\begin{equation*}
\partial(a_0\otimes a_1\otimes\dots\otimes a_d)
= \pm\sum_{i=0}^{d-1} a_0\otimes\dots\otimes\mu(a_{i},a_{i+1})\otimes\dots\otimes a_d,
\end{equation*}
where $\mu: A\otimes A\rightarrow A$ refers to the product of~$A$
and its extension to $A_+$.

The tensor product of commutative algebras $\COp_+(M)\otimes\COp_+(\Sigma M)$
inherits a natural commutative algebra structure.
Let $\partial: \COp_+(M)\otimes\COp_+(\Sigma M)\rightarrow\COp_+(M)\otimes\COp_+(\Sigma M)$
be the unique derivation of commutative algebras of degree $-1$
which vanishes on $\COp_+(M)$
and extends the canonical homomorphism $\sigma: \Sigma M\rightarrow M$
on $\COp_+(\Sigma M)$.
We have $\partial^2 = 0$ so that the derivation $\partial$
provides the commutative algebra $\COp_+(M)\otimes\COp_+(\Sigma M)$
with a new dg-structure.
The twisted object $K = (\COp_+(M)\otimes\COp_+(\Sigma M),\partial)$
is an analogue in $\Sigma_*$-modules of the usual Koszul complex.

For $r\in\NN$ and any $N\in\M$, we set $\COp_r(N) = (N^{\otimes r})_{\Sigma_r}$.
The differential of the Koszul complex
satisfies $\partial(\COp_p(M)\otimes\COp_q(\Sigma M))\subset\COp_{p+1}(M)\otimes\COp_{q-1}(\Sigma M)$.
We define (non-equivariant) maps $\nu: \COp_p(M)\otimes\COp_q(\Sigma M)\rightarrow\COp_{p-1}(M)\otimes\COp_{q+1}(\Sigma M)$
such that $\delta\nu+\nu\delta = \id$.

We have an identity:
\begin{multline*}
(\COp_p(M)\otimes\COp_q(\Sigma M))(r)\\
= \bigoplus_{\eset_*\amalg\fset_* = \{1,\dots,r\}}
(M(\eset_1)\otimes\dots\otimes M(\eset_p))\otimes(\Sigma M(\fset_1)\otimes\dots\otimes\Sigma M(\fset_q))/\equiv,
\end{multline*}
where the sum is divided out by the action of permutations $(u,v)\in\Sigma_p\times\Sigma_q$.
For a tensor
\begin{equation*}
\xi = (x_1\cdot\ldots\cdot x_p)\otimes(y_1\cdot\ldots\cdot y_q)
\in(M(\eset_1)\otimes\dots\otimes M(\eset_p))\otimes(\Sigma M(\fset_1)\otimes\dots\otimes\Sigma M(\fset_q)),
\end{equation*}
we set
\begin{equation*}
\nu(\xi) = \begin{cases}
\pm(x_1\cdot\ldots\widehat{x_i}\ldots\cdot x_p)\otimes(x_i\cdot y_1\cdot\ldots\cdot y_q),
& \text{if $1\in\eset_i$ for some $i$}, \\
0, & \text{otherwise}.
\end{cases}
\end{equation*}
The relation $\delta\nu+\nu\delta = \id$ follows from an easy verification.
The assumption $M(0)=0$ is used at this point,
because our chain-homotopy does not work when $\{1,\dots,r\}$
is reduced to the emptyset.

The tensor product $\COp_+(M)\otimes\nabla$,
where $\nabla$ is the morphism of the lemma,
defines a morphism of commutative algebras
\begin{multline*}
(\COp_+(M)\otimes\COp_+(\Sigma M),\partial)
\xrightarrow{\COp_+(M)\otimes\nabla}(\COp_+(M)\otimes T^c_+(\Sigma\COp(M)),\partial)\\
= B(\COp_+(M),\COp(M),\unit).
\end{multline*}
Both terms form quasi-free resolutions of the unit object $\unit$
in the category of left $\COp_+(M)$-modules.
By a standard result of homological algebra,
we conclude that the morphism
\begin{equation*}
\COp_+(\Sigma M)\xrightarrow{\nabla}(T^c_+(\Sigma\COp(M)),\partial) = B_+(\COp(M))
\end{equation*}
defines a weak-equivalence.

Use the natural splitting $A_+ = \unit\oplus A$
to get the result of the proposition.
\end{proof}

Proposition~\ref{UsualBarComplexes:FreeCommutativeAlgebra}
implies by an immediate induction:

\begin{prop}\label{UsualBarComplexes:IteratedFreeCommutativeAlgebra}
We have an isomorphism of commutative algebras in $\Sigma_*$-modules
\begin{equation*}
\nabla: \COp(\Sigma^n M)\xrightarrow{\simeq} H_*(B^n(\COp(M))),
\end{equation*}
for every $n\in\NN$,
whenever the $\Sigma_*$-module $M$
satisfies the requirements~(\ref{UsualBarComplexes:ConnectednessAssumption}-\ref{UsualBarComplexes:ModuleAssumption}).
\qed
\end{prop}

\subsubsection{Hopf algebras}\label{UsualBarComplexes:HopfAlgebraStructures}
Recall that the shuffle product commutes with the deconcatenation product of~$T^c(\Sigma M)$,
so that $T^c(\Sigma M)$ forms a commutative Hopf algebra,
for any $\Sigma_*$-module $M$,
and so does the bar complex $B(A) = (T^c(\Sigma A),\partial)$
when $A$ is a commutative algebra (see~\S\ref{IteratedBarModules:HopfStructure}).

The free commutative algebra $\COp(\Sigma M)$
inherits a Hopf algebra structure as well.
The diagonal $\Delta: \COp(\Sigma M)\rightarrow\COp(\Sigma M)\otimes\COp(\Sigma M)$
of a product $x_1\cdot\,\dots\,\cdot x_r\in\COp(\Sigma M)$
can be defined by the explicit formula
\begin{multline*}
\Delta(x_1\cdot\,\dots\,\cdot x_r) \\
= \sum_{\substack{p+q=r\\p,q>0}}\Bigl\{\sum_{w\in\Sh(p,q)} \pm x_{w(1)}\cdot\,\dots\,\cdot x_{w(p)}
\cdot\varpi(x_{w(p+1)}\cdot\,\dots\,\cdot x_{w(p+q)})\Bigr\},
\end{multline*}
where the inner sum ranges over the set of $(p,q)$-shuffles.
Note that we assume $p,q>0$
since we have removed units from the commutative algebra $\COp(\Sigma M)$.
The generators $\xi\in\Sigma M$
are primitive in the sense that $\Delta(\xi) = 0$ (since we do not consider unit).

The result of proposition~\ref{UsualBarComplexes:FreeCommutativeAlgebra}
can be improved to:

\begin{prop}\label{UsualBarComplexes:HopfAlgebraIso}
The morphism of proposition~\ref{UsualBarComplexes:FreeCommutativeAlgebra}
\begin{equation*}
\nabla: \COp(\Sigma M)\rightarrow B(\COp(M))
\end{equation*}
commutes with diagonals and yields an isomorphism
\begin{equation*}
\nabla_*: \COp(\Sigma M)\xrightarrow{\simeq} H_*(B(\COp(M)))
\end{equation*}
in the category of Hopf algebras in $\Sigma_*$-modules.
\end{prop}

\begin{proof}
The image of an element $\xi\in\Sigma M$
under the suspension $\sigma: \Sigma M\rightarrow B(\COp(M))$
defines clearly a primitive element in $B(\COp(M))$.
Hence $\nabla$ preserves the diagonal of generators of the commutative algebra $\COp(\Sigma M)$.
We conclude readily that $\nabla$ preserves the diagonal of any element of $\COp(\Sigma M)$
by using the commutation relation between products and coproducts
in Hopf algebras (without unit).
\end{proof}

The classical Milnor-Moore and Poincar\'e-Birkhoff-Witt theorems,
which give the structure of cocommutative Hopf algebras,
have a natural generalization in the context of $\Sigma_*$-modules (see~\cite{Stover})
and so do the dual statements
which apply to commutative Hopf algebras.
If we restrict ourself to $\Sigma_*$-objects such that $M(0) = 0$,
then the generalized Milnor-Moore and Poincar\'e-Birkhoff-Witt theorems hold over a ring,
unlike the classical statements.

For the Hopf algebra $T^c(\Sigma M)$,
we obtain:

\begin{fact}\label{UsualBarComplexes:Coenveloping}
The Hopf algebra $T^c(\Sigma M)$ is identified with the coenveloping coalgebra of the cofree Lie coalgebra $L^c(\Sigma M)$.
\end{fact}

This assertion is a direct consequence of the adjoint definition of coenveloping coalgebras
(see~\cite[\S 4.2]{FressePoisson},
see also~\cite[Chapitre II, \S 3]{Bourbaki} or \cite[\S\S 0.1-0.2]{Reutenauer}
for the usual dual statement about enveloping algebras of free Lie algebras).

Then:

\begin{fact}\label{UsualBarComplexes:MilnorMoore}\hspace*{2mm}
\begin{enumerate}
\item
The indecomposable quotient of $T^c(\Sigma M)$
under the shuffle product
is isomorphic to the cofree Lie coalgebra $L^c(\Sigma M)$.
\item
The Hopf algebra $T^c(\Sigma M)$ comes equipped with a fitration such that
\begin{equation*}
\gr_1 T^c(\Sigma M)\simeq L^c(\Sigma M)
\end{equation*}
and we have an isomorphism of graded commutative Hopf algebras in $\Sigma_*$-modules
\begin{equation*}
\COp(L^c(\Sigma M))\simeq\gr_* T^c(\Sigma M).
\end{equation*}
\end{enumerate}
\end{fact}

These assertions follow from the dual version of the generalized Milnor-Moore and Poincar\'e-Birkhoff-Witt theorems of~\cite{Stover}.
In their usual form,
these theorems are stated for Hopf algebras with units.
Again,
we can simply take augmentation ideals to obtain the objects
required by our non-unital setting.

Note that the distribution relation $F(M\circ P) = F(M)\circ P$
holds for $F(-) = L^c(\Sigma -)$
since the cofree Lie coalgebra $L^c(\Sigma M)$ is identified with a quotient of $T^c(\Sigma M)$.
In fact, the functor $L^c(\Sigma M)$ can be identified with a composite $\LOp^{\vee}\circ\Sigma M$,
where $\LOp^{\vee}$ is the $\kk$-dual of the Lie operad
(see the digression and the discussion about coalgebra structures in~\cite[\S\S 1.2.12-1.2.19]{FressePartitions}).
From this observation we also deduce that $L^c(\Sigma M)$
satisfies the K\"unneth isomorphism theorem of~\S\ref{Background:Kunneth}.

\subsubsection{The Harrison complex}\label{UsualBarComplexes:Harrison}
The identity $\Indec T^c(\Sigma M) = L^c(\Sigma M)$
implies that the indecomposable quotient of the bar complex of a commutative algebra
is given by a twisted $\Sigma_*$-module of the form:
\begin{equation*}
\Indec B(A) = (\Indec T^c(\Sigma A),\partial) = (L^c(\Sigma A),\partial).
\end{equation*}
The chain complex $(L^c(\Sigma A),\partial)$ is a generalization,
in the context of $\Sigma_*$-modules,
of the standard Harrison complex with trivial coefficients (see~\cite{Harrison}).

The next proposition is classical for the standard Harrison complex of a free commutative algebra
over a field of characteristic zero.

\begin{prop}\label{UsualBarComplexes:HarrisonKoszul}
The Harrison complex of the free commutative algebra $\COp(M)$ is acyclic
and we have $H_*(L^c(\Sigma\COp(M)),\partial) = \Sigma M$
whenever $M$ satisfies the requirements~(\ref{UsualBarComplexes:ConnectednessAssumption}-\ref{UsualBarComplexes:ModuleAssumption}).
\end{prop}

This result does not hold for a usual commutative algebra in dg-modules
if the ground ring is not a field
or is a field of positive characteristic.
To make the result hold in this setting,
we really need the assumption $M(0) = 0$
and the K\"unneth isomorphism of~\S\ref{Background:Kunneth}.

\begin{proof}
The complex $(L^c(\Sigma\COp(M)),\partial)$
can be identified with a composite $\Sigma_*$-module
\begin{equation*}
(L^c(\Sigma\COp(M)),\partial) = \Sigma K(I,\COp,\COp)\circ M,
\end{equation*}
where $K(I,\COp,\COp)$ is the Koszul complex of the commutative operad
(see~\cite[\S 6.2, \S 6.6]{FressePartitions}).
This complex $K(I,\COp,\COp)$
is acyclic because the commutative operad is Koszul (see \emph{loc. cit.}).
The results of~\cite[\S 2.3]{FressePartitions}
imply moreover that the weak-equivalence $K(I,\COp,\COp)\xrightarrow{\sim} I$
induces a weak-equivalence
\begin{equation*}
\Sigma K(I,\COp,\COp)\circ M\xrightarrow{\sim}\Sigma I\circ M = \Sigma M
\end{equation*}
under the assumptions~(\ref{UsualBarComplexes:ConnectednessAssumption}-\ref{UsualBarComplexes:ModuleAssumption}).
The proposition follows.
\end{proof}

We determine now the homology of the iterated bar complexes $B^n(M)$,
where $M$ is a $\Sigma_*$-module equipped with a trivial commutative algebra structure.
For $n = 1$,
the complex $B(M)$ has a trivial differential.
Therefore:

\begin{fact}\label{UsualBarComplexes:TrivialAlgebra}
We have an identity of commutative algebras in $\Sigma_*$-module $B(M) = T^c(\Sigma M)$,
where $T^c(\Sigma M)$
is equipped with the shuffle product of tensors.
\end{fact}

Thus we study the bar complex of a commutative algebra of the form $T^c(\Sigma M)$, $M\in\M$.

\begin{prop}\label{UsualBarComplexes:ShuffleAlgebra}
The morphism $\sigma_*: \Sigma T^c(\Sigma M)\rightarrow H_*(B(T^c(\Sigma M)))$
induced by the suspension admits a factorization:
\begin{equation*}
\xymatrix{ \Sigma T^c(\Sigma M)\ar[rr]^{\sigma_*}\ar[dr] && H_*(B(T^c(\Sigma M))) \\
& \Sigma L^c(\Sigma M)\ar@{.>}[ur]_{\overline{\sigma}_*} & },
\end{equation*}
and the morphism of commutative algebras in $\Sigma_*$-modules
associated to $\overline{\sigma}_*$
defines an isomorphism
\begin{equation*}
\nabla: \COp(\Sigma L^c(\Sigma M))\xrightarrow{\simeq} H_*(B(T^c(\Sigma M)))
\end{equation*}
whenever the requirements~(\ref{UsualBarComplexes:ConnectednessAssumption}-\ref{UsualBarComplexes:ModuleAssumption})
are fulfilled.
\end{prop}

\begin{proof}
The first assertion of the proposition follows from an observation of~\S\ref{IteratedBarModules:Indecomposables}
and the identity $\Indec T^c(\Sigma M) = L^c(\Sigma M)$.

The filtration of the Hopf algebra $T^c(\Sigma M)$
gives rise to a spectral sequence of commutative algebras $E^r\Rightarrow H_*(B(T^c(\Sigma M)))$
such that
\begin{equation*}
(E^0,d^0) = B(\gr_* T^c(\Sigma M))).
\end{equation*}
Since $\gr_* T^c(\Sigma M)\simeq\COp(L^c(\Sigma M))$,
we have $E^1 = \COp(\Sigma L^c(\Sigma M))$
by proposition~\ref{UsualBarComplexes:FreeCommutativeAlgebra}.
By a straightforward inspection of the construction,
we check that the morphism
\begin{equation*}
\Sigma L^c(\Sigma M)\rightarrow H_*(B(T^c(\Sigma M)))
\end{equation*}
restricts to the isomorphism $\Sigma L^c(\Sigma M)\xrightarrow{\simeq} E^1_1$
on the $E^1$-term of the spectral sequence.

This observation implies that all differentials of the spectral sequence vanish
since $E^1 = \COp(\Sigma L^c(\Sigma M))\Rightarrow H_*(B(T^c(\Sigma M)))$
forms a spectral sequence of commutative algebras.

Hence the spectral sequence degenerates at the $E^1$-level
and we conclude that our morphism
$\Sigma L^c(\Sigma M)\rightarrow H_*(B(T^c(\Sigma M)))$
gives rise to an isomorphism of commutative algebras
$\COp(\Sigma L^c(\Sigma M))\xrightarrow{\simeq} H_*(B(T^c(\Sigma M)))$.
\end{proof}

\begin{prop}\label{UsualBarComplexes:IteratedShuffleAlgebra}
For every $n>1$,
the natural morphism
\begin{equation*}
\overline{\sigma}^{n-1}_*: \Sigma^{n-1} L^c(\Sigma M)\rightarrow H_*(B^{n-1}(T^c(\Sigma M)))
\end{equation*}
which arises from the $(n-1)$-fold suspension
\begin{equation*}
\xymatrix{ \Sigma^{n-1} T^c(\Sigma M)\ar@/^2em/[rrr]^{\sigma^{n-1}_*}\ar[r]\ar[d] &
\Sigma^{n-2} H_*(B(T^c(\Sigma M)))\ar[r] & \dots\ar[r] & H_*(B^{n-1}(T^c(\Sigma M))) \\
\Sigma^{n-1} L^c(\Sigma M)\ar@{.>}[ur] &&& }
\end{equation*}
induces an isomorphism of commutative algebras in $\Sigma_*$-modules:
\begin{equation*}
\COp(\Sigma^{n-1} L^c(\Sigma M))\xrightarrow{\simeq} H_*(B^{n-1}(T^c(\Sigma M)))
\end{equation*}
whenever the requirements~(\ref{UsualBarComplexes:ConnectednessAssumption}-\ref{UsualBarComplexes:ModuleAssumption})
are fulfilled.
\end{prop}

\begin{proof}
For any commutative dg-algebra $A$,
we have a natural spectral sequence of commutative algebras $E^r\Rightarrow H_*(B(A))$
such that $(E^1,d^1) = B(H_*(A))$.
In the case $A = B^{n-2}(T^c(\Sigma M))$,
we have by induction $H_*(A)\simeq\COp(\Sigma^{n-2} L^c(\Sigma M))$
and proposition~\ref{UsualBarComplexes:FreeCommutativeAlgebra}
implies $E^2\simeq\COp(\Sigma^{n-1} L^c(\Sigma M))$.

By a straightforward inspection of the construction,
we check that the composite of the proposition
restricts to the isomorphism $\Sigma^{n-1} L^c(\Sigma M)\xrightarrow{\simeq} E^2_1$
on the $E^2$-term of the spectral sequence.

This observation implies again that all differentials
of the spectral sequence vanish
since $(E^1,d^1) = B(H_*(A))\Rightarrow H_*(B(A))$
forms a spectral sequence of commutative algebras.
The conclusion follows readily.
\end{proof}

\section[Homology of algebras over operads]{Homology of algebras over operads\\and\\operadic $\Tor$-functors}\label{OperadicHomology}
In this section,
we review the definition of the homology theory $H^{\ROp}_*$
associated to an operad $\ROp$.
Usually,
the homology module $H^{\ROp}_*$
is defined by a derived functor of indecomposables $L\Indec: \Ho({}_{\ROp}\E)\rightarrow\Ho(\E)$.
We apply theorems of~\cite[\S 15]{Bar0}
to observe that $H^{\ROp}_*$ is represented by a generalized $\Tor$-functor on $\ROp$-modules $\Tor^{\ROp}_*(I,-)$,
where $I$ is the composition unit of the category of $\Sigma_*$-modules.
We use this representation in the next sections in order to prove that iterated bar complexes
determine the homology theories associated to $E_n$-operads.

\subsubsection{Augmented operads}\label{EnHomology:OperadAugmentation}
The homology $H^{\ROp}_*$ is defined for certain operads $\ROp$
equipped with an augmentation over the composition unit of $\Sigma_*$-modules $I$.

The unit relation gives an isomorphism $I\circ I\simeq I$
which provides $I$ with an obvious operad structure.
The category of algebras ${}_{I}\E$ associated to this operad
is identitified with the underlying category $\E$.
Let $\ROp$ be an operad equipped with an augmentation $\epsilon: \ROp\rightarrow I$.
The restriction functor
\begin{equation*}
\EOp = {}_{I}\E\xrightarrow{\epsilon^*}{}_{\ROp}\E
\end{equation*}
identifies an object $E\in\E$ with an $\ROp$-algebra equipped with a trivial structure.
The functor of indecomposables $\Indec: {}_{\ROp}\E\rightarrow\E$
represents the left adjoint of this category embedding $\E\hookrightarrow{}_{\ROp}\E$
and can be identified with the extension functor
\begin{equation*}
{}_{\ROp}\E\xrightarrow{\epsilon_!}{}_{I}\E = \E
\end{equation*}
associated to the augmentation $\epsilon: \ROp\rightarrow I$.

The associative operad $\AOp$ and the commutative operad $\COp$
are canonically augmented over $I$
in such a way that the diagram
\begin{equation*}
\xymatrix{ \AOp\ar[rr]\ar[dr] && \COp\ar[dl] \\ & I & }
\end{equation*}
commutes.
In these examples $\ROp = \AOp,\COp$,
the indecomposable quotient of an $\ROp$-algebra $A$
is identified with the cokernel of the product $\mu: A\otimes A\rightarrow A$.
Thus we retrieve the usual definition of~\S\ref{IteratedBarModules:Indecomposables}.

The operads occurring in a nested sequence
\begin{equation*}
\EOp_1\rightarrow\dots\rightarrow\EOp_n\rightarrow\dots\rightarrow\colim_n\EOp_n = \EOp
\end{equation*}
inherit a canonical augmentation $\EOp_n\rightarrow I$,
since the $E_\infty$-operad $\EOp$ is supposed to be augmented over $\COp$.

\subsubsection{Homology of algebras over operads (standard definition)}\label{EnHomology:OperadHomology}
The category of $\ROp$-algebras ${}_{\ROp}\E$ inherits a semi-model structure
if the operad $\ROp$ is $\Sigma_*$-cofibrant.
The functors $\Indec = \epsilon_!: {}_{\ROp}\EOp\rightarrow\EOp$ and $\epsilon^*: \EOp\rightarrow{}_{\ROp}\EOp$
forms a Quillen pair and as such determine a pair of adjoint derived functors
\begin{equation*}
L\Indec: \Ho({}_{\ROp}\E)\rightleftarrows\Ho(\E) :\epsilon^*.
\end{equation*}
The derived functor of indecomposables maps an $\ROp$-algebra $A$ to the indecomposable quotient $\Indec Q_A$
of a cofibrant replacement $0\rightarrowtail Q_A\xrightarrow{\sim} A$
in ${}_{\ROp}\E$.

The homology of $A$ is defined by setting $H^{\ROp}_*(A) = H_*(\Indec Q_A)$.
We refer to~\cite[\S 13, \S 16]{Bar0}
for a more comprehensive account on this background.

\subsubsection{Homology of algebras over operads and generalized $\Tor$-functors}\label{EnHomology:TorDefinition}
We define the generalized $\Tor$-functor $\Tor^{\ROp}_*(M): {}_{\ROp}\E\rightarrow\E$
associated to a right $\ROp$-module $M$
by the homology of the functor $\Sym_{\ROp}(P_M): {}_{\ROp}\E\rightarrow\E$
associated to any cofibrant replacement $0\rightarrowtail P_M\xrightarrow{\sim} M$
in $\M{}_{\ROp}$.
Explicitly, we set:
\begin{equation*}
\Tor^{\ROp}_*(M,A) = H_*(\Sym_{\ROp}(P_M,A)).
\end{equation*}
This definition makes sense for any $\E$-cofibrant operad $\ROp$,
for any $\E$-cofibrant right $\ROp$-module $M$
and for every $\E$-cofibrant $\ROp$-algebra $A$:
the assertions of~\cite[Theorems 15.1.A]{Bar0}
imply that the homotopy type of $\Sym_{\ROp}(P_M,A)$
does not depend on the choice of the cofibrant replacement $0\rightarrowtail P_M\xrightarrow{\sim} M$;
moreover,
the map $\Tor^{\ROp}_*: (M,A)\mapsto\Tor^{\ROp}_*(M,A)$
defines a bifunctor which satisfies reasonable homotopy invariance properties
with respect to $(M,A)$.

Pick a cofibrant replacement $0\rightarrowtail Q_A\xrightarrow{\sim} A$
of~$A$ in $\E{}_{\ROp}$,
assuming that the operad $\ROp$ is $\Sigma_*$-cofibrant.
The weak-equivalences $P_M\xrightarrow{\sim} M$ and $Q_A\xrightarrow{\sim} A$
induce morphisms at the functor level
and we have a commutative diagram:
\begin{equation*}
\xymatrix{ \Sym_{\ROp}(P_M,Q_A)\ar[r]\ar[d] & \Sym_{\ROp}(P_M,A)\ar[d] \\
\Sym_{\ROp}(M,Q_A)\ar[r] & \Sym_{\ROp}(M,A) }.
\end{equation*}
If $M$ is also $\Sigma_*$-cofibrant and $A$ is $\E$-cofibrant,
then the left-hand vertical morphism and the upper horizontal morphism of the diagram
are weak-equivalences (see~\cite[Theorems 15.1.A-15.2.A]{Bar0}).
Thus, for a $\Sigma_*$-cofibrant module $M$,
we have a natural isomorphism:
\begin{equation*}
\Tor^{\ROp}_*(M,A) = H_*(\Sym_{\ROp}(P_M,A))\simeq H_*(\Sym_{\ROp}(M,Q_A)),
\end{equation*}
for every $\E$-cofibrant $\ROp$-algebra $A$.

By~\cite[Theorem 7.2.2]{Bar0},
the extension functor $\epsilon_!: {}_{\ROp}\EOp\rightarrow\EOp$
is identified with the functor $\epsilon_! = \Sym_{\ROp}(I)$,
where we use the augmentation $\epsilon: \ROp\rightarrow I$
to provide the composition unit $I$ with a right $\ROp$-module structure.
Since the unit object $I$ is obviously $\Sigma_*$-cofibrant,
we have an identity
\begin{equation*}
H_*(\Sym_{\ROp}(P_I,A)) = H_*(\Sym_{\ROp}(M,Q_A)),
\end{equation*}
for any $\E$-cofibrant $\ROp$-algebra $A$,
from which we deduce the relation:
\begin{equation*}
\Tor^{\ROp}_*(I,A) = H^{\ROp}_*(A).
\end{equation*}

\section[Homology of algebras over $E_n$-operads]{Iterated bar complexes\\and\\homology of algebras over $E_n$-operads}\label{EnHomology}
The goal of this section is to prove that the homology of the category of algebras over an $E_n$-operad~$\EOp_n$
is determined by the $n$-fold bar complex $B^n(A)$.
For this purpose,
we check that the $n$-fold bar module $B^n_{\EOp_n}$
defines a cofibrant replacement of $I$ in the category of right $\EOp_n$-modules
and we apply the interpretation of $H^{\EOp_n}_*$ in terms of operadic $\Tor$-functors.

The module $\Sigma^{-n} B^n_{\EOp_n}$ is cofibrant
by construction (see~\S\ref{QuasiFreeLifting:QuasiFreeModules}).
Our main task is to prove that it is acyclic.

We use a spectral sequence to reduce the problem
to the acyclicity of a chain complex of the form $E^1 = H_*(B^n(I))\circ H_*(\EOp_n)$,
where $B^n(I)$ is the iterated bar complex of the composition unit $I\in\M$,
viewed as a commutative algebra equipped with a trivial structure.
We focus on the case $n>1$.
The homology operad $H_*(\EOp_n)$ is determined in~\cite{Cohen}
and has a nice description as a composite of the commutative operad and a desuspension of the Lie operad.
This composite is usually called the Gerstenhaber operad
and is denoted by~$\GOp_n$.
The Gerstenhaber operad is Koszul (see~\cite{Ginot,Markl})
and this property gives the deep reason for the acyclicity
of the iterated bar module~$B^n_{\EOp_n}$.
For technical reasons,
we split the proof of the acyclicity of $E^1 = H_*(B^n(I))\circ H_*(\EOp_n)$
in two steps
and we rather use that the commutative operad and the Lie operad,
which fit into the decomposition of the Gerstenhaber operad,
are both Koszul operads.

\subsubsection{The augmentation morphisms}\label{EnHomology:IteratedBarModuleApplication}
We already know that $B^n_{\EOp_n}$ is cofibrant.
We also have a natural augmentation $\epsilon: \Sigma^{-n} B^n_{\EOp_n}\rightarrow I$
which is defined as follows.
For any commutative algebra $A$,
the natural morphisms
\begin{equation*}
T^c(\Sigma A) = \bigoplus_{d=1}^{\infty} (\Sigma A)^{\otimes d}\rightarrow\Sigma A\rightarrow\Sigma\Indec A
\end{equation*}
give a morphism of commutative algebras $\epsilon: B(A)\rightarrow\Sigma\Indec A$,
where the object $\Sigma\Indec A$ is identified with a commutative algebra
equipped with a trivial structure.
For the composition unit $I$, viewed as a trivial commutative algebra in $\Sigma_*$-modules,
we have $\Indec I = I$.
Hence, in this case,
we have a morphism of commutative algebras $\epsilon: B(I)\rightarrow\Sigma I$.
More generally,
we have a morphism of commutative algebras $\epsilon: B(\Sigma^{n-1} I)\rightarrow\Sigma^n I$,
for any suspension $\Sigma^{n-1} I$,
and we can iterate the construction to obtain a morphism $\epsilon: B^n(I)\rightarrow\Sigma^n I$
on the $n$-fold bar complex~$B^n(I)$.

The augmentation of an $E_n$-operad $\EOp_n\rightarrow I$
defines a morphism of $\EOp_n$-algebras in right $\EOp_n$-modules
and, by functoriality, gives rise to a morphism
at the level of the $n$-fold bar complex.
Note that we can view the composition unit $I$
as a commutative algebra in right $\EOp_n$-modules in the previous construction,
so that the augmentation $\epsilon: B^n(I)\rightarrow\Sigma^{n} I$
defines a morphism of right $\EOp_n$-modules and not only of $\Sigma_*$-modules.
The augmentation $\epsilon: \Sigma^{-n} B^n_{\EOp_n}\rightarrow I$
is defined by a desuspension of the composite:
\begin{equation*}
B^n_{\EOp_n} = B^n(\EOp_n)\rightarrow B^n(I)\rightarrow\Sigma^n I.
\end{equation*}

It remains to determine the homology of $B^n_{\EOp_n}$
in order to prove that the morphism $\epsilon: B^n_{\EOp_n}\rightarrow\Sigma^n I$
(or its desuspension) forms a weak-equivalence.
In the remainder of this section,
we focus on the case $n>1$,
because for $n = 1$ we obtain immediately:

\begin{prop}\label{EnHomology:EoneCase}
In the case $n = 1$,
we have weak-equivalences $B_{\EOp_1} = B(\EOp_1)\xrightarrow{\sim} B(\AOp)\xrightarrow{\sim}\Sigma I$.
\end{prop}

\begin{proof}
Recall that the operad $\EOp_1$
forms an $A_\infty$-operad
and is connected to the associative operad $\AOp$
by a weak-equivalence $\epsilon: \EOp_1\xrightarrow{\sim}\AOp$.
This augmentation induces a weak-equivalence
at the bar complex level:
\begin{equation*}
B_{\EOp_1} = B(\EOp_1)\xrightarrow{\sim} B(\AOp).
\end{equation*}

The bar complex $B(\AOp)$ is identified with the suspension of the Koszul complex~$K(I,\AOp,\AOp)$
of the associative operad~$\AOp$ (see~\cite[\S 5.2]{FressePartitions} and~\cite{GinzburgKapranov}).
This complex is acyclic,
because the associative operad is an instance of a Koszul operad (see~\cite{GinzburgKapranov}).
\end{proof}

The spectral sequence, used to reduce the calculation of $H_*(B^n_{\EOp_n})$ for $n>1$,
arises from the definition of $B^n_{\EOp_n}$
as a quasi-free module over the operad $\EOp_n$.

\subsubsection{The natural spectral sequence of a quasi-free module}\label{EnHomology:QuasiFreeFiltration}
Let $M = (K\circ\ROp,\partial_{\alpha})$
be a quasi-free module over an operad $\ROp$
such that $\ROp(0) = 0$.
Let $F_0 K\subset\dots\subset F_s K\subset\dots\subset K$
be the filtration of $K$
formed by the $\Sigma_*$-modules such that:
\begin{equation*}
F_s K(r) = \begin{cases} K(r), & \text{if $s\leq r$}, \\ 0, & \text{otherwise}. \end{cases}
\end{equation*}

\begin{claim}
If $\ROp(0) = 0$, then the twisting homomorphism $\partial_{\alpha}: K\circ\ROp\rightarrow K\circ\ROp$
satisfies $\partial_{\alpha}(F_s K\circ\ROp)\subset F_s K\circ\ROp$,
for every $s\in\NN$.
\end{claim}

\begin{proof}
Observe that $\xi\in K\circ\ROp(r)\Rightarrow\xi\in F_r K\circ\ROp(r)$.
Indeed,
for a composite $y(q_1,\dots,q_s)\in K\circ\ROp(r)$,
where $y\in K(s)$ and $q_1\in\ROp(n_1),\dots,q_s\in\ROp(n_s)$,
we have $n_1+\dots+n_s = r$.
Since $\ROp(0) = 0$,
we have necessarily $n_i>0$, for $i = 1,\dots,s$,
from which we deduce $s\leq n_1+\dots+n_s = r$.
Hence,
we obtain $y(q_1,\dots,q_s)\in F_r K\circ\ROp(r)$.

For a generating element $x\in K(r)$,
we have $\partial_{\alpha}(x)\in(K\circ\ROp)(r)\Rightarrow\partial_{\alpha}(x)\in(F_r K\circ\ROp)(r)$.
For a composite $x(p_1,\dots,p_r)\in K\circ\ROp$,
we still have $\partial_{\alpha}(x(p_1,\dots,p_r))\in F_r K\circ\ROp$
since $\partial_{\alpha}(x(p_1,\dots,p_r)) = \partial_{\alpha}(x)(p_1,\dots,p_r)$.
Therefore
we conclude $\xi\in F_s K\circ\ROp\Rightarrow\partial_{\alpha}(\xi)\in F_s K\circ\ROp$.
\end{proof}

The relation $\partial_{\alpha}(F_s K\circ\ROp)\subset F_s K\circ\ROp$
implies that the quasi-free module $M = (K\circ\ROp,\partial_{\alpha})$
has a filtration by submodules
such that:
\begin{equation*}
F_s M = (F_s K\circ\ROp,\partial_{\alpha}).
\end{equation*}
The spectral sequence $E^r(M)\Rightarrow H_*(M)$
associated to a quasi-free module $M$
is the spectral sequence defined by this filtration $F_0 M\subset\dots\subset F_s M\subset\dots\subset M$.

\medskip
We have the easy observation:

\begin{obsv}\label{EnHomology:QuasiFreeModuleSpectralSequence}
Each term of the spectral sequence $E^r(M)$ inherits a natural right $H_*(\ROp)$-action.
This right $H_*(\ROp)$-action is preserved by the differential $d^r: E^r(M)\rightarrow E^r(M)$
so that $E^r(M)\Rightarrow H_*(M)$
defines a spectral sequence of right $H_*(\ROp)$-modules.
\end{obsv}

\subsubsection{Functoriality of the spectral sequence}\label{EnHomology:FiltrationFunctoriality}
Let $\phi: M\rightarrow N$
be any morphism between quasi-free modules $M = (K\circ\ROp,\partial)$ and $N = (L\circ\ROp,\partial)$.

For a generating element $x\in K(r)$, we have $\phi(x)\in K\circ\ROp(r)\Rightarrow\phi(x)\in F_r K\circ\ROp(r)$.
From this assertion we deduce immediately that $\phi$ preserves the filtrations
of the spectral sequence of~\S\ref{EnHomology:QuasiFreeFiltration}.
Thus we obtain that $\phi: M\rightarrow N$
induces a morphism of spectral sequences $E^r(\phi): E^r(M)\rightarrow E^r(N)$.

\subsubsection{The $E^0$-term of the spectral sequence}\label{EnHomology:QuasiFreeEzero}
We adopt the notation $E^0_s = F_s/F_{s-1}$
for the subquotient of any filtration.
We have an obvious isomorphism $E^0_s(M) = (E^0_s(K)\circ\ROp,\partial_{\alpha})$,
where
\begin{equation*}
E^0_s(K)(r) = \begin{cases} K(s), & \text{if $r = s$}, \\ 0, & \text{otherwise}. \end{cases}
\end{equation*}
As a consequence,
the projection $F_s M\rightarrow E^0_s M$ has an obvious section $E^0_s M\rightarrow F_s M$
and we have a natural isomorphism $E^0(M)\simeq K\circ\ROp$.

\medskip
We apply the spectral sequence $E^0(M)\Rightarrow H_*(M)$
to the iterated bar module $B^n_{\EOp_n} = (T^n\circ\EOp_n,\partial_{\epsilon})$.
In one argument (lemma~\ref{EnHomology:HarrisonDifferential}),
we also use the spectral sequence of the $m$-fold bar complex
\begin{equation*}
B^m_{\EOp_n} = B^m_{\EOp_m}\circ_{\EOp_m}\EOp_n = (T^m\circ\EOp_n,\partial_{\epsilon}),
\end{equation*}
for $m<n$,
and the morphism of spectral sequences
$E^r(\sigma): E^r(\Sigma^{-m} B^m_{\EOp_n})\rightarrow E^r(\Sigma^{-n} B^n_{\EOp_n})$
induced by the composite suspension
\begin{equation*}
\Sigma^{-m} B^m_{\EOp_n}\xrightarrow{\sigma}\Sigma^{-m-1} B^{m+1}_{\EOp_n}\xrightarrow{\sigma}\dots
\xrightarrow{\sigma}\Sigma^{-n} B^n_{\EOp_n}.
\end{equation*}

We go back to the inductive definition of the twisting homomorphism
$\partial_{\epsilon}: T^n\circ\EOp_n\rightarrow T^n\circ\EOp_n$
in order to determine the differentials $d^r$ of $E^r(B^n_{\EOp_n})$
for $r = 0,1$.
In the case $r = 0$,
we obtain the following result:

\begin{prop}\label{EnHomology:EzeroTerm}
The chain complex $(E^0(B^n_{\EOp_n}),d^0)$ is isomorphic to the composite $B^n(I)\circ\EOp_n$,
where $B^n(I) = ((T^c\Sigma)^n(I),\partial)$
is the $n$-fold bar complex of the composition unit of $\Sigma_*$-modules,
viewed as a trivial commutative algebra in $\Sigma_*$-modules.
\end{prop}

\begin{proof}
The twisting homomorphism of the iterated bar complex $B^n(\COp)$
has a splitting $\partial_{\gamma} = \partial^0_{\gamma}+\partial^1_{\gamma}$
such that:
\begin{itemize}
\item
the component $\partial^1_{\gamma}: (T^c\Sigma)^n(I)\circ\COp\rightarrow(T^c\Sigma)^n(I)\circ\COp$
is yielded by the bar differential of the first factor of the composite
\begin{equation*}
B^n(\COp) = \BFact_n\circ\dots\circ\BFact_2\circ\BFact_1(\COp),
\end{equation*}
\item
the component $\partial^0_{\gamma}: (T^c\Sigma)^n(I)\circ\COp\rightarrow(T^c\Sigma)^n(I)\circ\COp$
is yielded by the bar differential of factors $2,\dots,n$ of $B^n(\COp)$.
\end{itemize}
The component $\partial^0_{\gamma}$
is identified with the homomorphism of free right $\COp$-modules
\begin{equation*}
\partial^0\circ\COp: (T^c\Sigma)^n(I)\circ\COp\rightarrow(T^c\Sigma)^n(I)\circ\COp,
\end{equation*}
where $\partial^0: (T^c\Sigma)^n(I)\rightarrow (T^c\Sigma)^n(I)$
is the differential of the iterated bar complex $B^n(I)$.

The twisting homomorphism $\partial_{\gamma}: T^c(\Sigma I)\circ\COp\rightarrow T^c(\Sigma I)\circ\COp$
of the bar complex $B(\COp) = (T^c(\Sigma I)\circ\COp,\partial_{\gamma})$
satisfies clearly $\partial_{\gamma}(F_s T^c(\Sigma I))\subset F_{s-1} T^c(\Sigma I)\circ\COp$.
This relation implies that the induced differential on $B^n(\COp)$
satisfies
\begin{equation*}
\partial_{\gamma}(F_s(T^c\Sigma)^n(I)\circ\COp)\subset F_{s-1}(T^c\Sigma)^n(I)\circ\COp
\end{equation*}
and vanishes in $E^0_s(B^n_{\COp})$.

Recall that the homomorphism $\epsilon_*: G^n\rightarrow T^n\circ\EOp_n$
which determines the twisting homomorphism of $B^n_{\EOp_n}$
is a sum $\epsilon = \sum_{m=0}^{\infty}\epsilon_m$,
whose terms $\epsilon_m$
are determined inductively by formulas of the form:
\begin{equation*}
\epsilon_0(\xi) = \tilde{\iota}\gamma(\xi)
\quad\text{and}
\quad\epsilon_m(\xi) = \sum_{p+q=m-1}\tilde{\nu}\partial_{\epsilon_p}\epsilon_q(\xi),
\end{equation*}
for any generating element $\xi\in G^n(r)$.
The homomorphisms $\tilde{\iota}$ and $\tilde{\nu}$ defined in~\S\ref{QuasiFreeLifting:TwistingCochainLifting}
satisfy clearly $\tilde{\iota}(F_s B^n_{\EOp_n})\subset F_s B^n_{\EOp_n}$
and $\tilde{\nu}(F_s B^n_{\EOp_n})\subset F_s B^n_{\EOp_n}$.
From this observation,
we deduce the relation
\begin{equation*}
\epsilon_0(\xi)\equiv\tilde{\iota}\partial^0(\xi)\mod F_{r-1} B^n_{\EOp_n}
\end{equation*}
and we conclude easily that the splitting $E^r(B^n_{\EOp_n})\simeq T^n\circ\EOp_n$
identifies the class of $d^0(\xi) = \epsilon_*(\xi) \mod F_{r-1} B^n_{\EOp_n}$
with $\partial^0(\xi)\in (T^c\Sigma)^n(I)$.

For a composite $\xi(p_1,\dots,p_r)\in T^n\circ\EOp_n$, where $\xi\in G^n(\eset)$,
we have:
\begin{equation*}
d^0(\xi(p_1,\dots,p_r))
= \partial^0(\xi)(p_1,\dots,p_r) + \sum_{i=1}^{r} \pm\xi(p_1,\dots,\delta(p_i),\dots,p_r),
\end{equation*}
where $\delta$ refers to the internal differential of~$\EOp_n$.
Hence,
we conclude that the differential $d^0: E^0(B^n_{\EOp_n})\rightarrow E^0(B^n_{\EOp_n})$
is identified with the natural differential
of the composite
\begin{equation*}
B^n(I)\circ\EOp_n = ((T^c\Sigma)^n(I),\partial^0)\circ\EOp_n.
\end{equation*}
The proposition follows.
\end{proof}

\subsubsection{The homology of $E_n$-operads and the Gerstenhaber operad}\label{EnHomology:GerstenhaberOperad}
We study the homology of the factors of the composite $(E^0,d^0) = B^n(I)\circ\EOp_n$
in order to determine the $E^1$-term of the spectral sequence $E^1(B^n_{\EOp_n})\Rightarrow H_*(B^n_{\EOp_n})$.
We still focus on the case $n>1$.

The homology of $B^n(I)$ is given by the result of proposition~\ref{UsualBarComplexes:IteratedShuffleAlgebra}.

The results of~\cite{Cohen} imply that the homology of~$\EOp_n$, $n>1$,
is isomorphic to the Gerstenhaber operad $\GOp_n$.
The structure of a $\GOp_n$-algebra consists of a commutative algebra $A$
equipped with a Lie bracket $\lambda_{n-1}: A\otimes A\rightarrow A$ of degree $n-1$
which satisfies a distribution relation with respect to the product $\mu: A\otimes A\rightarrow A$.
This description of the structure of a $\GOp_n$-algebra
reflects a definition of~$\GOp_n$ by generators and relations (see~\cite{GetzlerJones,Markl}).

The element $\mu\in\GOp_n$ which represents the commutative product of $\GOp_n$-algebras
generates a suboperad of $\GOp_n$
isomorphic to the commutative operad $\COp$.
The element $\lambda_{n-1}\in\GOp_n$ which represents the Lie product of $\GOp_n$-algebras
generates a suboperad of $\GOp_n$
isomorphic to the $(n-1)$-desuspension $\Lambda^{1-n}\LOp$ of the Lie operad $\LOp$.
The embeddings $\COp\hookrightarrow\GOp_n$ and $\Lambda^{1-n}\LOp\hookrightarrow\GOp_n$
assemble to an isomorphism $\COp\circ\Lambda^{1-n}\LOp\simeq\GOp_n$.

\begin{prop}\label{EnHomology:EoneTerm}
For $n>1$,
we have an isomorphism
\begin{equation*}
E^1(B^n_{\EOp_n})\simeq\COp(\Sigma^{n-1} L^c(\Sigma I))\circ\GOp_n
\simeq\COp(\Sigma^{n-1} L^c(\Sigma\GOp_n)).
\end{equation*}
\end{prop}

\begin{proof}
In~\S\ref{Background:Kunneth},
we recall that the K\"unneth morphism $H_*(M)\circ H_*(N)\rightarrow H_*(M\circ N)$
is an isomorphism for any composite $M\circ N$
such that the $\Sigma_*$-modules $M,N$ and their homology $H_*(M),H_*(N)$
consist of projective $\kk$-modules,
under the connectedness assumption $N(0) = 0$.

The $E_n$-operad $\EOp_n$ is $\Sigma_*$-cofibrant (and hence $\kk$-projective) by assumption.
The iterated bar complex $B^n(I)$ is $\Sigma_*$-cofibrant (and hence $\kk$-projective as well)
by the result of proposition~\ref{QuasiFreeLifting:SigmaFreeStructure}.

We have $H_*(B^n(I)) = \COp(\Sigma^{n-1} L^c(\Sigma I))$
by proposition~\ref{UsualBarComplexes:IteratedShuffleAlgebra}
and we just recall that $H_*(\EOp_n) = \GOp_n = \COp\circ\Lambda^{1-n}\LOp$.
The Lie operad $\LOp$ consists of free $\kk$-modules (see for instance~\cite[Chapitre II, \S 2]{Bourbaki} or \cite[Corollary 0.10]{Reutenauer})
and so does the commutative operad $\COp$.
The result of~\cite[Lemma 1.3.9]{FressePartitions} shows that a composite $M\circ N$,
where $M,N$ are $\kk$-projective $\Sigma_*$-modules,
is still $\kk$-projective
under the connectedness assumption $N(0) = 0$.
Thus, the $\Sigma_*$-modules $H_*(B^n(I))$ and $H_*(\EOp_n)$
are both $\kk$-projective.

From these observations,
we conclude that the K\"unneth morphism $H_*(M)\circ H_*(N)\rightarrow H_*(M\circ N)$
yields, in the case $M = B^n(I)$ and $N = \EOp_n$,
an isomorphism
\begin{equation*}
\COp(\Sigma^{n-1} L^c(\Sigma I))\circ\GOp_n
\simeq H_*(B^n(I))\circ H_*(\EOp_n)
\xrightarrow{\simeq} H_*(E^0(B^n_{\EOp_n}),d^0).
\end{equation*}
The proposition follows.
\end{proof}

The differential of certain particular elements of $E^1(B^n_{\EOp_n})$
can easily be determined:

\begin{lemm}\label{EnHomology:HarrisonDifferential}
The restriction of the differential $d^1: E^1\rightarrow E^1$
to the summand
\begin{equation*}
\Sigma^{n-1} L^c(\Sigma\GOp_n)\subset\COp(\Sigma^{n-1} L^c(\Sigma\GOp_n))
\end{equation*}
is identified with the differential of the Harrison complex $(L^c(\Sigma\GOp_n),\partial)$.
\end{lemm}

\begin{proof}
In~\S\ref{EnHomology:QuasiFreeFiltration},
we observe that the suspension
\begin{equation*}
B(\EOp_n)\xrightarrow{\sigma} B^n(\EOp_n) = B^n_{\EOp_n}
\end{equation*}
determines a morphism of spectral sequences
\begin{equation*}
E^r(\sigma): E^r(B(\EOp_n))\rightarrow E^r(B^n(\EOp_n)).
\end{equation*}
For the spectral sequence $E^r(B(\EOp_n))\Rightarrow H_*(B(\EOp_n))$,
we also have $(E^0,d^0) = T^c(\Sigma\EOp_n)$
since the differential of the $1$-fold bar complex decreases filtrations
(see proof of proposition~\ref{EnHomology:EzeroTerm}).
Hence, we obtain:
\begin{equation*}
E^1(B(\EOp_n)) = T^c(\Sigma H_*(\EOp_n)) = T^c(\Sigma\GOp_n).
\end{equation*}

By an immediate inspection of constructions,
we see that the morphism $\Sigma^{n-1} L^c(\Sigma\GOp_n)\rightarrow E^1(B^n_{\EOp_n})$
fits a commutative diagram
\begin{equation*}
\xymatrix{ & \Sigma^{n-1} T^c(\Sigma\GOp_n)\ar[]!LD;[ddl]_{=}\ar[]!RD;[ddr]\ar[d] & \\
& \Sigma^{n-1} L^c(\Sigma\GOp_n)\ar@{.>}[]!LD;[dl]\ar@{.>}[]!RD;[dr] & \\
\Sigma^{n-1} E^1(B(\EOp_n))\ar[rr]_{E^1(\sigma)} && E^1(B^n(\EOp_n)) }.
\end{equation*}
Thus we are reduced to determine the bar differential
of the representative $\xi = p_1\otimes\dots\otimes p_s\in T^c(\Sigma\EOp_n)$
of an element of $T^c(\Sigma\GOp_n)$.

The definition of the differential of $B(\EOp_n)$
gives immediately:
\begin{align*}
\partial(p_1\otimes\dots\otimes p_s)
& = \sum_{t=2}^{s}\Bigl\{\sum_{i=1}^{s-t+1} \pm p_1\otimes\dots\otimes\mu_t(p_{i},\dots,p_{i+t-1})\otimes
\dots\otimes p_s\Bigr\} \\
& \equiv\sum_{i=1}^{s-1} \pm p_1\otimes\dots\otimes\mu_2(p_{i},p_{i+1})\otimes\dots\otimes p_s
\mod F_{s-2} B(\EOp_n),
\end{align*}
where $\mu_2\in\EOp_2(2)$ is a representative of the product $\mu\in\GOp_2(2)$.
Hence we obtain that $d^1(p_1\otimes\dots\otimes p_s)$ is given by the differential
of the bar complex $B(\GOp_n)$.
The conclusion follows.
\end{proof}

\begin{lemm}\label{EnHomology:ChevalleyDifferential}
For an element of the form
\begin{equation*}
\xi = x(e_1)\cdot x(e_2)\in\COp(\Sigma^{n-1} L^c(\Sigma I))(\{e_1,e_2\}),
\end{equation*}
where $x$ represents the canonical generator of~$\Sigma I(1) = L^c(\Sigma I)(1)$,
we have
\begin{equation*}
d^1(\xi) = \lambda_{n-1}(e_1,e_2)\in\GOp_n(\{e_1,e_2\}),
\end{equation*}
where $\lambda_{n-1}$ is the operation of $\GOp_n(2)$
which represents the Lie bracket of $\GOp_n$-algebras.
\end{lemm}


\begin{proof}
In this proof,
it is convenient to adopt the notation $\otimes_m$
to refer to a tensor product in the $m$th factor of the composite
\begin{equation*}
(T^c\Sigma)^n(I) = \TFact_n\circ\dots\circ\TFact_1(I).
\end{equation*}
Observe that
\begin{equation*}
(T^c\Sigma)^n(I)(\{e_1,e_2\})
= \bigoplus_{m=1}^{n}\Bigl\{\kk x(e_1)\otimes_m x(e_2)\oplus\kk x(e_2)\otimes_m x(e_1)\Bigl\},
\end{equation*}
where $x$ refers to a generator of $I(1) = \kk$.
The element $\xi$ is represented by the product of $x(e_1)$ and $x(e_2)$
in $B^n(\EOp_n)$:
\begin{equation*}
x(e_1)\cdot x(e_2) = x(e_1)\otimes_n x(e_2) + \pm x(e_2)\otimes_n x(e_1)\in T^c(\Sigma(T^c\Sigma)^{n-1}(I)).
\end{equation*}

\textbf{1)}
First,
we prove inductively that the differential in the iterated bar complex $B^n(\EOp_n)$
of an element of the form $x(e_1)\otimes_m x(e_2)$
is given by a sum:
\begin{equation*}
\partial_{\epsilon}(x(e_1)\otimes_m x(e_2)) = \partial^0(x(e_1),x(e_2)) + \upsilon_m(e_1,e_2)
\end{equation*}
where
\begin{multline*}
\partial^0(x(e_1),x(e_2)) = x(e_1)\otimes_{m-1} x(e_2) + \pm x(e_1)\otimes_{m-1} x(e_2)\\
\in(T^c\Sigma)^{m-1}(I)(\{e_1,e_2\})
\end{multline*}
is the shuffle product of $x(e_1)$ and $x(e_2)$ in $(T^c\Sigma)^{m-1}(I)$
and
\begin{equation*}
\upsilon_m(e_1,e_2)\in\EOp_n(\{e_1,e_2\})\subset(T^c\Sigma)^n(\EOp_n)(\{e_1,e_2\})
\end{equation*}
is a representative of the $\cup_m$-product.

Recall that the operad $\EOp$
is supposed to be equipped with a chain-homotopy $\nu: \EOp\rightarrow\EOp$
so that $\delta(\nu) = \id - \iota\epsilon$
for some fixed section $\iota: \COp\rightarrow\EOp$
of the augmentation $\epsilon: \EOp\rightarrow\COp$.
The $\cup_m$-products $\upsilon_m$
are defined inductively by $\upsilon_0 = \iota\mu$,
where $\mu\in\COp(2)$ represents the product of commutative algebras,
and $\upsilon_m = \nu(\upsilon_{m-1}+\pm\tau\upsilon_{m-1})$,
where $\tau\in\Sigma_2$ denotes the transposition of $(1,2)$.
Our assumptions on $E_\infty$-operads
ensure that $\upsilon_m\in\EOp_n(2)$ for $m<n$.
The cocycle $\upsilon_{n-1}+\pm\tau\upsilon_{n-1}\in\EOp_n(2)$
defines a representative of the operation $\lambda_{n-1}\in\GOp_n(2)$.

By equivariance, we can assume $e_1<e_2$.
Observe that $x(e_1)\otimes_m x(e_2)$
belongs to the generating $\Sigma_*$-module of $B^n(\EOp_n)$.
By definition of the twisting homomorphism $\partial_{\epsilon}$,
the differential of $x(e_1)\otimes_m x(e_2)$ in $B^n(\EOp_n)$
has an expansion of the form
\begin{equation*}
\partial_{\epsilon}(x(e_1)\otimes_m x(e_2)) = \sum_{r=0}^{\infty}\epsilon_r(x(e_1)\otimes_m x(e_2)).
\end{equation*}
For $m=1$,
the definition of~\S\ref{QuasiFreeLifting:TwistingCochainLifting}
returns
\begin{equation*}
\epsilon_0(x(e_1)\otimes_1 x(e_2)) = \iota\mu(e_1,e_2)\in\EOp_n(\{e_1<e_2\}).
\end{equation*}
For $m>1$,
we obtain
\begin{equation*}
\epsilon_0(x(e_1)\otimes_m x(e_2)) = x(e_1)\otimes_{m-1} x(e_2) + \pm x(e_2)\otimes_{m-1} x(e_1),
\end{equation*}
the shuffle product of $x(e_1)$ and $x(e_2)$ in $(T^c\Sigma)^{m-1}(I)$.

From the definition of~\S\ref{QuasiFreeLifting:TwistingCochainLifting},
we see by an easy induction on $m$ and $r$ that the term $\epsilon_r(x(e_1)\otimes_m x(e_2))$
vanishes for $0<r<m-1$.

For $r = m-1$,
we obtain:
\begin{align*}
\epsilon_{m-1}(x(e_1)\otimes_m x(e_2))
& = \tilde{\nu}\cdot\partial_{\epsilon_{m-2}}\cdot\partial_{\epsilon_{0}}(x(e_1)\otimes_m x(e_2))\\
& = \tilde{\nu}\partial_{\epsilon_{m-2}}(x(e_1)\otimes_{m-1} x(e_2) + \pm x(e_2)\otimes_{m-1} x(e_1)),
\end{align*}
where $\tilde{\nu}$ arises from a natural extension of the chain-homotopy $\nu$.
This identity gives by an immediate induction
\begin{equation*}
\epsilon_{m-1}(x(e_1)\otimes_m x(e_2)) = \upsilon_m(e_1,e_2)\in\EOp_n(\{e_1<e_2\}),
\end{equation*}
and the terms $\epsilon_r(x(e_1)\otimes_m x(e_2))$
are trivial for $r>m-1$.

\textbf{2)}
From the result
\begin{equation*}
\partial_{\epsilon}(x(e_1)\otimes_n x(e_2))
= x(e_1)\otimes_{n-1} x(e_2) + x(e_2)\otimes_{n-1} x(e_1) + \upsilon_n(e_1,e_2),
\end{equation*}
we deduce:
\begin{multline*}
\begin{aligned} \partial_{\epsilon}(x(e_1)\cdot x(e_2))
= & \partial_{\epsilon}(x(e_1)\otimes_n x(e_2)) + \pm\partial_{\epsilon}(x(e_2)\otimes_n x(e_1))\\
= & (x(e_1)\otimes_{n-1} x(e_2) + \pm\pm x(e_1)\otimes_{n-1} x(e_2))\\
+ & (\pm x(e_2)\otimes_{n-1} x(e_1) + \pm x(e_2)\otimes_{n-1} x(e_1)) \end{aligned} \\
+ (\upsilon_n(e_1,e_2) + \pm\upsilon_n(e_2,e_1)).
\end{multline*}
In our verification,
we have not specified any sign,
but the coherence of dg-algebra ensures that:
\begin{itemize}
\item
the terms $x(e_1)\otimes_{n-1} x(e_2)$
cancel each other in the expansion of $\partial_{\epsilon}(x(e_1)\cdot x(e_2))$,
and so do the terms $x(e_2)\otimes_{n-1} x(e_1)$;
\item
the terms $\upsilon_n(e_1,e_2)$ and $\upsilon_n(e_2,e_1)$
assemble to the Lie bracket $\lambda_{n-1}(e_1,e_2) = \upsilon_n(e_1,e_2) + \pm\upsilon_n(e_2,e_1)$.
\end{itemize}
Thus we obtain
\begin{equation*}
\partial_{\epsilon}(x(e_1)\cdot x(e_2)) = \lambda_{n-1}(e_1,e_2).
\end{equation*}
The conclusion about the differential of~$x(e_1)\cdot x(e_2)$ in~$E^1(B^n_{\EOp_n})$ follows.
\end{proof}

We use the coalgebra structure of the iterated bar complexes
to determine the $E^2$-term of the spectral sequence $E^1(B^n_{\EOp_n})\Rightarrow H_*(B^n_{\EOp_n})$
from the partial calculations of lemmas~\ref{EnHomology:HarrisonDifferential}-\ref{EnHomology:ChevalleyDifferential}.

\begin{claim}\label{EnHomology:CoalgebraStructures}
The filtration of~\S\ref{EnHomology:QuasiFreeFiltration}
satisfies the relation
$\Delta(F_r B^n_{\EOp_n})\subset\sum_{s+t = r} F_s B^n_{\EOp_n}\otimes F_t B^n_{\EOp_n}$
with respect to the coproduct of~$B^n_{\EOp_n}$.
\end{claim}

\begin{proof}
Indeed,
the isomorphism $(T^c\Sigma)^n(M)\simeq(T^c\Sigma)^n(I)\circ M = T^n\circ M$
defines an isomorphism of coalgebras
in the sense that the diagonal of~$(T^c\Sigma)^n(M)$
corresponds to the composite of the morphism
\begin{equation*}
\Delta\circ M: T^n\circ M\rightarrow(T^n\otimes T^n)\circ M
\end{equation*}
induced by the diagonal of $T^n = (T^c\Sigma)^n(I)$
with the distribution isomorphism
\begin{equation*}
(T^n\otimes T^n)\circ M\simeq T^n\circ M\otimes T^n\circ M.
\end{equation*}
The diagonal of $T^n$
satisfies
$\Delta(T^n(r))\subset\bigoplus_{s+t=r}\Sigma_r\otimes_{\Sigma_s\times\Sigma_t} T^n(s)\otimes T^n(t)$
by definition of the tensor product of $\Sigma_*$-modules.
Hence,
we obtain $\Delta(F_r T^n)\subset\sum_{s+t=r} F_s T^n\otimes F_t T^n$,
from which we deduce $\Delta(F_r T^n\circ\EOp_n)\subset\sum_{s+t=r} F_s T^n\circ\EOp_n\otimes F_t T^n\circ\EOp_n$
and the claim follows.
\end{proof}

This observation
implies that the spectral sequence of~\S\ref{EnHomology:QuasiFreeFiltration}
defines a spectral sequence of coalgebras $E^r(B^n_{\EOp_n})\Rightarrow H_*(B^n_{\EOp_n})$.

Recall that any commutative algebra $\COp(M)$
is equipped with a natural Hopf algebra structure
so that $M$ is primitive in $\COp(M)$.
The result of proposition~\ref{EnHomology:EoneTerm}
can be improved to:

\begin{lemm}\label{EnHomology:HopfEoneTerm}
The isomorphism of proposition~\ref{EnHomology:EoneTerm}
\begin{equation*}
E^1(B^n_{\EOp_n})\simeq\COp(\Sigma^{n-1} L^c(\Sigma\GOp_n))
\end{equation*}
defines an isomorphism of coalgebras.
\end{lemm}

\begin{proof}
Let $M = \Sigma^{n-1} L^c(\Sigma I)$.
The distribution isomorphism
\begin{equation*}
\COp(M\circ\GOp_n)\simeq\COp(M)\circ\GOp_n
\end{equation*}
maps the natural diagonal of the commutative algebra
to the composite
\begin{equation*}
\COp(M)\circ\GOp_n
\xrightarrow{\Delta\circ\GOp_n}(\COp(M)\otimes\COp(M))\circ\GOp_n
\simeq\COp(M)\circ\GOp_n\otimes\COp(M)\circ\GOp_n,
\end{equation*}
where $\Delta$ refers to the diagonal of $\COp(M)$.

The isomorphism $H_*(B^n(I))\simeq\COp(\Sigma^{n-1} L^c(\Sigma I))$
is an isomorphism of coalgebras by proposition~\ref{UsualBarComplexes:HopfAlgebraIso}
and
the K\"unneth morphism $H_*(T^n)\circ H_*(\EOp_n)\rightarrow H_*(T^n\circ\EOp_n)$
forms clearly a morphism of coalgebras
with respect to the coalgebra structure yielded by $T^n$.
The conclusion follows immediately.
\end{proof}

\subsubsection{Coderivations on cofree coalgebras}\label{EnHomology:CofreeCoderivations}
Recall (see~\S\ref{UsualBarComplexes:HopfAlgebraStructures})
that the diagonal of a monomial $\xi = x_1\cdot\,\dots\,\cdot x_r\in\COp(M)$
in the free commutative algebra $\COp(M)$
is defined by the formula:
\begin{multline*}
\Delta(x_1\cdot\,\dots\,\cdot x_r) \\
= \sum_{\substack{p+q=r\\p,q>0}}\Bigl\{\sum_{w\in\Sh(p,q)}\pm(x_{w(1)}\cdot\,\dots\,\cdot x_{w(p)})
\otimes(x_{w(p+1)}\cdot\,\dots\,\cdot x_{w(p+q)})\Bigr\},
\end{multline*}
where the inner sum ranges over the set of $(p,q)$-shuffles in $\Sigma_r$.
Recall again that $\COp(M)$
refers to a free commutative algebra without unit.
Therefore
the expansion of~$\Delta(x_1\cdot\,\dots\,\cdot x_r)$
runs over pairs $(p,q)$ such that $p,q>0$.
Note that the $(p,q)$-shuffle
$(x_{w(1)}\cdot\,\dots\,\cdot x_{w(p)})\otimes(x_{w(p+1)}\cdot\,\dots\,\cdot x_{w(p+q)})$
includes a permutation of the inputs $\eset_i$ of the elements $x_i\in M(\eset_i)$.

Let $\COp_r(M) = (M^{\otimes r})_{\Sigma_r}$.
The $n$-fold diagonal $\Delta^n: \COp(M)\rightarrow\COp(M)^{\otimes n}$
is defined inductively by $\Delta^n = \Delta^{n-1}\otimes\id\cdot\Delta$
(as in the unital setting since the coassociativity relation does not change).
The composite of the $n$-fold diagonal $\Delta^n$
with the canonical projection $\COp(M)\rightarrow\COp_1(M) = M$
vanishes over the summands $\COp_s(M)$, $s\not=r$,
and is identified with the trace morphism on the summand $\COp_r(M)$.
Recall that the trace morphism $\Tr_{\Sigma_r}$
is the morphism
\begin{equation*}
(M^{\otimes r})_{\Sigma_r}\xrightarrow{\Tr_{\Sigma_r}} (M^{\otimes r})^{\Sigma_r}\subset M^{\otimes r}
\end{equation*}
defined by the sum of all tensor permutations $w_*: M^{\otimes r}\rightarrow M^{\otimes r}$, $w\in\Sigma_r$.

For a $\Sigma_*$-module $M$ such that $M(0) = 0$
the trace $\Tr_{\Sigma_*}$ is an isomorphism
(see~\cite[\S 1.1]{FresseDividedPowers} or~\cite[\S\S 1.2.12-1.2.19]{FressePartitions})
and this observation implies that $\COp(M) = \bigoplus_{n=1} (M^{\otimes r})_{\Sigma_r}$
satisfies the universal property of a cofree coalgebra (without counit):
any morphism $f: \Gamma\rightarrow M$,
where $\Gamma$ is a coalgebra in $\Sigma_*$-modules,
admits a unique lifting
\begin{equation*}
\xymatrix{ & \COp(M)\ar[d] \\ \Gamma\ar[r]_{f}\ar@{.>}[ur]^{\psi_f} & M }
\end{equation*}
such that $\psi_f$ is a morphism of coalgebras.

Similarly,
any homomorphism $\varpi: \COp(M)\rightarrow M$
has a unique lifting
\begin{equation*}
\xymatrix{ & \COp(M)\ar[d] \\ \COp(M)\ar[r]_{\varpi}\ar@{.>}[ur]^{\theta_{\varpi}} & M }
\end{equation*}
such that $\theta_{\varpi}$ is a coderivation.
This coderivation is given on products $x_1\cdot\,\dots\,\cdot x_r\in\COp_r(M)$
by the formula:
\begin{multline*}
\theta_{\varpi}(x_1\cdot\,\dots\,\cdot x_r) \\
= \sum_{\substack{p+q=r\\p,q>0}}\Bigl\{\sum_{w\in\Sh(p,q)} \pm x_{w(1)}\cdot\,\dots\,\cdot x_{w(p)}
\cdot\varpi(x_{w(p+1)}\cdot\,\dots\,\cdot x_{w(p+q)})\Bigr\}.
\end{multline*}
Note that $\varpi$ is recovered from the associated coderivation $\theta_{\varpi}$
by the composite of $\theta_{\varpi}: \COp(M)\rightarrow\COp(M)$
with the canonical projection $\epsilon: \COp(M)\rightarrow M$.

Lemma~\ref{EnHomology:HopfEoneTerm}
implies that the differential $d^1$ of the spectral sequence $E^1(B^n_{\EOp_n})\Rightarrow H_*(B^n_{\EOp_n})$
is given by a formula of this form.

\subsubsection{The natural spectral sequence of a cofree coalgebra}\label{EnHomology:CofreeFiltration}
The cofree coalgebra $\COp(M)$
admits a canonical filtration
defined by $F_s\COp(M) = \bigoplus_{r\leq s}\COp_r(M)$.

Let $\Gamma = (\COp(M),\partial)$
be a quasi-cofree cofree coalgebra
so that the twisting homomorphism $\partial: \COp(M)\rightarrow\COp(M)$
is a coderivation.
Thus we have $\partial = \theta_{\varpi}$,
for a certain homomorphism $\varpi: \COp(M)\rightarrow M$.

The formula of~\S\ref{EnHomology:CofreeCoderivations}
implies that $\theta_{\varpi}(F_s\COp(M))\subset F_s\COp(M)$.
From this observation,
we deduce that the quasi-cofree coalgebra $\Gamma = (\COp(M),\partial)$
has a filtration by twisted dg-modules
such that:
\begin{equation*}
F_s\Gamma = (F_s\COp(M),\partial).
\end{equation*}
Let $D^r(\Gamma)\Rightarrow H_*(\Gamma)$
be the spectral sequence defined by this filtration.

Note that $\Delta(F_r\COp(M))\subset\sum_{p+q = r} F_p\COp(M)\otimes F_q\COp(M)$.
This relation implies that the spectral sequence $D^r(\Gamma)\Rightarrow H_*(\Gamma)$
forms a spectral sequence of coalgebras.

We have an obvious isomorphism $D^0(\Gamma)\simeq\COp(M)$.

Recall that $\varpi$ is determined from $\partial = \theta_{\varpi}$
by the composite
\begin{equation*}
\COp(M)\xrightarrow{\partial}\COp(M)\xrightarrow{\epsilon} M,
\end{equation*}
where $\epsilon$ refers to the canonical projection of $\COp(M)$ onto the summand $M\subset\COp(M)$.
Let $\varpi_r: \COp_{r+1}(M)\rightarrow M$
be the restriction of $\varpi: \COp(M)\rightarrow M$
to the summand $\COp_{r+1}(M)\subset\COp(M)$, for $r\in\NN$.
Let $\partial^r = \theta_{\varpi_r}$
be the coderivation associated to this homomorphism.
We have clearly $\varpi = \sum_r\varpi_r\Rightarrow\partial = \sum_{r=0}^{\infty} \partial^r$.

The formula of~\S\ref{EnHomology:CofreeCoderivations}
implies moreover $\partial^r(F_s\COp(M)) = \theta_{\varpi_r}(F_s\COp(M))\subset F_{s-r}\COp(M)$.
As a consequence,
the differential $d^0$ of the spectral sequence $D^r(\Gamma)$
is identified with the coderivation $\partial^0 = \theta_{\varpi_0}: \COp(M)\rightarrow\COp(M)$.
Note that the assumption $\varpi_0(M)\subset M$
implies that $\partial^0 = \theta_{\varpi_0}$
defines a differential on $M$
and the quasi-cofree coalgebra $(D^0,d^0) = (\COp(M),\partial^0)$
is isomorphic to the cofree coalgebra $(D^0,d^0) = \COp(M,\partial^0)$
cogenerated by the twisted $\Sigma_*$-object $(M,\partial^0)$.

\medskip
We apply the spectral sequence of cofree coalgebras
to the quasi-cofree coalgebra:
\begin{equation*}
E^1(B^n_{\EOp_n}) = \COp(\Sigma^{n-1} L^c(\Sigma\GOp_n)).
\end{equation*}
We have in this case:

\begin{lemm}\label{EnHomology:DzeroDifferential}
The $D^0$-term of the spectral sequence $D^r(E^1(B^n_{\EOp_n}))$
is identified with the cofree coalgebra on the $(n-1)$-fold suspension
of the Harrison complex $(L^c(\Sigma\GOp_n),\partial)$.
\end{lemm}

\begin{proof}
By observations of~\S\ref{EnHomology:CofreeFiltration},
the differential $d^0: D^0\rightarrow D^0$
is the coderivation of $E^1$
associated to the component
\begin{equation*}
\xymatrix{ \COp(\Sigma^{n-1} L^c(\Sigma\GOp_n))\ar[r]^{d^1} &
\COp(\Sigma^{n-1} L^c(\Sigma\GOp_n))\ar[d] \\
\Sigma^{n-1} L^c(\Sigma\GOp_n)\ar@{.>}[r]_{\partial^0}\ar[u] & \Sigma^{n-1} L^c(\Sigma\GOp_n) }
\end{equation*}
of the differential of $E^1$.
In proposition~\ref{EnHomology:HarrisonDifferential},
we observe that this component of $d^1$
is identified with the Harrison differential.
The conclusion follows.
\end{proof}

From this lemma,
we deduce:

\begin{prop}\label{EnHomology:DoneTerm}
We have $D^1(E^1(B^n_{\EOp_n})) = \COp(\Sigma^n\Lambda^{1-n}\LOp)$.
\end{prop}

\begin{proof}
The composite $\GOp_n = \COp\circ\Lambda^{1-n}\LOp$
is identified with the free commutative algebra $\GOp_n = \COp(\Lambda^{1-n}\LOp)$.
Thus we have $H_*(L^c(\Sigma\GOp_n),\partial) = \Lambda^{1-n}\LOp$
by proposition~\ref{UsualBarComplexes:HarrisonKoszul}.

Since the $\Sigma_*$-modules $L^c(\Sigma\GOp_n)$ and $\Lambda^{1-n}\LOp$
are both $\kk$-projective and connected,
the K\"unneth morphism
\begin{equation*}
\COp\circ H_*(L^c(\Sigma\GOp_n),\partial)\rightarrow H_*(\COp\circ(L^c(\Sigma\GOp_n),\partial))
\end{equation*}
is an isomorphism (according to~\S\ref{Background:Kunneth}).
Equivalently,
we have an isomorphism
\begin{equation*}
\COp(H_*(L^c(\Sigma\GOp_n),\partial))
\xrightarrow{\simeq} H_*(\COp(L^c(\Sigma\GOp_n),\partial))
= H_*(D^0(\Gamma),d^0),
\end{equation*}
from which we conclude $D^1(\Gamma) = \COp(\Sigma^n\Lambda^{1-n}\LOp)$.
\end{proof}

For the next stage of the spectral sequence,
we obtain:

\begin{lemm}\label{EnHomology:DoneDifferential}
For $\Gamma = E^1(B^n_{\EOp_n})$,
the chain complex $(D^1(\Gamma),d^1(\Gamma))$
is identified with the Chevalley-Eilenberg complex of the free Lie algebra $\LOp(\Sigma^{n-1} I)$.
\end{lemm}

\begin{proof}
Note that $\Sigma^{n-1}\Lambda^{1-n}\LOp = \LOp(\Sigma^{n-1} I)$
by definition of operadic suspensions.

By observations of~\S\ref{EnHomology:CofreeFiltration},
the component
\begin{equation*}
\xymatrix{ \COp(\Sigma^{n-1} L^c(\Sigma\GOp_n))\ar[r]^{\partial} &
\COp(\Sigma^{n-1} L^c(\Sigma\GOp_n))\ar[d] \\
\COp_{r+1}(\Sigma^{n-1} L^c(\Sigma\GOp_n))\ar@{.>}[r]_{\varpi_{r+1}}\ar[u] & \Sigma^{n-1} L^c(\Sigma\GOp_n) }
\end{equation*}
of the differential of $E^1$
determines a coderivation $\partial^{r}$
such that $\partial^{r}(F_s E^1)\subset F_{s-r} E^1$.
Thus,
to determine the differential $d^1: D^1\rightarrow D^1$,
we have to determine the restriction of~$d^1: E^1\rightarrow E^1$
to the submodule
\begin{equation*}
D^1_2 = \COp_2(\Sigma^{n}\Lambda^{1-n}\LOp)\subset\COp_2(\Sigma^{n-1} L^c(\Sigma\GOp_n)).
\end{equation*}

In lemma~\ref{EnHomology:ChevalleyDifferential},
we prove that the differential in $E^1$
of an element of the form
\begin{equation*}
x(e_1)\cdot x(e_2)\in\COp_2(\Sigma^{n} I)(\{e_1,e_2\})
\end{equation*}
is given by $d^1(x(e_1)\cdot x(e_2)) = \lambda_{n-1}(e_1,e_2)$.
Since $d^1: E^1\rightarrow E^1$ is a differential of right $\GOp_n$-modules,
we have
$d^1(p_1\cdot p_2) = \lambda_{n-1}(p_1,p_2)$
for any product of elements $p_1,p_2\in\LOp(\Sigma^{n-1} I)$.

Hence,
we have $\varpi_2(p_1\cdot p_2) = \lambda(p_1,p_2)$
and the coderivation associated to $\varpi_2$
has an expansion of the form:
\begin{equation*}
\partial^1(p_1\cdot\,\dots\,\cdot p_r)
= \sum\pm p_1\cdot\,\dots\,\widehat{p_i}\,\dots\,\widehat{p_j}\,\dots\,\cdot\lambda_{n-1}(p_i,p_j).
\end{equation*}
The terms on the right-hand side belong to:
\begin{equation*}
D^1 = \COp_r(\Sigma^{n}\Lambda^{1-n}\LOp)\subset\COp_2(\Sigma^{n-1} L^c(\Sigma\GOp_n)).
\end{equation*}
Thus the expansion of $d^1: D^1\rightarrow D^1$
is given by the same formula,
which is also identified with the expression of the Chevalley-Eilenberg differential.
The lemma follows.
\end{proof}

Then:

\begin{prop}\label{EnHomology:DtwoTerm}
We have $D^2(E^1(B^n_{\EOp_n})) = \Sigma^n I$.
\end{prop}

\begin{proof}
We have the distribution relation
\begin{equation*}
(\COp(\Sigma\LOp(\Sigma^{n-1} I)),\partial) = (\COp(\Sigma\LOp(I)),\partial)\circ(\Sigma^{n-1} I),
\end{equation*}
and the Chevalley-Eilenberg complex $(\COp(\Sigma\LOp(\Sigma^{n-1} I)),\partial)$
is identified with a composite
\begin{equation*}
(\COp(\Sigma\LOp(\Sigma^{n-1} I)),\partial) = \Sigma K(I,\LOp,\LOp)\circ(\Sigma^{n-1} I)
\end{equation*}
where $K(I,\LOp,\LOp)$ is the Koszul complex of the Lie operad
(see~\cite[\S\S 6.3-6.4]{FressePartitions}).
This complex $K(I,\LOp,\LOp)$
is acyclic because the Lie operad is Koszul (see \emph{loc. cit.}).
Therefore,
we obtain
\begin{equation*}
D^2 = H_*(\COp(\Sigma\LOp(\Sigma^{n-1} I)),\partial) = \Sigma^n I.
\end{equation*}
\end{proof}

This result implies immediately:

\begin{prop}\label{EnHomology:EtwoTerm}
The spectral sequence $D^r(E^1,d^1)\Rightarrow H_*(E^1,d^1)$ degenerates at $D^2$
and gives $H_*(E^1,d^1) = \Sigma^n I$.\qed
\end{prop}

Thus, we have
\begin{equation*}
E^2(B^n_{\EOp_n}) = H_*(E^1,d^1) = \Sigma^n I,
\end{equation*}
from which we obtain:

\begin{prop}\label{EnHomology:SpectralSequenceResult}
The spectral sequence $E^r(B^n_{\EOp_n})\Rightarrow H_*(B^n_{\EOp_n})$
degenerates at the $E^2$-stage and returns $H_*(B^n_{\EOp_n}) = \Sigma^n I$.\qed
\end{prop}

Finally,
we conclude:

\begin{thm}\label{EnHomology:IteratedBarModuleHomotopy}
The augmentation $\epsilon: B^n_{\EOp_n}\rightarrow\Sigma^n I$
defines a weak-equivalence of right $\EOp_n$-modules
and the bar module $B^n_{\EOp_n}$ produced by theorem~\ref{EnDefinition:TwistingCochainRestriction}
defines, after desuspension, a cofibrant replacement of the composition unit~$I$
in the category of right $\EOp_n$-modules.\qed
\end{thm}

As we explain in~\S\ref{EnHomology:TorDefinition},
this theorem implies:

\begin{thm}\label{EnHomology:Conclusion}
Let $\EOp$ be an $E_\infty$-operad satisfying all requirements of~\S\ref{EnDefinition:KCellOperads}
so that the conclusions of theorem~\ref{EnDefinition:TwistingCochainRestriction}
and theorem~\ref{EnDefinition:Result} hold.
We have $H_*(\Sigma^{-n} B^n(A)) = \Tor^{\EOp_n}_*(I,A) = H^{\EOp_n}_*(A)$
for every $\EOp_n$-algebra $A\in{}_{\EOp_n}\E$
which forms a cofibrant object in the underlying category $\E$.\qed
\end{thm}

\section{Infinite bar complexes}\label{EinfinityHomology}
The goal of this section is to extend the result of theorem~\ref{EnHomology:Conclusion}
to the case $n = \infty$.
First of all,
we have to define an infinite bar complex $\Sigma^{-\infty} B^{\infty}(A)$
on the category of algebras over an $E_\infty$-operad $\EOp$.

\subsubsection{The infinite bar complex and the infinite bar module}\label{EinfinityHomology:infiniteBarModule}
Throughout this section,
we use the letter $\ROp$ to refer either to the commutative operad $\ROp = \COp$ or to an $E_\infty$-operad $\ROp = \EOp$.
In our construction,
we observe that the iterated bar modules $B^n_{\ROp}$
are connected by suspension morphisms
so that the diagram
\begin{equation*}
\xymatrix{ \Sigma^{-1} B_{\EOp}\ar[d]_{\sim}\ar[r]^(0.6){\sigma_{\EOp}} & \ar@{.}[r] &
\ar[r]^(0.3){\sigma_{\EOp}} & \Sigma^{1-n} B^{n-1}_{\EOp}\ar[d]_{\sim}\ar[r]^(0.55){\sigma_{\EOp}} &
\Sigma^{-n} B^{n}_{\EOp}\ar[d]_{\sim}\ar[r]^(0.6){\sigma_{\EOp}} & \ar@{.}[r] & \\
\Sigma^{-1} B_{\COp}\ar[r]^(0.6){\sigma_{\COp}} & \ar@{.}[r] &
\ar[r]^(0.3){\sigma_{\COp}} & \Sigma^{1-n} B^{n-1}_{\COp}\ar[r]^(0.55){\sigma_{\COp}} &
\Sigma^{-n} B^{n}_{\COp}\ar[r]^(0.6){\sigma_{\COp}} & \ar@{.}[r] & }
\end{equation*}
commutes.
For short, we may write $\sigma = \sigma_{\ROp}$.
These suspension morphisms induce natural transformations at the functor level
\begin{equation*}
\Sigma^{-1} B(A)\xrightarrow{\sigma}\cdots
\xrightarrow{\sigma}\Sigma^{1-n} B^{n-1}(A)
\xrightarrow{\sigma}\Sigma^{-n} B^{n}(A)
\xrightarrow{\sigma}\cdots.
\end{equation*}

Set
\begin{equation*}
\Sigma^{-\infty}B^{\infty}_{\ROp} = \colim_n\Sigma^{-n} B^{n}_{\ROp}
\end{equation*}
and form the infinite bar complex
\begin{equation*}
\Sigma^{-\infty}B^{\infty}(A) = \colim_n\Sigma^{-n} B^{n}(A).
\end{equation*}
Since the functor $(M,A)\mapsto\Sym_{\ROp}(M,A)$
commutes with colimits in $M$,
we have the relation
\begin{equation*}
\Sigma^{-\infty}B^{\infty}(A) = \Sym_{\ROp}(\Sigma^{-\infty}B^{\infty}_{\ROp},A)
\end{equation*}
for every $\ROp$-algebra $A$.

The relationship $\Sigma^{-n} B^n_{\EOp}\circ_{\EOp}\COp\simeq\Sigma^{-n} B^n_{\COp}$
extends to $n = \infty$,
because the suspension morphisms satisfy the coherence relation $\sigma_{\EOp}\circ_{\EOp}\COp = \sigma_{\COp}$,
and the diagram
\begin{equation*}
\xymatrix{ {}_{\EOp}\E\ar@{.>}[dr]_{\Sigma^{-n}B^{n} = \Sym_{\EOp}(\Sigma^{-n}B^{n}_{\EOp})} &&
{}_{\COp}\E\ar[dl]^{\Sigma^{-n}B^{n}}\ar[ll]_{\epsilon^*} \\
& \E & }
\end{equation*}
is still commutative at $n = \infty$.

\medskip
We prove in proposition~\ref{QuasiFreeLifting:SuspensionCofibration}
that the suspension morphisms
are all cofibrations in the category of right $\ROp$-modules.
As a corollary:

\begin{prop}\label{EinfinityHomology:CofibrantColimit}
The infinite bar modules $\Sigma^{-\infty}B^{\infty}_{\ROp}$
form cofibrant objects in the category of right $\ROp$-modules
for $\ROp = \EOp,\COp$.

The morphisms $\epsilon: \Sigma^{-n}B^{n}_{\EOp}\xrightarrow{\sim}\Sigma^{-n}B^{n}_{\COp}$, $n\in\NN$,
yield a weak-equivalence at $n = \infty$.\qed
\end{prop}

\medskip
We adapt the arguments of~\S\ref{EnHomology}
to prove that the infinite bar complex determines the homology of $\EOp$-algebras.
We prove that the infinite bar modules $\Sigma^{-\infty}B^{\infty}_{\ROp}$, where $\ROp = \EOp,\COp$,
define a cofibrant replacement of~$I$
in the category of right $\ROp$-modules
to obtain the relation $\Tor^{\ROp}(I,A) = H_*(\Sym_{\ROp}(\Sigma^{-\infty}B^{\infty}_{\ROp},A)$.
We can address the cases $\ROp = \COp$ and $\ROp = \EOp$
in parallel.
In fact,
proposition~\ref{EinfinityHomology:CofibrantColimit} implies that $\Sigma^{-\infty}B^{\infty}_{\EOp}$
forms a cofibrant replacement of $\Sigma^{-\infty}B^{\infty}_{\COp}$
in the category of right $\EOp$-modules.
Therefore it is sufficient to gain the result at the level of the commutative operad $\COp$.
This situation contrasts with the case of finite iterations of the bar construction,
addressed in~\S\ref{EnHomology},
where we can not avoid the study of extended bar modules $B^n_{\EOp_n}$.

\medskip
The right $\ROp$-module $\Sigma^{-\infty}B^{\infty}_{\ROp}$
comes equipped with a natural augmentation $\epsilon: \Sigma^{-\infty}B^{\infty}_{\ROp}\rightarrow I$
yielded by the morphisms of~\S\ref{EnHomology:IteratedBarModuleApplication}:
\begin{equation*}
\Sigma^{-n} B^n_{\ROp} = \Sigma^{-n} B^n(\ROp)\rightarrow\Sigma^{-n} B^n(I)\rightarrow I.
\end{equation*}
We already observed that $\Sigma^{-\infty}B^{\infty}_{\ROp}$ forms a cofibrant object.
It remains to check:

\begin{lemm}\label{EinfinityHomology:ColimitHomology}
We have $H_*(\Sigma^{-\infty}B^{\infty}_{\ROp}) = I$.
\end{lemm}

\begin{proof}
We have $H_*(\Sigma^{-\infty}B^{\infty}_{\EOp}) = H_*(\Sigma^{-\infty}B^{\infty}_{\COp})$
by proposition~\ref{EinfinityHomology:CofibrantColimit}.
Therefore we are reduced to prove the vanishing of $H_*(\Sigma^{-\infty}B^{\infty}_{\ROp}) = \colim_n H_*(\Sigma^{-n}B^{n}_{\ROp})$
for the commutative operad $\ROp = \COp$.

By proposition~\ref{UsualBarComplexes:FreeCommutativeAlgebra},
we have weak-equivalences
\begin{equation*}
\COp(\Sigma^n I)\xrightarrow{\sim} B(\COp(\Sigma^{n-1} I))\xrightarrow{\sim}\cdots
\xrightarrow{\sim} B^n(\COp(I)) = B^n_{\COp}.
\end{equation*}

In~\S\ref{IteratedBarModules:Indecomposables},
we observe that the morphism
\begin{equation*}
\sigma_*: \Sigma A\rightarrow H_*(B(A))
\end{equation*}
induced by the suspension factors through the indecomposables $\Indec A$,
for every commutative algebra $A$.
In the case $A = \COp(\Sigma^{n-1} I)$,
we have $\Indec A = \Sigma^{n-1} I$ and we obtain a commutative diagram
\begin{equation*}
\xymatrix{ & \Sigma\COp(\Sigma^{n-1} I)\ar[dd]_{\sigma}\ar[r]^{\sim}\ar[dl] & \Sigma B^{n-1}_{\COp}\ar[dd]^{\sigma} \\
\Sigma^n I\ar@{.>}[dr] && \\
& B(\COp(\Sigma^{n} I))\ar[r]_{\sim} & B^{n}_{\COp} },
\end{equation*}
from which we deduce that the morphism
\begin{equation*}
H_*(\Sigma^{1-n} B^{n-1}_{\COp})\xrightarrow{\sigma_*} H_*(\Sigma^{-n} B^{n}_{\COp})
\end{equation*}
admits a factorization
\begin{equation*}
\xymatrix{ H_*(\Sigma^{1-n} B^{n-1}_{\COp})\ar[rr]^{\sigma_*}\ar@{.>}[dr] && H_*(\Sigma^{-n} B^{n}_{\COp}) \\
& I\ar@{.>}[ur] & }.
\end{equation*}
Hence we conclude
\begin{equation*}
H_*(\Sigma^{-\infty}B^{\infty}_{\EOp}) = H_*(\Sigma^{-\infty}B^{\infty}_{\COp}) = \colim_n H_*(\Sigma^{-n} B^{n}_{\COp}) = I.
\end{equation*}
\end{proof}

To summarize:

\begin{prop}\label{EinfinityHomology:UnitCofibrantReplacement}
The morphism $\epsilon: \Sigma^{-\infty} B^{\infty}_{\ROp}\rightarrow I$
defines a weak-equivalence of right $\ROp$-modules,
so that $\Sigma^{-\infty} B^{\infty}_{\ROp}$
forms a cofibrant replacement of the composition unit $I$
in the category of right $\ROp$-modules for $\ROp = \EOp,\COp$.
\end{prop}

Hence,
we have
\begin{equation*}
H_*(\Sigma^{-\infty} B^{\infty}(A)) = H_*(\Sym_{\ROp}(\Sigma^{-\infty} B^{\infty}_{\ROp},A)) = \Tor^{\ROp}_*(I,A),
\end{equation*}
for $\ROp = \EOp,\COp$,
and from the identity $H^{\EOp}_*(A) = \Tor^{\EOp}_*(I,A)$
we conclude:

\begin{thm}\label{EinfinityHomology:Result}
For an $\EOp$-algebra~$A\in{}_{\EOp}\E$,
we have the identity
\begin{equation*}
H_*(\Sigma^{-\infty} B^{\infty}(A)) = H^{\EOp}_*(A)
\end{equation*}
as long as $A$ forms a cofibrant object in the underlying category $\E$.\qed
\end{thm}

Recall that the identity $H^{\ROp}_*(A) = \Tor^{\ROp}_*(I,A)$
makes sense for a $\Sigma_*$-cofibrant operad only,
though the $\Tor$-functor $\Tor^{\ROp}_*(I,A)$
can be defined as long as the operad $\ROp$ is $\C$-cofibrant.
In the case of a commutative algebra,
we have a restriction relation
\begin{equation*}
\Sym_{\EOp}(\Sigma^{-\infty} B^{\infty}_{\EOp},\epsilon^* A)
\simeq\Sym_{\COp}(\Sigma^{-\infty} B^{\infty}_{\EOp}\circ_{\EOp}\COp,A)
\simeq\Sym_{\COp}(\Sigma^{-\infty} B^{\infty}_{\COp},A)
\end{equation*}
from which we deduce
\begin{equation*}
H^{\EOp}_*(\epsilon^* A)\simeq\Tor^{\EOp}_*(I,\epsilon^* A)\simeq\Tor^{\COp}_*(I,A)
\end{equation*}
(we prefer to mark the restriction functor $\epsilon^*: {}_{\COp}\E\rightarrow{}_{\EOp}\E$ in these formulas).

\subsubsection{Remarks: relationship with $\Gamma$-homology of commutative algebras}\label{EinfinityHomology:GammaHomology}
In the case of an $E_\infty$-operad $\EOp$,
the homology theory $H^{\EOp}_*(A)$,
defined abstractly in~\S\ref{EnHomology:OperadHomology},
represents the $\Gamma$-homology of $A$ over $\kk$
with trivial coefficients $\kk$.
The usual notation for this homology theory is $\HGamma_*(A|\kk,\kk) = H^{\EOp}_*(A)$.

The article~\cite{Robinson} gives another chain complex $\CXi_*(A|R,E)$
(rather denoted by $\Xi_*(A|R,E)$ in \emph{loc. cit.})
which determines the $\Gamma$-homology $\HGamma_*(A|R,E)$
in the case where $A$ is a commutative algebra over another commutative algebra $R$,
and for any coefficient $E$ in the category of $A$-modules.
The author of~\cite{Robinson}
deals with unital commutative algebras over $R$,
but a unital commutative algebra $A_+$
can be replaced by a quotient $A = A_+/R$
to give a normalized chain complex $\NXi_*(A|R,E)$ equivalent to $\CXi_*(A_+|R,E)$.

In the case $E = R = \kk$,
the normalized chain complex $\NXi_*(A|\kk,\kk)$
can be identified with the functor $\NXi_*(A|\kk,\kk) = \Sym_{\COp}(\NXi_{\COp},A)$
associated to a particular cofibrant replacement of $I$
in the category of right $\COp$-modules.
This cofibrant replacement $\NXi_{\COp}$ is defined over the ring $\kk = \ZZ$.

Observe that any category of dg-modules over a ring $R$
forms a symmetric monoidal category $\E = \C_{R}$
over the base category $\C = \C_{\kk}$
of dg-modules over~$\kk$.
An augmented algebra over $R$
is equivalent to a (non-unital non-augmented) commutative algebra in $\C_{R}$
and the extended functor $\Sym_{\COp}(\NXi_{\COp}): \C_{R}\rightarrow\C_{R}$
satisfies
\begin{equation*}
\NXi_*(A|R,R) = \Sym_{\COp}(\NXi_{\COp},A),
\end{equation*}
for any $A\in{}_{\COp}\C_{R}$.

The chain complex of~\cite{Robinson} inherits a grading from $\NXi_{\COp}$
and an internal grading from $A$,
so that $\NXi_*(A|R,R)$
forms naturally a bigraded object.
So does the infinite bar complex
of commutative algebras.
The morphism
\begin{equation*}
\underbrace{\Sym_{\COp}(\Sigma^{-\infty} B^{\infty}_{\COp},A)}_{= \Sigma^{-\infty} B^{\infty}(A)}
\xrightarrow{\Sym_{\COp}(\psi,A)}
\underbrace{\Sym_{\COp}(\NXi_{\COp},A)}_{\NXi_*(A|\kk,\kk)}
\end{equation*}
associated to a lifting
\begin{equation*}
\xymatrix{ & \NXi_{\COp}\ar[d]^{\sim} \\
\Sigma^{-\infty} B^{\infty}_{\COp}\ar[r]_{\sim}\ar@{.>}[ur]^{\psi} & I }.
\end{equation*}
preserves clearly bigradings.
Hence, in the context of commutative algebras over a ring $R$,
our theorem gives a natural isomorphism of bigraded objects
\begin{equation*}
H_{*,*}(\Sigma^{-\infty} B^{\infty}(A)) = \HGamma_{*,*}(A|R,R)
\end{equation*}
for every commutative algebra $A$ which is cofibrant (or simply flat) in the underlying category of dg-modules over $R$.

\subsubsection{Remarks: relationship with Koszul duality and applications to the Lie operad}\label{EinfinityHomology:LieOperad}
The Koszul duality of operads gives another quasi-free complex $K(I,\COp,\COp)$
together with an acyclic fibration $\epsilon: K(I,\COp,\COp)\xrightarrow{\sim} I$.
In the proof of proposition~\ref{UsualBarComplexes:HarrisonKoszul},
we already recalled that this Koszul complex
is identified with a desuspension of the Harrison complex of the commutative operad $\COp$,
viewed as a commutative algebra in the category of right modules over itself.
Equivalently,
the Koszul complex $K(I,\COp,\COp)$ is a quasi-free module such that
\begin{equation*}
K(I,\COp,\COp) = (\Sigma^{-1} L^c(\Sigma\COp),\partial) = (\Lambda^{-1}\LOp^{\vee}\circ\COp,\partial),
\end{equation*}
where $\LOp^{\vee}$ refers to the $\kk$-dual of the Lie operad $\LOp$.
This quasi-free module is not cofibrant
since the Lie operad does not form a cofibrant $\Sigma_*$-module.

Nevertheless, we can pick a lifting in the diagram
\begin{equation*}
\xymatrix{ & K(I,\COp,\COp)\ar@{->>}[d]^{\sim} \\
\Sigma^{-\infty} B^{\infty}_{\COp}\ar@{.>}[ur]\ar[r]_{\sim} & I }
\end{equation*}
to obtain a weak-equivalence $\kappa: \Sigma^{-\infty} B^{\infty}_{\COp}\xrightarrow{\sim} K(I,\COp,\COp)$.
By~\cite[Theorem 2.1.14]{FressePartitions} or~\cite[Theorem 15.3.A]{Bar0},
the quasi-free structure is sufficient to ensure that the morphism $\kappa\circ_{\COp} I$
is still a weak-equivalence.
Thus,
since $\Sigma^{-\infty} B^{\infty}_{\COp}\circ_{\COp} I = \Sigma^{-\infty} B^{\infty}(I)$
and $K(I,\COp,\COp)\circ_{\COp} I = \Lambda^{-1}\LOp^{\vee}$,
we have a weak-equivalence of $\Sigma_*$-modules
\begin{equation*}
\bar{\kappa}: \Sigma^{-\infty} B^{\infty}(I)\xrightarrow{\sim}\Lambda^{-1}\LOp^{\vee}.
\end{equation*}

Note that $\Sigma^{-\infty} B^{\infty}(I)$
forms a cofibrant $\Sigma_*$-module.
Thus the object $\Sigma^{-\infty} B^{\infty}(I)$
defines a cofibrant replacement of the $\kk$-dual of the Lie operad $\Lambda^{-1}\LOp^{\vee}$
in the category of $\Sigma_*$-modules.
From this observation,
we deduce an identity
\begin{equation*}
H_*(\Sym(\Sigma^{-\infty} B^{\infty}(I),\kk)) = \bigoplus_{r=0}^{\infty}\Tor^{\Sigma_r}_*(\Lambda^{-1}\LOp^{\vee}(r),\kk),
\end{equation*}
where $\Sym(\Sigma^{-\infty} B^{\infty}(I),\kk)$
denotes the image of the free module of rank $1$
under the symmetric tensor functor associated to $\Sigma^{-\infty} B^{\infty}(I)$.
The object $\Sym(\Sigma^{-\infty} B^{\infty}(I),\kk)$
is identified with the iterated bar complex $\Sigma^{-\infty} B^{\infty}(A)$
of a trivial algebra $A = \kk e$
(which represents the non-unital algebra associated to a unital exterior algebra in one generator),
since we have the restriction relation:
\begin{equation*}
\Sym(\Sigma^{-\infty} B^{\infty}(I),\kk)
= \Sym(\Sigma^{-\infty} B^{\infty}_{\COp}\circ_{\COp} I,\kk)
= \Sym_{\COp}(\Sigma^{-\infty} B^{\infty}_{\COp},\kk e) = \Sigma^{-\infty} B^{\infty}(\kk e).
\end{equation*}

The arguments of~\cite{Cartan}
give the homology of the infinite bar complex $\Sigma^{-\infty} B^{\infty}(A)$
for this particular commutative algebra,
for every ground ring $\kk$.
Thus our result relates the calculation of $\Tor^{\Sigma_r}_*(\Lambda^{-1}\LOp^{\vee}(r),\kk)$
to classical homological computations.
Such $\Tor$-functors are determined by other methods in~\cite{AroneKankaanrinta,AroneMahowald}
in the case $\kk = \FF_p$
(see also~\cite{BetleySlominska} for another approach to this computation).

\part*{Afterword: applications to the cohomology of iterated loop spaces}\label{Conclusion}

The goal of this concluding part is to explain the applications of our main theorems
to the cohomology of iterated loop spaces.

Let $\bar{N}^*(X)$ denote to the reduced normalized cochain complex
of a simplicial set $X$.
By~\cite{HinichSchechtman} (see also~\cite{BergerFresse,McClureSmith}),
the cochain complex $\bar{N}^*(X)$ inherits an action of an $E_\infty$-operad $\EOp$
so that the map $N^*: X\mapsto\bar{N}^*(X)$
defines a functor from the category of simplicial sets $\Simp$
to the category of $\EOp$-algebras in dg-modules ${}_{\EOp}\C$.
Moreover,
all functorial actions of an $E_\infty$-operad $\EOp$
on the functor $\bar{N}^*: \Simp\rightarrow\C$
are homotopy equivalent.
In~\cite{BergerFresse},
we prove by an explicit construction that the Barratt-Eccles operad
is an instance of an $E_\infty$-operad
which acts on cochain complexes.

\medskip
In~\cite{Bar1},
we prove that the $n$-fold bar complex $B^n(\bar{N}^*(X))$ of a cochain complex $\bar{N}^*(X)$
determines the cohomology of the $n$-fold loop space of $X$
under reasonable finiteness and completeness assumptions.
Thus,
with the new results of the present article,
we obtain (in the case $\kk = \FF_p$):

\begin{mainthm}\label{Conclusion:TopologicalApplications}
Let $X$ be a pointed simplicial set whose cohomology modules $H^*(X,\FF_p)$ are degreewise finitely generated.
Let $\bar{N}^*(X)$ be the reduced cochain complex of $X$
with coefficients in $\kk = \FF_p$.

Let $\EOp$ be any $E_\infty$-operad which, like the Barratt-Eccles operad,
acts on cochain complexes of spaces
and fulfil the requirements of~\S\ref{EnDefinition:KCellOperads} so that the conclusions of theorem~\ref{EnDefinition:Result}
and theorem~\ref{EnHomology:Conclusion} hold for this operad.
Then we have identities
\begin{equation*}
H^{\EOp_n}_*(\bar{N}^*(X)) = H_*(\Sigma^{-n} B^n(\bar{N}^*(X))) = \colim_s\bar{H}^*(\Sigma^n\Omega^n R_s X),
\end{equation*}
where $R_s X$ refers to Bousfield-Kan' tower of $X$.\qed
\end{mainthm}

A similar result can be stated in the case $\kk = \ZZ$ or in the case $\kk = \QQ$ (assuming in this case that $X$ is $n$-connected).

\medskip
The explicit construction of the $n$-fold bar complex
implies the existence of a spectral sequence
\begin{equation*}
B^n(H^*(X))\Rightarrow H^*(B^n(X)).
\end{equation*}
We refer to~\cite{SmirnovChain,SmirnovHomology}
for another definition of a similar spectral sequence
converging to $H_*(\Omega^n X)$.
We conjecture that our spectral sequence is isomorphic (from the $E^2$-stage and up to $\kk$-duality)
to the $H_*(-,\kk)$-version of the spectral sequence of~\cite{AhearnKuhn},
defined from Goodwillie's approximations of the functor $\Sigma^{\infty}\Map(S^n,X)_+$.
We prove in~\cite{Bar3} that $E_n$-operads are (up to operadic suspension) self-dual
in the sense of Koszul duality of operads.
We deduce from this result another representation of the homology theory $H^{\EOp_n}_*(A)$
which relates the $n$-fold bar complex $B^n(A)$
to the $\kk$-dual of the $E_n$-operad $\EOp_n$.
We conjecture that this relationship reflects the occurrence of the little $n$-cubes operad
in~\cite{AhearnKuhn}.

The reference~\cite{Cartan} gives the homology of the iterated bar complexes
of many usual commutative algebras,
like exterior algebras, polynomial algebras, divided power algebras, abelian group algebras.
These results could be used to determine $E^2$-terms in the spectral sequence $B^n(H^*(X))\Rightarrow H^*(B^n(X))$.
Note that the calculations of~\cite{Cartan}
are performed over $\ZZ$.

\end{document}